\documentclass[3p,times]{elsarticle}
\usepackage{amssymb}
\usepackage{amsmath}
\usepackage{mathtools}
\usepackage{mathabx}
\usepackage{xcolor}
\usepackage{lineno}
\usepackage{subcaption}
\usepackage{adjustbox}
\graphicspath{{figures/}} 
\usepackage{graphicx}
\usepackage{float}
\newdefinition{rem}{Remark}
\usepackage{amsthm}
\usepackage{cancel}
\usepackage{csquotes}

\newif\ifrevisions
\revisionsfalse      

\ifrevisions
\else
    \colorlet{red}{black}
    \colorlet{blue}{black}
\fi

\usepackage{placeins}
\usepackage{etoolbox}

\pretocmd{\subsection}{\FloatBarrier}{}{}

\journal{CMAME}

\begin{document}

\begin{frontmatter}

\title{Local $h$-, $p$-, and $k$-Refinement Strategies for the Isogeometric Shifted Boundary Method Using THB-Splines}

\author[TUM,OTH]{Christoph Hollweck}
\author[UPC,CIMNE]{Andrea Gorgi}
\author[UPC,CIMNE]{Nicolò Antonelli}
\author[OTH]{Marcus Wagner}
\author[TUM]{Roland Wüchner}


\affiliation[TUM]{organization={Technische Universität München, Lehrstuhl für Statik und Dynamik},
            addressline={Arcisstr. 21}, 
            city={80333 München},
            country={Germany}}

\affiliation[OTH]{organization={OTH Regensburg, Labor Finite-Elemente-Methode},
            addressline={Galgenbergstraße 30}, 
            city={93053 Regensburg},
            country={Germany}}

\address[UPC]{Universidad Politécnica de Cataluña (UPC), Campus Norte UPC, 08034 Barcelona, Spain}
\address[CIMNE]{Centre Internacional de Mètodes Numèrics en Enginyeria (CIMNE), Campus Norte UPC, 08034 Barcelona, Spain}


\begin{abstract} 
The concept of trimming, embedding, or immersing geometries into a computational background mesh has gained considerable attention in recent years, particularly in isogeometric analysis (IGA). In this approach, the physical domain is represented independently from the computational mesh, allowing the latter to be generated more easily compared with body-fitted meshes. While this facilitates the treatment of complex geometries, it also introduces challenges, such as ill-conditioning of the stiffness matrix caused by small cut elements and difficulties in accurately enforcing boundary conditions.
A recently proposed technique to address these issues is the Shifted Boundary Method (SBM), which represents the computational domain solely through uncut elements and enforces boundary conditions via a Taylor expansion from a surrogate boundary to the true boundary. Previous studies have shown that, for Neumann boundary conditions, the flux evaluation requires additional derivatives in the Taylor expansion, effectively reducing the order of convergence by one.
In this work, we investigate for the first time the performance of SBM combined with Truncated Hierarchical B-splines (THB-splines) under various local refinement strategies. In particular, we {\color{red}study} local $p$- and $k$-refinement for THB-splines and compare them with local $h$-refinement and the unmodified SBM. Furthermore, we propose an enhanced shift operator based on a tensor-product index set, which includes derivatives of total order up to $2p$, whereas the classical directional Taylor operator is limited to total order $p$. The study assesses accuracy, stability, and computational efficiency for benchmark problems on trimmed domains. The results highlight how different refinement strategies affect convergence behavior in trimmed IGA formulations using SBM and demonstrate that targeted degree elevation can mitigate the Neumann boundary limitations of the standard method.

\end{abstract}

\begin{highlights}
    \item Local $h$-, $p$- and $k$-refinement strategies based on THB-splines are investigated in the context of the Shifted Boundary Method (SBM).
    \item A framework for combining the SBM with locally refined THB-splines is presented.
    \item An enhanced shift operator based on a tensor-product index set is introduced.
    \item The approach is assessed on benchmark problems with varying geometric complexity.
    \item Local refinement improves convergence for Neumann boundary conditions and extends the regime in which convergence remains close to optimal, with $p$- and $k$-refinement providing the strongest improvements.
\end{highlights}


\end{frontmatter}

\section{Introduction} \label{sec:introduction}

Computer-Aided Design (CAD) and Computer-Aided Engineering (CAE) have traditionally developed with distinct goals and methods, resulting in fundamentally different geometric representations. CAD systems commonly use Non-Uniform Rational B-splines (NURBS) combined with trimming operations to accurately model complex shapes~\cite{piegl_nurbs_1997}. In contrast, numerical simulation techniques such as the finite element method (FEM) rely on body-fitted meshes constructed from simpler polynomial basis functions optimized for analysis rather than design~\cite{Bathe}. This disconnect necessitates costly and error-prone geometry conversions when transitioning from CAD to analysis models.

Isogeometric Analysis (IGA) was introduced to bridge the long-standing gap between CAD and numerical analysis by employing the same smooth spline basis functions for both geometry representation and field approximation~\cite{hughes2005isogeometric, cottrell2009isogeometric}.
This unified framework enables an accurate description of CAD geometries and provides globally smooth basis functions with higher continuity across element interfaces, which can significantly improve the accuracy per degree of freedom compared to conventional finite element methods~\cite{cottrell2006isogeometric}.


Despite these strengths, the practical use of IGA is often limited by geometric complexity. The generation of analysis-suitable parameterizations from intricate CAD models remains a major challenge. Body-fitted approaches, which require reparametrization and remeshing, are computationally demanding and inherit many of the meshing difficulties familiar from traditional finite element analysis~\cite{rorris2019ansa, rorris2020latest}.

To avoid these issues, embedded or immersed methods have been developed. These techniques embed the physical geometry into a simpler, structured background mesh and remove the void regions outside the physical domain via trimming procedures. This approach allows the exact geometry to be represented without the need for body-fitted meshing. Two popular examples of this approach include the Isogeometric B-Rep Analysis (IBRA)~\cite{breitenberger_analysis_2015, teschemacher_realization_2018, teschemacher2022concepts} and the Finite Cell Method (FCM)~\cite{parvizian2007finite, rank2012geometric}. 

{\color{red}
The use of higher-order and higher-continuity spline bases in numerical methods is often associated with IGA. However, the idea of employing splines instead of Lagrange polynomials for the discretization of variational problems is neither new nor exclusive to IGA. Early investigations on spline-based discretizations can be found prior to the emergence of IGA, see, e.g., \cite{Hoellig2003, Kwok2001}.
In many immersed and embedding approaches, the geometric description serves primarily to define the physical boundary, for instance for numerical integration on cut elements and for the imposition of boundary conditions, and does not require a CAD-based description. In this context, the term IGA is often used in a broader sense to denote spline-based approximation spaces, rather than a strict CAD-to-analysis workflow.}

From a computational standpoint, trimmed elements require specialized numerical integration schemes that account for irregular cut shapes~\cite{breitenberger_analysis_2015, parvizian2007finite, rank2012geometric, Scholz2019ErrorCorrection}. In arbitrary geometries, small cut cells naturally arise, causing numerical challenges such as ill-conditioning of system matrices, degraded convergence rates, and, in explicit dynamics, the emergence of unstable nodal motions often referred to as "flying nodes" or "light control points"~\cite{leidinger_explicit_nodate}.

To address these numerical difficulties, several stabilization techniques have been developed, including preconditioning and ghost penalty stabilization~\cite{dePrenter2023}. Recently, the Shifted Boundary Method (SBM) was introduced in FEM to overcome the small cut cell problem~\cite{main2018shifted, atallah2021shifted, li2025shifted}. The SBM tackles these difficulties by replacing the physical boundary with a surrogate boundary defined by the background mesh edges. Boundary conditions are imposed on this surrogate boundary by shifting them from the true boundary via a Taylor expansion, thus eliminating the need for complicated integration over cut elements.

Although this method reliably enforces Dirichlet conditions with optimal convergence both in FEM and IGA, it struggles with Neumann boundary conditions~\cite{antonelli2024shifted, antonelli2025isogeometric}.
The truncated Taylor expansion used in SBM limits the accuracy of Neumann flux approximations due to incomplete gradient recovery~\cite{antonelli2024shifted}. 
A  recent study achieves optimal convergence in linear finite elements via accurate evaluation of the variational integrals in the gap between the surrogate and the true boundaries ~\cite{collins2025gap}.
Mixed formulations introducing additional variables can improve Neumann treatment, but they increase system size and complexity. 
To address this challenge, local refinement strategies inherent to IGA can be exploited.

IGA provides various refinement options, including $h$-refinement (knot insertion), $p$-refinement (degree elevation), and $k$-refinement (a combination of degree elevation and knot insertion)~\cite{hughes_isogeometric_2009}. Local $h$-refinement using hierarchical splines, such as Truncated Hierarchical B-splines (THB-splines), has gained popularity due to its simple construction and the preservation of essential spline properties, for example linear independence and the Partition of Unity~\cite{giannelli_thB-splines_2012, giannelli_strongly_2014}.
 In contrast, local $p$- and $k$-refinement remain less explored in practical scenarios despite their theoretical potential~\cite{schillinger_2012_diss}.

This work proposes a novel approach to improve the enforcement of Neumann boundary conditions in the SBM by incorporating local $h$-, $p$-, and $k$-refinement within an isogeometric framework using THB-splines. The key insight is that the order of the Taylor expansion of boundary fluxes is limited by the polynomial degree $p$ of the basis functions. While the solution approximation can reach order $p+1$, higher derivatives required for accurate gradient representation vanish in standard bases. By locally elevating the spline degree in regions with Neumann conditions, higher-order Taylor expansions of gradients become feasible, aiming to restore optimal convergence behavior. THB-splines enable this selective degree elevation without increasing the polynomial degree globally. 

Furthermore, we propose an enhanced shift operator based on a tensor-product multi-index expansion. 
While the classical directional Taylor operator retains all derivatives of total order $p$, 
the enhanced formulation includes derivatives with a total order up to $2p$.
{\color{red}The effectiveness of the proposed method is systematically assessed through numerical experiments.} The results demonstrate that combining SBM with localized degree elevation improves Neumann boundary treatment in immersed isogeometric discretizations. This strategy opens new possibilities for accurate and efficient simulation of problems that involve complex geometries and boundary conditions.
\newpage

\section{Preliminaries} \label{sec:preliminaries}
\newcommand{\uTilde}{\widetilde{\mathbf{u}}}
\newcommand{\OmegaTilde}{\widetilde{\Omega}}
\newcommand{\GammaTilde}{\widetilde{\Gamma}}
\newcommand{\ntilde}{\widetilde{\mathbf{n}}}
\newcommand{\Bu}{\mathbf{u}}
\newcommand{\Bv}{\mathbf{v}}
\newcommand{\Bd}{\mathbf{d}}
\newcommand{\xTilde}{\widetilde{\mathbf{x}}}
\newcommand{\Bx}{\mathbf{x}}
\newcommand{\Bn}{\mathbf{n}}

{\color{red}
Isogeometric analysis (IGA), originally introduced in~\cite{hughes2005isogeometric, cottrell2009isogeometric}, is commonly associated with the use of spline bases for the numerical solution of partial differential equations. While its original motivation was to bridge CAD and finite element analysis, the use of splines as approximation spaces is not restricted to a CAD-to-analysis workflow. Unlike classical finite element methods based on Lagrangian polynomials, IGA employs spline-based basis functions that enable higher-order approximations with increased smoothness. In the present work, B-splines are used exclusively as approximation functions for the solution field. No geometric parametrization is introduced. Instead, the physical domain is identified with the parametric domain, such that the physical coordinate $\mathbf{x}$ coincides with the parametric coordinate. This identification allows the spline basis functions, originally defined on the parametric space, to be evaluated directly in physical coordinates without introducing a mapping.

Let $\Xi = \{\xi_1, \xi_2, \ldots, \xi_{n+p+1}\}$ denote a non-decreasing knot vector, where $p$ is the polynomial degree and $n$ denotes the number of basis functions associated with the knot vector. The knots partition the parametric domain into knot spans $[\xi_i, \xi_{i+1}]$ with $\xi_i \neq \xi_{i+1}$. The univariate B-spline basis functions $N_{i,p}(\xi)$ are defined recursively via the Cox--de Boor formula,

\begin{equation}
N_{i,0}(\xi) =
\begin{cases}
1 & \xi_i \leq \xi < \xi_{i+1}, \\
0 & \text{otherwise},
\end{cases}
\end{equation}

and for $p>0$,
\begin{equation}
N_{i,p}(\xi) =
\frac{\xi - \xi_i}{\xi_{i+p} - \xi_i} N_{i,p-1}(\xi)
+
\frac{\xi_{i+p+1} - \xi}{\xi_{i+p+1} - \xi_{i+1}} N_{i+1,p-1}(\xi),
\end{equation}
where terms with zero denominators are treated as zero. A knot of multiplicity $m$ reduces the continuity to $C^{p-m}$ at that location, allowing the regularity of the approximation space to be controlled locally. Open knot vectors, in which the first and last knots are repeated $p+1$ times, ensure interpolation at the domain boundaries.

For problems defined in $d$ spatial dimensions, multivariate B-spline basis functions are constructed as tensor products of the corresponding univariate bases defined in each parametric direction. The resulting multivariate basis functions are denoted by $N_A$, where the multi-index $A$ collects the univariate indices. For notational simplicity, the dependence on the individual parametric coordinates is suppressed in the following.

The discrete solution reads
\begin{equation}
u_h(\mathbf{x}) = \sum_A N_A(\mathbf{x}) \, u_A ,
\end{equation}
where $u_A$ denote the degrees of freedom. Since B-spline bases are generally non-interpolatory, these coefficients cannot be interpreted as nodal values. The associated control points should therefore be understood as carriers of the discrete unknowns rather than as geometric entities.

B-spline based discretizations offer several advantageous properties: the basis functions are non-negative, possess local convex support, and form a partition of unity. Most notably, B-splines enable up to $C^{p-1}$ continuity across element interfaces, in contrast to the $C^0$ continuity of standard Lagrange bases. This increased smoothness leads to superior convergence rates and higher accuracy in the approximation of smooth solutions~\cite{cottrell2007refinement}.

The practical deployment of IGA is often limited by the complexity of constructing boundary-fitted mappings. Immersed and unfitted approaches address this limitation by decoupling the discretization space from the geometric representation and employing an unfitted background mesh to define the active computational domain. Such approaches provide the methodological foundation for the SBM considered in this work.}

\clearpage

\section{Methodology} \label{sec:methodology}

This section outlines the theoretical and numerical foundations of our work. We first introduce the strong and weak forms of the Poisson problem, including Nitsche’s method for the weak enforcement of Dirichlet boundary conditions. We then describe the SBM, which enables robust handling of complex geometries on unfitted meshes. {\color{red} Finally, THB-splines are introduced as the underlying discretization technology, and local $h$-, $p$-, and $k$-refinement strategies are considered within this framework to provide flexibility in the numerical analysis.}

\subsection{Strong form}
Let $\Omega \subset \mathbb{R}^2$ be a bounded, connected Lipschitz domain
with boundary $\Gamma$. The boundary is decomposed as
\[
\Gamma = \Gamma_D \cup \Gamma_N,
\qquad
\Gamma_D \cap \Gamma_N = \emptyset.
\]
The classical Poisson problem reads
\begin{equation}
-\,\Delta u = f \quad \text{in } \Omega,
\end{equation}
with boundary conditions
\begin{equation}
u = u_D \quad \text{on } \Gamma_D,
\qquad
\partial_n u = t_N \quad \text{on } \Gamma_N,
\end{equation}
where $f \in L^2(\Omega)$ and
\(
\partial_n u := \nabla u \cdot \mathbf{n},
\)
with $\mathbf{n}$ denoting the outward unit normal on $\Gamma$.
{\color{red}The corresponding weak formulation is incorporated directly into the Nitsche-based framework described in the following section.
}


\begin{rem}[Strong vs.\ weak Dirichlet enforcement]
For conforming $C^{0}$ finite elements, the condition $u = u_D$ on $\Gamma_D$ is typically imposed \emph{strongly} by restricting the trial and test function spaces. In IGA, however, the non-interpolatory nature of B-splines and NURBS makes strong enforcement impractical in most applications. Therefore, a \emph{weak} enforcement such as Nitsche’s method is commonly used.
\end{rem}

\subsection{Nitsche’s weak enforcement of Dirichlet conditions}

Within the Nitsche framework, we seek $u \in H^1(\Omega)$ such that
\begin{equation}
a(u,v) = \ell(v) \qquad \forall\, v \in H^{1}(\Omega),
\label{eq:nitsche-form}
\end{equation}
with bilinear and linear forms
\begin{equation}
\begin{aligned}
a(u,v) &:= (\nabla u,\nabla v)_{\Omega}
          - \langle \partial_n u , v \rangle_{\Gamma_D}
          - \theta\,\langle \partial_n v , u \rangle_{\Gamma_D}
          + \Bigl\langle \frac{\alpha}{h}\,u , v \Bigr\rangle_{\Gamma_D},
\\[4pt]
\ell(v) &:= (f,v)_{\Omega}
          + \langle t_N , v \rangle_{\Gamma_N}
          - \theta\,\langle \partial_n v , u_D \rangle_{\Gamma_D}
          + \Bigl\langle \frac{\alpha}{h}\,u_D , v \Bigr\rangle_{\Gamma_D}.
\end{aligned}
\label{eq:nitsche-forms}
\end{equation}
Here, $\alpha > 0$ is chosen sufficiently large to ensure coercivity of the symmetric formulation,
$h$ denotes a local length scale, and
$\theta \in \{-1,1\}$ determines the formulation type:
$\theta = -1$ corresponds to the skew-symmetric variant
(penalty-free for the Poisson problem with $\alpha = 0$~\cite{burman2012penalty,hu2018skew}),
while $\theta = 1$ corresponds to the symmetric variant.
Nitsche’s method provides a consistent framework for the weak imposition of Dirichlet conditions
on both boundary-fitted and immersed meshes and will be employed throughout this work.

\subsection{THB-splines}
THB-splines are built on a hierarchical structure of nested spline spaces and meshes, enabling local refinement. The refined spaces are generated through standard global B-spline refinement strategies. We consider a nested sequence of multivariate spline spaces
$$
V^0 \subset V^1 \subset \dots \subset V^L,
$$
where each level $V^l$ is defined by a spline degree $p^l$ and knot vector $\Xi^l$.
Traditionally, level $l+1$ is obtained from level $l$ via uniform $h$-refinement, i.e., by bisecting all non-zero knot spans. However, it is equally possible to apply $p$-refinement (increasing the polynomial degree without modifying the mesh) or to combine $p$- and $h$-refinement sequentially, here referred to as $k$-refinement. This generalization offers greater flexibility in locally controlling the approximation properties. Each level $l$ is associated with a mesh $\mathcal{M}^l$ and a corresponding set of basis functions $\mathcal{N}^l$, defined by the knot vectors and degrees at that level. Starting from the coarsest level $l=0$, local refinement regions $\Omega^{l+1}$ are introduced, subject to the constraint $\Omega^{l+1} \subseteq \Omega^l$. These regions are chosen such that at least one basis function of level $l$ is fully supported within $\Omega^{l+1}$. By marking basis functions (rather than elements) for refinement, the refinement domain $\Omega^{l+1}$ is implicitly defined as the union of the supports of all marked functions.

The hierarchical mesh \( \hat{\mathcal{M}}\) is defined as the set of active elements \( e \) given by

\begin{equation}
\hat{\mathcal{M}} := \left\{ e \in \mathcal{M}^l : e \subseteq \Omega^l \land e \not\subseteq \Omega^{l+1}, \; l = 0, \dots, L - 1 \right\}.
\label{eq:meshthb}
\end{equation}
The THB basis \( \mathcal{T} \) is defined recursively as follows:
\begin{equation}
\begin{cases}
\mathcal{T}^0 := \mathcal{N}^0, \\

\mathcal{T}^{l+1} := \left\{ \text{trunc}(T_i^l): T_i^l \in \mathcal{T}^l \land \text{supp}(T_i^l) \nsubseteq \Omega^{l+1} \right\} \cup \\

\quad \quad \quad \, \, \, \, \left\{ \text{child}(T_i^l)_j \in \mathcal{N}^{l+1} : \text{supp}(T_i^l) \subseteq \Omega^{l+1} \right\} \quad l = 0, \dots, L - 1, \\
\mathcal{T} := \mathcal{T}^L
\end{cases}
\label{eq:recursive}
\end{equation}
where the truncation operator is given by:
\begin{equation}
\text{trunc}(T_i^l) = \sum_{\substack{j : \text{child}(T_i^l)_j \in \mathcal{N}^{l+1} \\ \text{supp}(\text{child}(T_i^l)_j) \nsubseteq \Omega^{l+1}}} \lambda_{ij} \cdot \text{child}(T_i^l)_j.
\label{eq:trunc}
\end{equation}
where the coefficients $\lambda_{ij}$ arise from the two-scale relation. The two-scale relation states that each basis function $N_i^l$ can be expressed as a linear combination of its children from $\mathcal{N}^{l+1}$. Children are the functions of $\mathcal{N}^{l+1}$ with support contained in $\text{supp}(N_i^l)$. {\color{red} According to the truncation procedure, a parent function $T_i^l$
is re-expressed via the two-scale relation, while coefficients
$\lambda_{ij}$ corresponding to already active child functions
$T_j^{l+1}$ are set to zero. 

Originally for the classical THB-spline formulation, the activation and deactivation of basis functions were primarily governed by their support \cite{giannelli_thB-splines_2012, giannelli_strongly_2014, Buffa2022_AdaptiveIGA}. The procedure adopted here, in which basis functions are activated and deactivated based on the mother–child relation, does not follow the classical THB-spline definition. Instead, it is based on the construction proposed in \cite{BuffaGarau2016RefinableSpaces}, which itself builds upon earlier concepts introduced in \cite{Grinspun2002CHARMS} and \cite{Krysl2003NaturalHierarchical}. 
Since it is known \emph{a priori} which mother splines contain which child splines, both the data structure and the refinement procedure become particularly simple. From a mathematical perspective, the resulting space retains the same approximation properties as the standard formulation while exhibiting a lower dimension. Implementation details can be found in \cite{GARAU201858}.}

{\color{red} Early work on THB-splines primarily focused on local $h$-refinement based on a hierarchy of globally $h$-refined tensor-product spaces. However, different refinement approaches can be employed; see, e.g., \cite{cottrell2009isogeometric} for a detailed discussion in the context of IGA.

\begin{itemize}
    \item \emph{$h$-refinement} (knot insertion) increases the local resolution by inserting additional knots while keeping the polynomial degree unchanged. In contrast to classical FEM, where $h$-refinement is typically associated with $C^0$ continuity, knot insertion in spline spaces preserves the continuity determined by the knot multiplicities. A $C^0$ analogue of classical FEM $h$-refinement is recovered only if newly inserted knots are repeated $p$ times.

    \item \emph{$p$-refinement} (degree elevation) increases the polynomial degree $p$ while preserving the existing continuity and remain the element sizes unchanged. When the initial discretization is $C^0$, degree elevation coincides with classical FEM $p$-refinement.

    \item \emph{$k$-refinement}, as introduced in the IGA literature~\cite{cottrell2009isogeometric}, consists of performing degree elevation prior to knot insertion. This results in newly inserted knots attaining the highest continuity permitted by the elevated polynomial degree. In the classical global setting, $k$-refinement is often associated with obtaining maximum continuity across all element interfaces. In practice, however, once a mesh contains multiple elements, the continuity at existing knots is fixed by their multiplicities and cannot be increased by subsequent refinement. Consequently, maximum continuity can only be achieved at newly inserted knots.
\end{itemize}

In the present work, $k$-refinement is applied locally within the THB-spline hierarchy. More precisely, a single local $k$-refinement step is defined as the sequence of a local $p$-refinement (degree elevation) followed by a local $h$-refinement (knot insertion) step within the marked region. Consequently, newly inserted knots attain the highest continuity 
permitted by the elevated polynomial degree, while existing knot lines retain their original continuity. Owing to the nested-space structure of THB-splines, all refinement operations generate nested approximation spaces.} On the level of basis functions, this implies that any basis function before refinement can be expressed as a linear combination of basis functions after refinement. The basis functions involved in this relation are called children of the corresponding coarse (parent or mother) function. This two-scale relation is expressed as
\begin{equation}
\label{eq:twoscale_sum}
N_i^l = \sum_{j=1}^n \lambda_{ij} \, N_j^{l+1},
\end{equation}
where a necessary condition for a nonzero coefficient \(\lambda_{ij}\) is that the child basis function \(N_j^{l+1}\) has support contained within the support of the mother basis function \(N_i^l\) of level $l$. In general, a basis function from the set $\mathcal{N}^{l+1}$ can be a child of multiple basis functions from $\mathcal{N}^l$. By collecting the basis functions into vectors
\[
\mathbf{N}^l := \bigl[\,N_1^l,\dots,N_{n_l}^l\,\bigr]^{\top},
\qquad
\mathbf{N}^{l+1} := \bigl[\,N_1^{l+1},\dots,N_{n_{l+1}}^{l+1}\,\bigr]^{\top},
\]
and arranging the two-scale coefficients \(\lambda_{ij}\) in the refinement matrix
\[
\mathbf{\Lambda} := \bigl[\lambda_{ij}\bigr]\in\mathbb{R}^{\,n_l\times n_{l+1}},
\]
the two-scale relation can be written compactly as
\begin{equation}
\label{eq:twoscale_matrix}
\mathbf{N}^l = \mathbf{\Lambda}\,\mathbf{N}^{l+1}.
\end{equation}


The extension to the two-dimensional case is straightforward: each parametric direction is refined independently by computing $\mathbf{\Lambda}_{\xi}$ and $\mathbf{\Lambda}_{\eta}$, and the global refinement matrix is obtained using the Kronecker product,
\begin{equation}
\mathbf{\Lambda} = \mathbf{\Lambda}_{\xi} \otimes \mathbf{\Lambda}_{\eta}.
\end{equation}

Now we focus on how to obtain the weights for the two-scale relation. For $h$-refinement, various algorithms can be found in the literature, see \cite{piegl_nurbs_1997} and references therein. For $p$-refinement, several steps are required. First, the curve is split into separate Bézier elements by simple knot insertion. Second, the degree of each Bézier element is increased by one. Finally, the Bézier elements are reassembled to reconstruct the degree-elevated curve. This approach, described in \cite{piegl_nurbs_1997}, is known to be very efficient per element, but loops over the entire knot vector and computes the weights and new control points on the fly.  

For local $p$-refinement, where degree elevation of the entire curve is unnecessary, the straightforward but less efficient method of solving a linear system is applied as needed.
 While the computational cost per basis function is higher, this approach provides the flexibility to select individual basis functions of interest and avoids looping over the entire global knot vector. Especially in the case of equidistant knot vectors, many basis functions share the same shape. By using a canonical representation, in which a basis function is shifted to the origin, it is possible to easily compare basis functions and efficiently store the corresponding two-scale coefficients in a lookup table.

\begin{figure}[H]
    \newlength{\imgheight}
    \setlength{\imgheight}{2.3cm}
    \centering
    \begin{tabular}{@{\hspace{30pt}} c @{\hspace{30pt}} c @{\hspace{30pt}} c} 

        \subfloat[][]{
            \includegraphics[height=\imgheight]{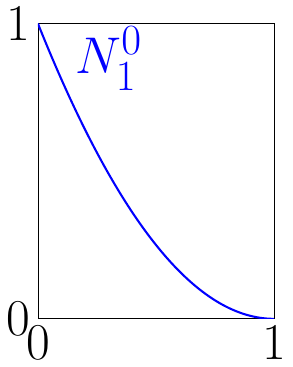}
        } &
        \subfloat[][]{
            \includegraphics[height=\imgheight]{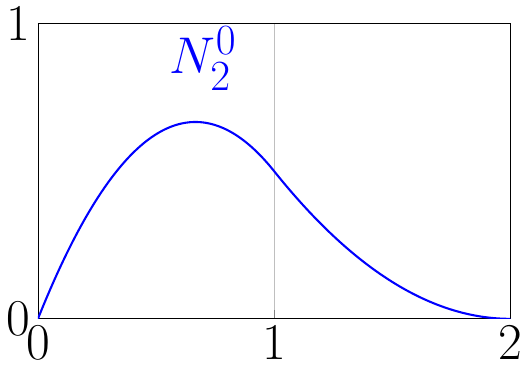}
        } &
        \subfloat[][]{
            \includegraphics[height=\imgheight]{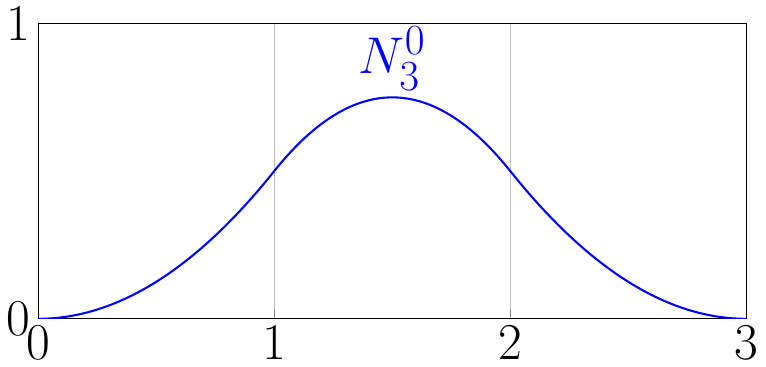}
        } \\[25pt]

        \subfloat[][]{
            \includegraphics[height=\imgheight]{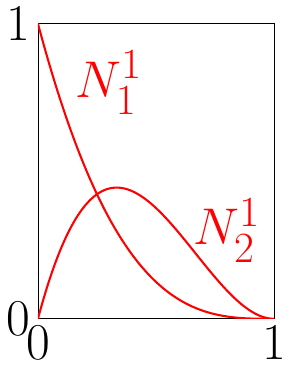}
        } &
        \subfloat[][]{
            \includegraphics[height=\imgheight]{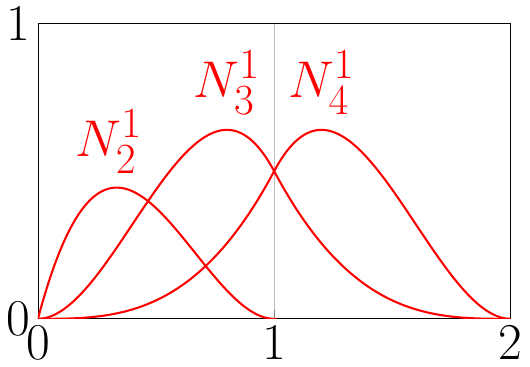}
        } &
        \subfloat[][]{
            \includegraphics[height=\imgheight]{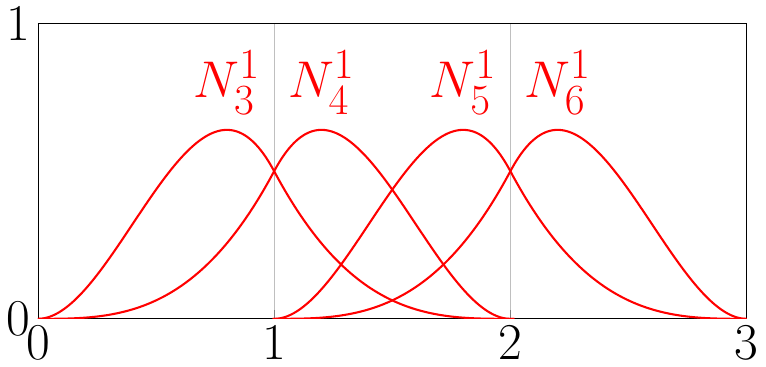}
        } \\[25pt]

        \subfloat[][]{
            \includegraphics[height=\imgheight]{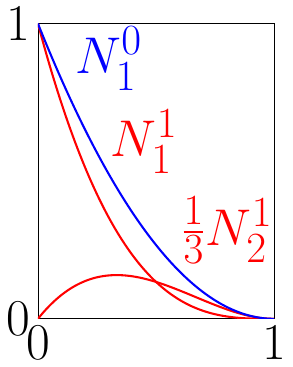}
        } &
        \subfloat[][]{
            \includegraphics[height=\imgheight]{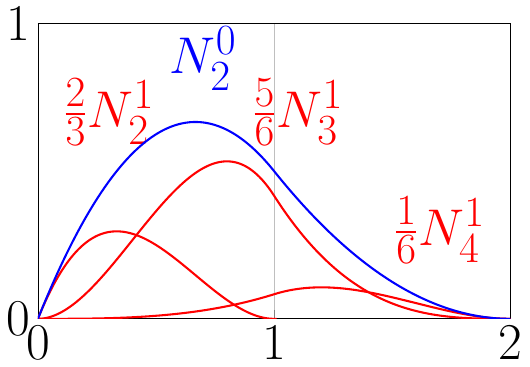}
        } &
        \subfloat[][]{
            \includegraphics[height=\imgheight]{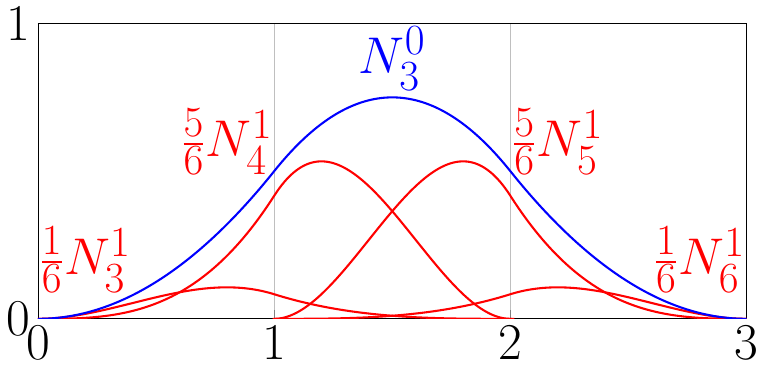}
        }
    \end{tabular}
\caption{Two-scale relation for different $p$-refined basis functions. The first row shows basis functions on level 0. The second row shows their corresponding children after $p$-refinement. The third row illustrates the scaled children whose weighted sum reproduces the parent basis function.}
    \label{fig:twoscale}
\end{figure}

Fig.~\ref{fig:twoscale} shows the two-scale relation for different basis functions for $p=2$ using $p$-refinement. The basis functions in (a), (b), and (c) are computed from the level-0 local knot vectors $\Xi_1^0 = [0,0,0,1]$, $\Xi_2^0 = [0,0,1,2]$, and $\Xi_3^0 = [0,1,2,3]$, each containing $p+2 = 4$ knots. $p$-refinement is performed by increasing the multiplicity of every knot by one. The refined basis functions are still defined by $p+2$ knots, but since the degree is elevated to $p=3$, this now corresponds to five knots per function. The intermediate knot vector $\widetilde\Xi_1^0 = [0,0,0,0,1,1]$ is split into new local knot vectors: $\Xi_1^1 = [0,0,0,0,1]$ and $\Xi_2^1 = [0,0,0,1,1]$, as shown in (d). Similarly, in (e), the intermediate vector $\widetilde\Xi_2^0 = [0,0,0,1,1,2,2]$ gives rise to the refined vectors $\Xi_2^1 = [0,0,0,1,1]$, $\Xi_3^1 = [0,0,1,1,2]$, and $\Xi_4^1 = [0,1,1,2,2]$. In (f), the refinement continues with $\widetilde\Xi_3^0 = [0,0,1,1,2,2,3,3]$ generating $\Xi_3^1 = [0,0,1,1,2]$, $\Xi_4^1 = [0,1,1,2,2]$, $\Xi_5^1 = [1,1,2,2,3]$, and $\Xi_6^1 = [1,2,2,3,3]$. In (g), (h), and (i), the two-scale relation is illustrated, showing that the weighted sum of the child basis functions reproduces the parent basis function. The $n$ coefficients in Eq.~\ref{eq:twoscale_sum} can be computed by solving a linear system of equations, obtained by evaluating the involved functions at $n$ distinct locations within the support. 
In multilevel refinement, any basis function at level $l = \hat{l}$ can be represented in terms of basis functions at levels $l > \hat{l}$. The two-scale relation forms a hierarchical tree, and the coefficients are computed by multiplying the weights along the paths from parent to descendants. Different types of refinement (e.g., $h$-, $p$-, or $k$-refinement) can be applied in sequence within this hierarchy.

\begin{figure}[h]
    \centering
    \begin{tabular}{@{\hspace{3pt}} c @{\hspace{3pt}} c @{\hspace{3pt}} c} 

        \subfloat[][Level 0]{
            \includegraphics[width=0.3\textwidth]{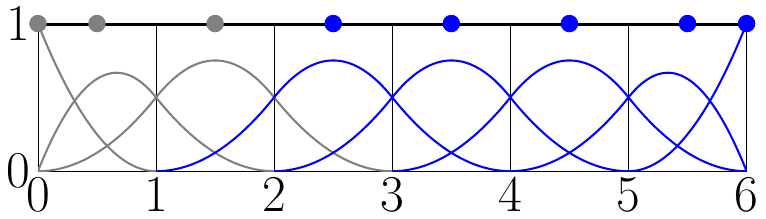}
        } &
        \subfloat[][Level 0]{
            \includegraphics[width=0.3\textwidth]{images/lvl0.pdf}
        } &
        \subfloat[][Level 0]{
            \includegraphics[width=0.3\textwidth]{images/lvl0.pdf}
        } \\[25pt] 

        \subfloat[][Level 1 after $h$-refinement]{
            \includegraphics[width=0.3\textwidth]{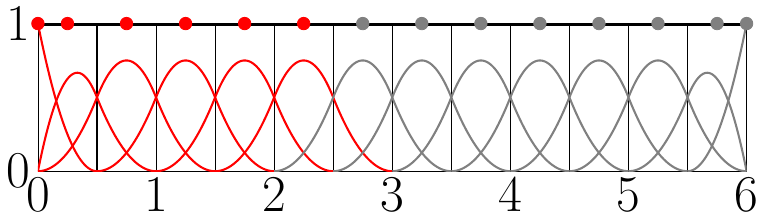}
        } &
        \subfloat[][Level 1 after $p$-refinement]{
            \includegraphics[width=0.3\textwidth]{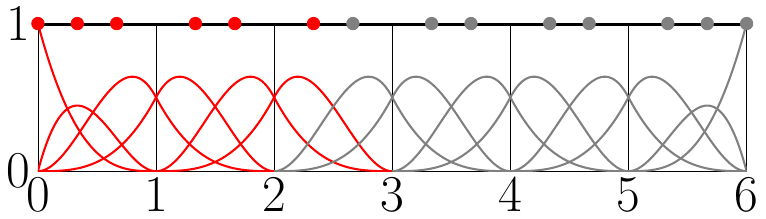}
        } &
        \subfloat[][Level 1 after $k$-refinement]{
            \includegraphics[width=0.3\textwidth]{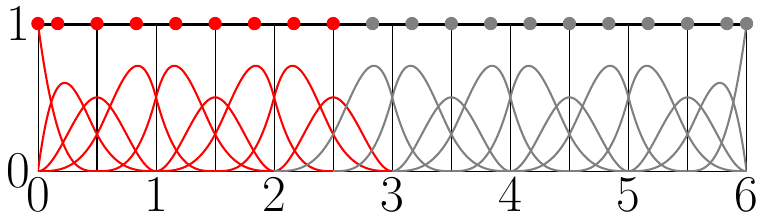}
        } \\[25pt] 

        \subfloat[][$h$-refined THB-spline basis]{
            \includegraphics[width=0.3\textwidth]{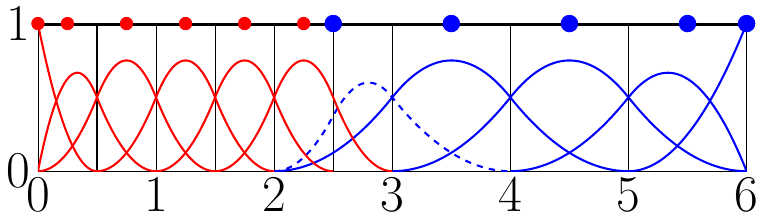}
        } &
        \subfloat[][$p$-refined THB-spline basis]{
            \includegraphics[width=0.3\textwidth]{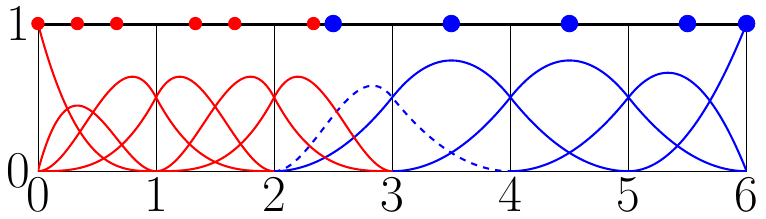}
        } &
        \subfloat[][$k$-refined THB-spline basis]{
            \includegraphics[width=0.3\textwidth]{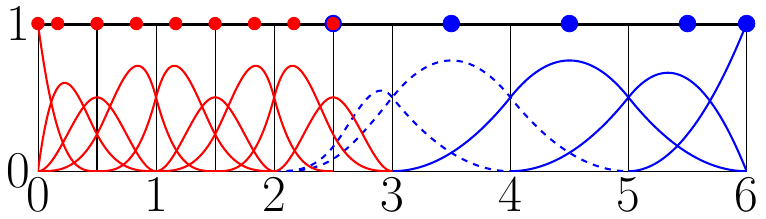}
        }
    \end{tabular}
    \caption{Comparison of refinement strategies: $h$-, $p$-, and $k$-refinement. The first row shows the initial basis with $p=2$ and respective control points for all strategies. The second row illustrates the basis functions after global refinement for each method. The third row displays the final THB-spline basis after local refinement of the leftmost three basis functions (grey), which are replaced by finer-level child basis functions (red) fully contained in their support. Level-0 basis functions partially overlapping the refined region are candidates for truncation. Truncation is performed when a child basis function is already active, ensuring partition of unity.}
    \label{fig:refinements_grid}
\end{figure}

With the two-scale relation at hand, the conceptual idea of THB-splines is identical for all refinement types. In contrast to classical FEM, IGA shifts the focus from the element to the basis function. Consequently, we refine basis functions rather than elements, which implicitly leads to element refinement. In Fig.~\ref{fig:refinements_grid}, columns 1, 2, and 3 illustrate $h$-, $p$-, and $k$-refinement, respectively. We aim to refine the three leftmost basis functions (highlighted in gray in the first row) corresponding to level 0. The first step is to replace these level-0 basis functions with level-1 basis functions that are their children, resulting in an HB-spline. However, this procedure loses the partition-of-unity property and introduces significant overlap of the coarse functions into the refined region. To restore the partition of unity, the concept of truncation is introduced. Basis functions that partially overlap the refined domain qualify for truncation. A basis function is truncated if at least one of its children is already active. In that case, the parent basis function is re-expressed via the two-scale relation, with the weights corresponding to the active children set to zero (i.e., excluded from the two-scale relation). Truncated basis functions are shown by dotted lines in (g), (h), and (i).

\begin{figure}[H]
  \centering
  \parbox{0.3\textwidth}{
    \includegraphics[width=\linewidth]{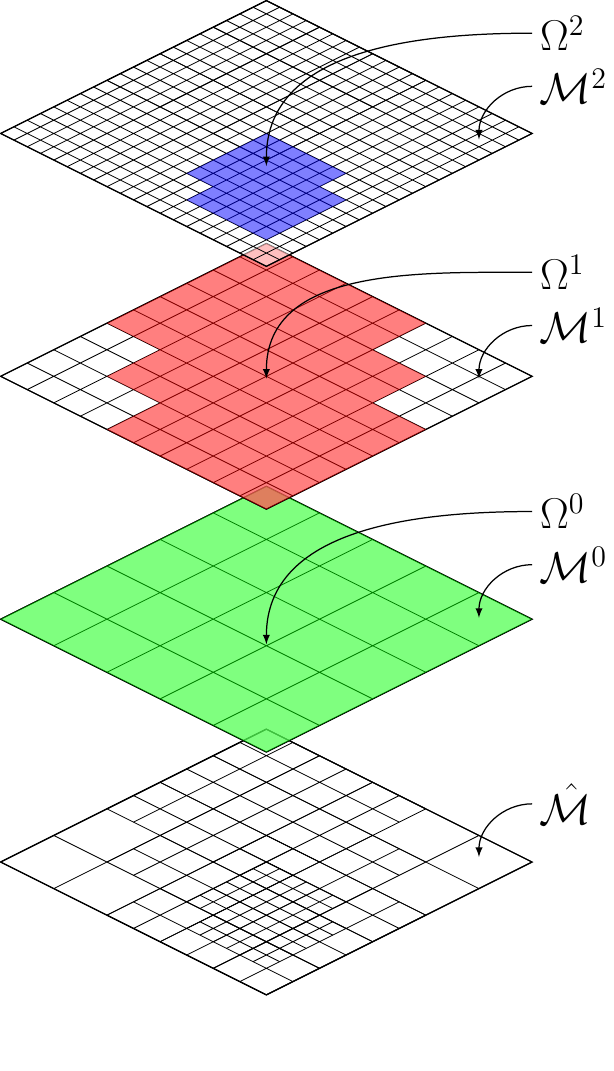}
  }
  \hspace{0.05\textwidth}
  \parbox{0.45\textwidth}{
    \centering
    \includegraphics[width=0.6\linewidth]{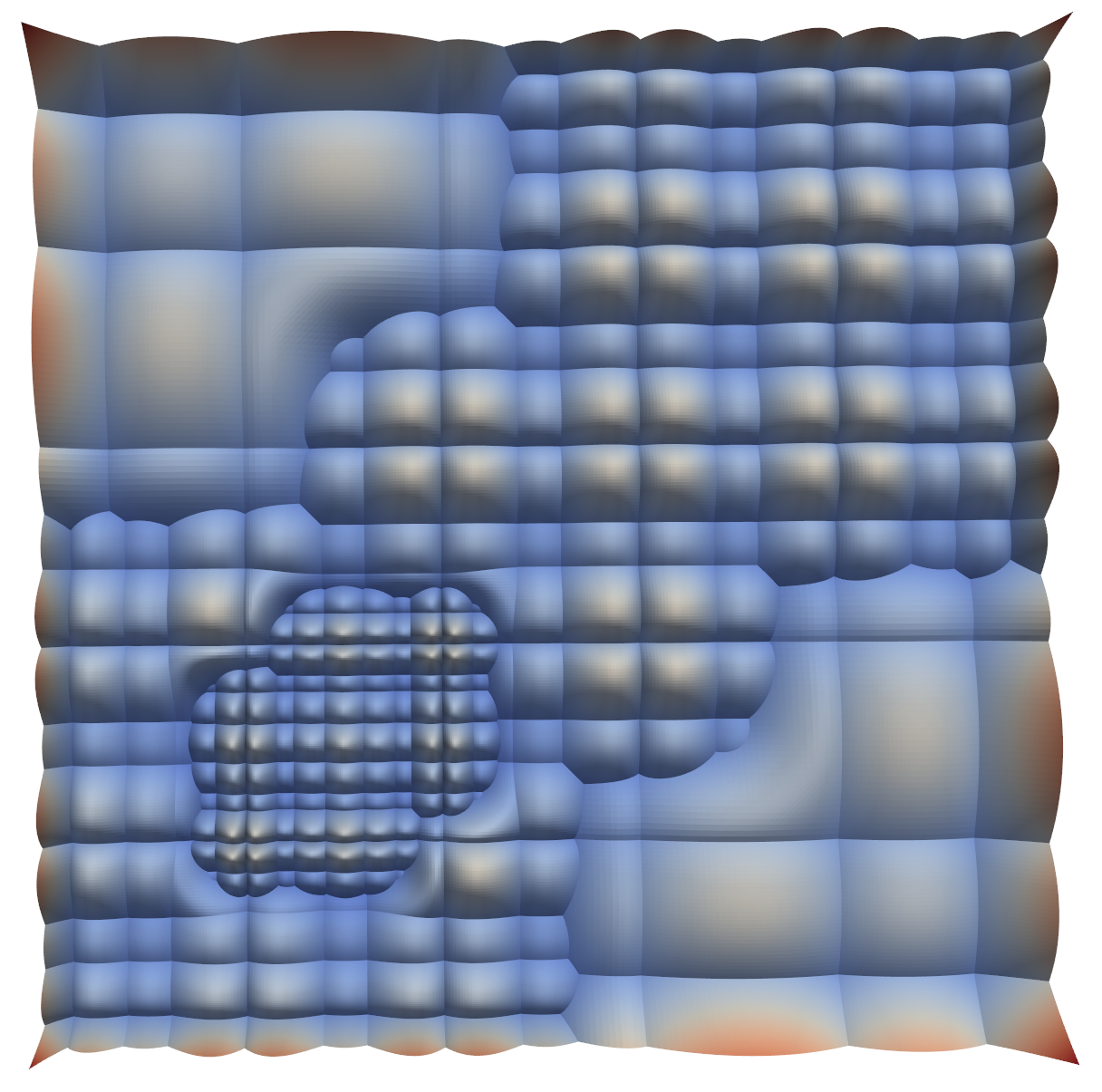}
    \\[0.5em] 
    \includegraphics[width=0.8\linewidth]{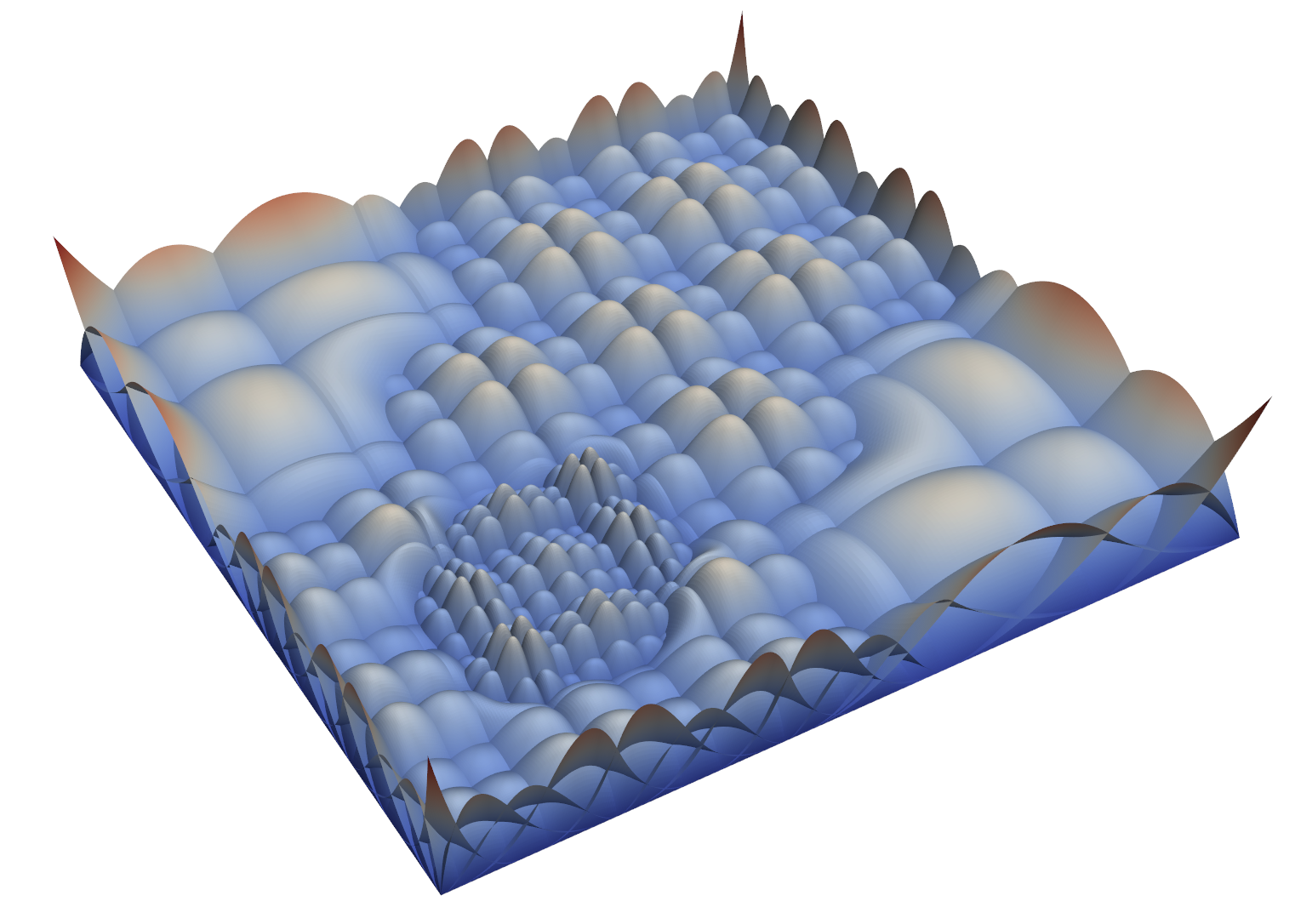}
  }
  \caption{Refinement domains and final mesh $\hat{\mathcal{M}}$ for local $k$-refinement (left) and THB basis functions after local $k$-refinement (right). The coarsest level has $p=2$ with $C^1$ continuity. The second level uses $p=3$, retaining $C^1$ continuity at existing knots and achieving $C^2$ at newly inserted ones. The third level employs $p=4$ with knot continuities ranging from $C^1$ to $C^3$.}
  \label{fig:THBmesh}
\end{figure}

Refinement in two dimensions follows the same conceptual principles as in one dimension. However, the truncation process can lead to non-convex basis function supports. Unlike the 1D case, situations may arise in 2D where a function $N_j^{l+1}$ is fully contained within the refined domain but is not a child of any deactivated function. By applying the recursive algorithm defined in Eq.~\ref{eq:recursive}, we ensure that a function at level $l > 0$ is only activated if its parent was deactivated beforehand. As a result, such functions are excluded from the final basis. For further details on efficient algorithms and implementations, we refer to \cite{DAngellaReali2020THBExtraction}. Fig.~\ref{fig:THBmesh} illustrates examples of the hierarchical mesh and the corresponding THB-basis after local $k$-refinement. Note that basis functions may be truncated multiple times, particularly in the presence of sharp transitions between refinement levels.

\subsection{The Shifted Boundary Method for the Poisson Problem} \label{SBM}
The SBM is an unfitted, embedded domain technique designed to efficiently handle complex geometries without requiring the mesh to conform to the domain boundary~\cite{main2018shifted, atallah2021shifted}. Rather than directly enforcing boundary conditions on the true boundary $\Gamma$ of a physical domain $\Omega$, SBM introduces a surrogate boundary $\widetilde{\Gamma}$, embedded within a background mesh and chosen to closely approximate $\Gamma$ while aligning with element edges. This surrogate boundary eliminates the need for geometry-conforming meshing and avoids integration over arbitrarily small cut elements, which are known to induce stability and accuracy issues in traditional immersed boundary methods.

The overall SBM procedure consists of the following steps: 
\begin{itemize}
    \item A structured background mesh is generated such that it fully contains the true domain $\Omega$.
    \item A surrogate domain $\widetilde{\Omega}$ is defined as the union of elements that best approximate $\Omega$. Its boundary, $\widetilde{\Gamma}$, is formed by the edges of the active elements that lie near $\Gamma$.
    \item A variational formulation is posed on $\widetilde{\Omega}$, and boundary conditions are imposed not on $\Gamma$, but on $\widetilde{\Gamma}$, using a Taylor expansion to shift the data from the true to the surrogate boundary.
\end{itemize}

{\color{red}
More specifically, the surrogate domain $\widetilde{\Omega}$ is constructed as follows: 
First, we identify the background elements intersected by the true boundary $\Gamma$ (cut elements). 
All elements that are strictly inside $\Omega$ and not intersected by $\Gamma$ are marked as active and are therefore included in $\widetilde{\Omega}$. 
For cut elements, we adopt the standard \emph{optimal} criterion used in SBM~\cite{yang2023optimal}: an element is added to the active set if the portion of its area belonging to $\Omega$ exceeds the inactive portion; otherwise it is discarded. 
The surrogate boundary $\widetilde{\Gamma}$ is then defined as the union of the outer edges of the active elements adjacent to the boundary. 
Further implementation details can be found in~\cite{main2018shifted,antonelli2024shifted}. In the present work, the body force is defined analytically on the entire background domain through the manufactured solution. 
}

\begin{figure}[H]
    \centering
    \vspace{0.5em}
      \begin{minipage}{0.9\textwidth}
        \centering
        \includegraphics[width=0.40\linewidth]{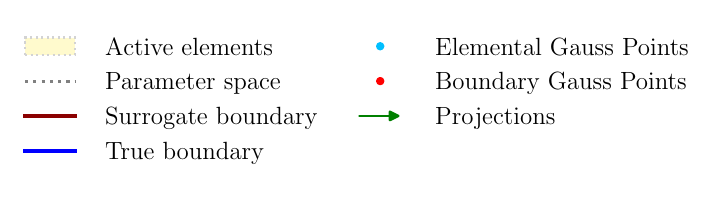}
      \end{minipage}
      
    \begin{subfigure}[b]{0.45\textwidth}
        \centering
        \includegraphics[width=\textwidth, trim={0.5cm, 0.5cm, 0.5cm, 0.5cm}, clip]
        {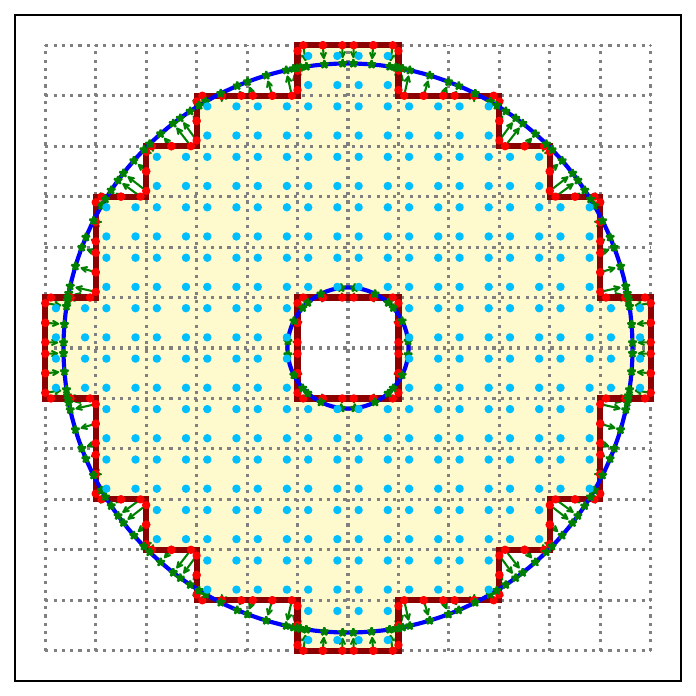}
        \caption{Surrogate and true boundaries with active elements and Gauss points.}
        \label{fig:sbm-full}
    \end{subfigure}
    \hspace{0.05\textwidth}
    \begin{subfigure}[b]{0.45\textwidth}
        \centering
        \includegraphics[width=\textwidth]{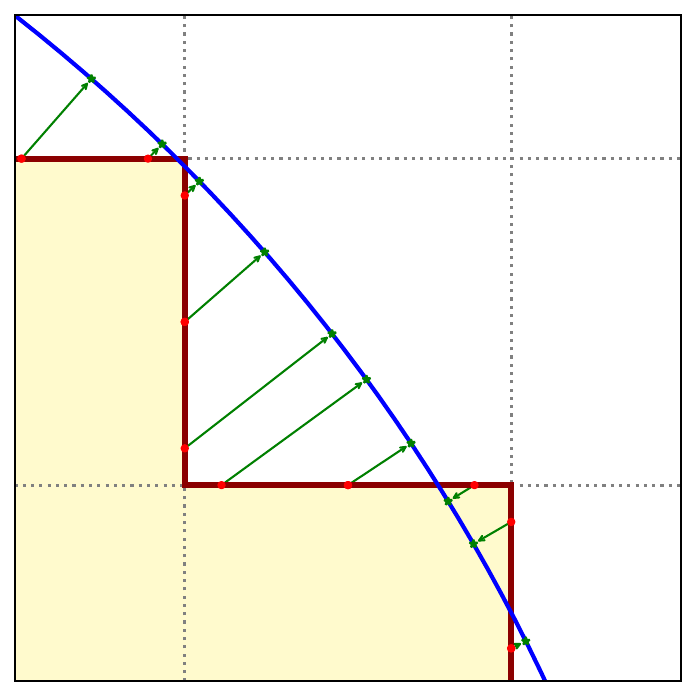}
        \caption{Zoomed view of projections from surrogate to true boundary.}
        \label{fig:sbm-zoom}
    \end{subfigure}
    \caption{Visualization of the SBM for B-splines of order $p=1$. (a) Full domain with active elements and boundary geometry. (b) Detail showing projection of boundary Gauss points onto the true geometry.}
    \label{fig:sbm-overview}
\end{figure}

Fig.~\ref{fig:sbm-overview} illustrates the geometry of the method. The surrogate boundary $\widetilde{\Gamma}$ (in red) lies entirely within the background mesh, while the true boundary $\Gamma$ (in blue) may intersect elements arbitrarily. The projection from $\widetilde{\Gamma}$ to $\Gamma$ is defined via a closest point mapping.

We now describe the SBM formulation for the scalar Poisson problem, where the unknown field $u$ represents temperature. Let the true domain be $\Omega \subset \mathbb{R}^2$, with boundary $\Gamma = \Gamma_D \cup \Gamma_N$, where Dirichlet and Neumann conditions are imposed. The surrogate domain and boundary are denoted by $\widetilde{\Omega}$ and $\widetilde{\Gamma}$, respectively.
We define the closest-point projection mapping
\begin{align}
    \mathbb{M} : \widetilde{\Gamma} \to \Gamma, \\
    \xTilde \mapsto \Bx,
\end{align}
which maps each point $\xTilde \in \widetilde{\Gamma}$ to its closest counterpart $\Bx \in \Gamma$. 
{\color{red} The projection distance vector is
\begin{equation}
    \bold{d}(\xTilde) = \Bx - \xTilde = [\mathbb{M} - \mathbb{I}](\xTilde).
\end{equation}

To impose boundary conditions on $\widetilde{\Gamma}$ instead of $\Gamma$, we approximate the pullback of the solution and its gradient using Taylor expansions. Then, the shift operators of order $m$ for the solution and its gradient are defined as:
\begin{align} 
    \mathbb{S}_{\bold{d}}^m(u(\xTilde)) 
    &= u(\xTilde) + \sum_{i = 1}^{m} \frac{D^i_{\bold{d}} u(\xTilde)}{i!}
    \approx u(\mathbb{M}(\xTilde)) = u(\Bx), \label{eq:shift_u_orig} \\
    \mathbb{S}_{\bold{d}}^m(\nabla u(\xTilde)) 
    &= \nabla u(\xTilde) + \sum_{i = 1}^{m} \frac{D^i_{\bold{d}} \nabla u(\xTilde)}{i!}
    \approx \nabla u(\mathbb{M}(\xTilde)) = \nabla u(\Bx), \label{eq:shift_grad_orig}
\end{align}
where $D^i_{\bold{d}}(\cdot)$ denotes the $i$-th order directional derivative. Let $\bold{d} = (d_1,d_2) \in \mathbb{R}^2$ and let $\boldsymbol{\alpha} = (\alpha_1,\alpha_2) \in \mathbb{N}^2$ be a multi-index with the notation $\bold{d}^{\boldsymbol{\alpha}} = d_1^{\alpha_1} d_2^{\alpha_2}$, $|\boldsymbol{\alpha}| = \alpha_1 + \alpha_2$, and $\boldsymbol{\alpha}! = \alpha_1!\alpha_2!$. The directional derivative operator is explicitly given by:
\begin{equation} \label{eq:D(.)}
    D^i_{\bold{d}}(\cdot) = \sum_{|\boldsymbol{\alpha}| = i} \frac{i!}{\boldsymbol{\alpha}!} \frac{\partial^{i} (\cdot) }{\partial \Bx^{\boldsymbol{\alpha}}} \bold{d}^{\bold{\alpha}}.
\end{equation}
By inserting \eqref{eq:D(.)} into the definitions \eqref{eq:shift_u_orig} and \eqref{eq:shift_grad_orig}, the $i!$ terms cancel out. Note that for $|\boldsymbol{\alpha}|=0$, the term $\frac{\bold{d}^{\boldsymbol{\alpha}}}{\boldsymbol{\alpha}!} \frac{\partial^0 u(\xTilde)}{\partial \Bx^{\boldsymbol{\alpha}}}$ recovers the original function $u(\xTilde)$. Thus, the shift operators are written in the simplified multi-index form:
\begin{align} 
    \mathbb{S}_{\bold{d}}^m(u(\xTilde)) 
    &= \sum_{|\boldsymbol{\alpha}| \le m} \frac{\bold{d}^{\boldsymbol{\alpha}}}{\boldsymbol{\alpha}!} \frac{\partial^{|\boldsymbol{\alpha}|} u(\xTilde)}{\partial \Bx^{\boldsymbol{\alpha}}} 
    \approx u(\Bx), \label{eq:shift_u_simplified} \\
    \mathbb{S}_{\bold{d}}^m(\nabla u(\xTilde)) 
    &= \sum_{|\boldsymbol{\alpha}| \le m} \frac{\bold{d}^{\boldsymbol{\alpha}}}{\boldsymbol{\alpha}!} \frac{\partial^{|\boldsymbol{\alpha}|} \nabla u(\xTilde)}{\partial \Bx^{\boldsymbol{\alpha}}} 
    \approx \nabla u(\Bx), \label{eq:shift_grad_simplified}
\end{align}
with $\boldsymbol{\alpha} = (\alpha_1,\alpha_2) \in \mathbb{N}_0^2$ also including zero.
These operators allow us to approximate the solution behavior on $\Gamma$ using only information defined on $\widetilde{\Gamma}$.

For Dirichlet boundary conditions, the maximum admissible shift order $m$
coincides with the polynomial degree $p$ of the underlying spline basis,
since derivatives of order higher than $p$ vanish identically.
Hence, we set $m=p$ and impose the Dirichlet data $u_D$ as
\begin{equation}
    u_D(\mathbf{x})
    = \mathbb{S}_{\mathbf{d}}^{p}(u(\widetilde{\mathbf{x}}))
      + R^{p}(u,\mathbf{d})(\widetilde{\mathbf{x}})
    \approx \mathbb{S}_{\mathbf{d}}^{p}(u(\widetilde{\mathbf{x}})),
\end{equation}
where the remainder satisfies
\[
|R^{p}(u,\mathbf{d})|
= \mathcal{O}\bigl(\|\mathbf{d}\|^{p+1}\bigr).
\]
Since the optimal convergence rate of an IGA discretization of degree $p$
is $\mathcal{O}(h^{p+1})$ and the projection distance satisfies
\[
\|\mathbf{d}\| = \mathcal{O}(h) \quad \text{as } h \to 0,
\]
it follows that
\[
R^{p}(u,\mathbf{d})
= \mathcal{O}(h^{p+1}).
\]
Therefore, the SBM imposition of Dirichlet boundary conditions preserves
the optimal convergence rate. Consequently, employing locally refined
THB-splines in this setting is not expected to provide a substantial
accuracy improvement.

For Neumann conditions, the maximal admissible shift order in the discrete expansion is $m = p-1$. Given the prescribed flux $t_N$ on $\Gamma$, we obtain
\begin{equation}
    t_N(\mathbf{x})
    =
    \bigl[
        \mathbb{S}_{\mathbf{d}}^{p-1}(\nabla u(\widetilde{\mathbf{x}}))
        +
        R^{p-1}(\nabla u,\mathbf{d})(\widetilde{\mathbf{x}})
    \bigr]
    \cdot \mathbf{n}
    \approx
    \bigl[
        \mathbb{S}_{\mathbf{d}}^{p-1}(\nabla u(\widetilde{\mathbf{x}}))
    \bigr]
    \cdot \mathbf{n},
\end{equation}
where $\mathbf{n}$ denotes the outward unit normal on $\Gamma$. The corresponding remainder satisfies
\[
|R^{p-1}(\nabla u,\mathbf{d})|
=
\mathcal{O}\bigl(\|\mathbf{d}\|^{p}\bigr).
\]
Since $\|\mathbf{d}\| = \mathcal{O}(h)$, it follows that
\[
R^{p-1}(\nabla u,\mathbf{d})
=
\mathcal{O}(h^{p}),
\]
which suggests a one-order reduction of the asymptotic convergence rate
compared to the optimal $\mathcal{O}(h^{p+1})$ behaviour.}

{\color{red} 
\begin{rem}
In practice, for each quadrature point 
$\widetilde{\mathbf{x}} \in \widetilde{\Gamma}$, 
the mapping $\mathbb{M}$ is constructed by computing a closest point on the true boundary $\Gamma$ with respect to the Euclidean distance, 
in accordance with the standard SBM formulation~\cite{main2018shifted, antonelli2024shifted}. For sufficiently smooth boundaries, the closest-point projection is locally unique. Configurations in which the projection is not unique are not analyzed in detail in the present study.
\end{rem}
}

{\color{red}
\begin{rem}
For sufficiently smooth solutions, the above argument highlights a
structural difference between Dirichlet and Neumann boundary conditions
in the SBM formulation.
For Dirichlet data, the shift operator can be applied up to order $p$,
which formally yields a Taylor remainder of order
$\mathcal{O}(\|\mathbf{d}\|^{p+1})$
under the standard shift expansion.
In contrast, for Neumann conditions, the maximal admissible shift order
in the discrete expansion is $p-1$. Consequently, the corresponding
Taylor remainder scales as
$\mathcal{O}(\|\mathbf{d}\|^{p})$.
Since $\|\mathbf{d}\| = \mathcal{O}(h)$ under mesh refinement,
this indicates that the boundary-consistency contribution in the
Neumann case may scale as
$\mathcal{O}(h^{p})$,
compared to
$\mathcal{O}(h^{p+1})$
for Dirichlet data.
This suggests a potential one-order reduction in the asymptotic
convergence rate relative to the body-fitted case or to the Dirichlet
SBM formulation, unless additional correction strategies are introduced.
\end{rem}
}

\subsection{THB-splines for the SBM} \label{THB for SBM}

THB-splines are widely used for adaptive refinement in isogeometric analysis~\cite{Buffa2022_AdaptiveIGA,GIANNELLI2016,gu_adaptive_2018,patrizi_adaptive_2020}, typically by refining the basis functions supported in elements exhibiting large error indicators or estimators~\cite{kumar_simple_2015}.  

Within the SBM, the primary source of error {\color{red} for smooth problems} is the Taylor‐series approximation on the surrogate boundary. Accordingly, we can target refinement a priori: the only basis functions that need to be refined are those whose supports intersect the surrogate boundary. Refining one such basis function affects all knot spans in its support, which for a B-spline of degree $p$, $C^{p-1}$ continuous, form a rectangular patch of size $(p+1)\times(p+1)$.  

\begin{figure}[h]
    \centering
    \begin{minipage}[b]{0.4\textwidth}
        \centering
        \includegraphics[width=\linewidth,
                         trim={0.5cm 0.5cm 0.5cm 0.5cm}, clip]
                         {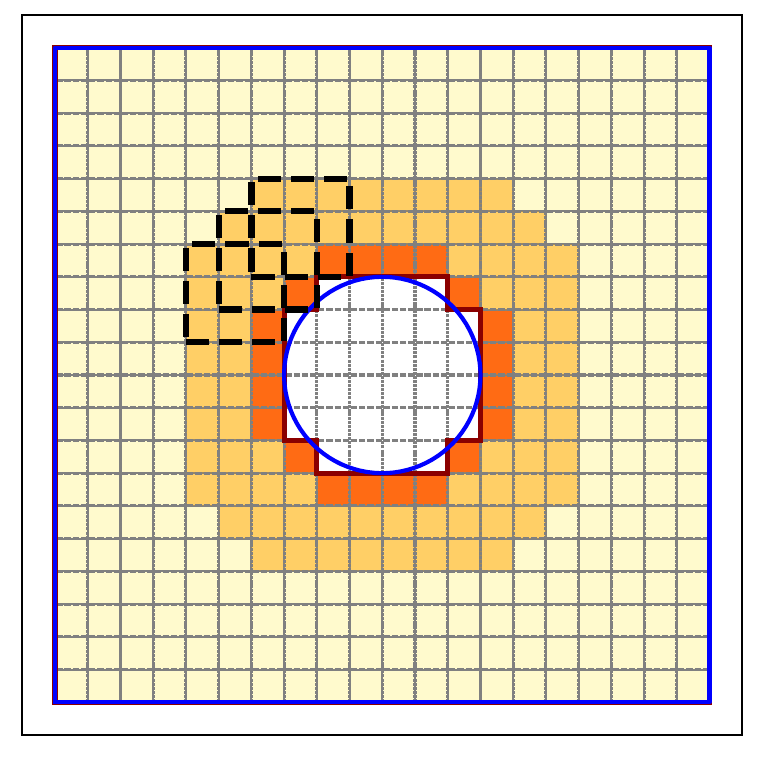}
    \end{minipage}
    \hspace{0.05\textwidth}   
    \begin{minipage}[b]{0.4\textwidth}
        \centering
        \includegraphics[width=\linewidth,
                         trim={0.5cm 0.5cm 0.5cm 0.5cm}, clip]
                         {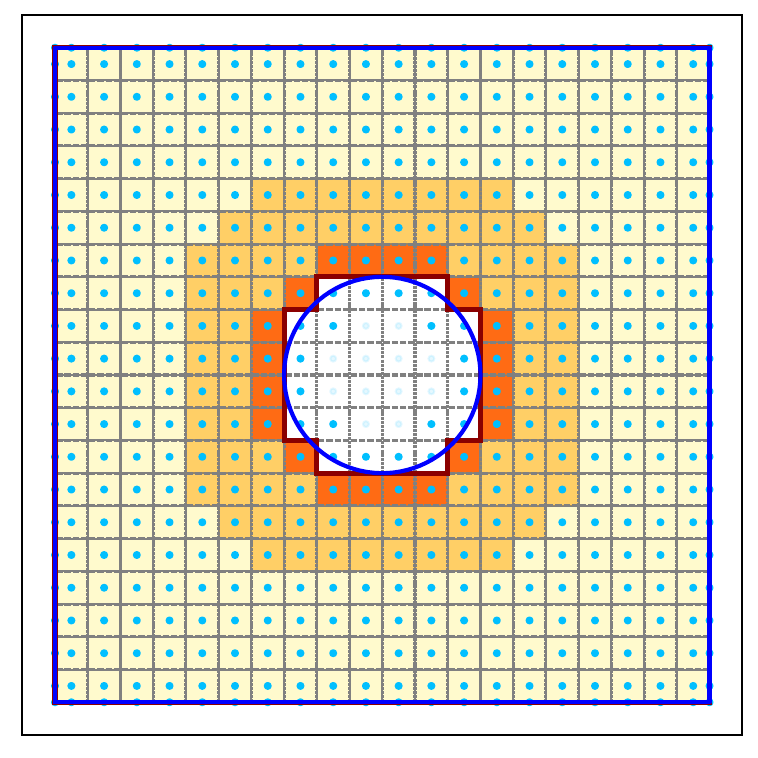}
    \end{minipage}
    \caption{Surrogate domain and boundary for an immersed circle with
             standard B-splines of degree $p=2$.  
             Left: dashed squares indicate the joint support of some
             basis functions active on the surrogate boundary that will
             be refined.  
             Right: same configuration with the corresponding
             active control points (azure). The true boundary (blue), active knot
             spans (light yellow), boundary knot spans (dark orange) and
             the combined support region (light orange) are highlighted.}
    \label{fig:sbm-thb}
\end{figure}

\begin{figure}[h]
    \centering
    \setlength{\tabcolsep}{2pt}
    \renewcommand{\arraystretch}{1.0}
    \begin{tabular}{ccc}
        \subcaptionbox{Local $h$-refinement}{%
            \includegraphics[width=0.30\linewidth,
                             trim={0.5cm 0.5cm 0.5cm 0.5cm}, clip]
                             {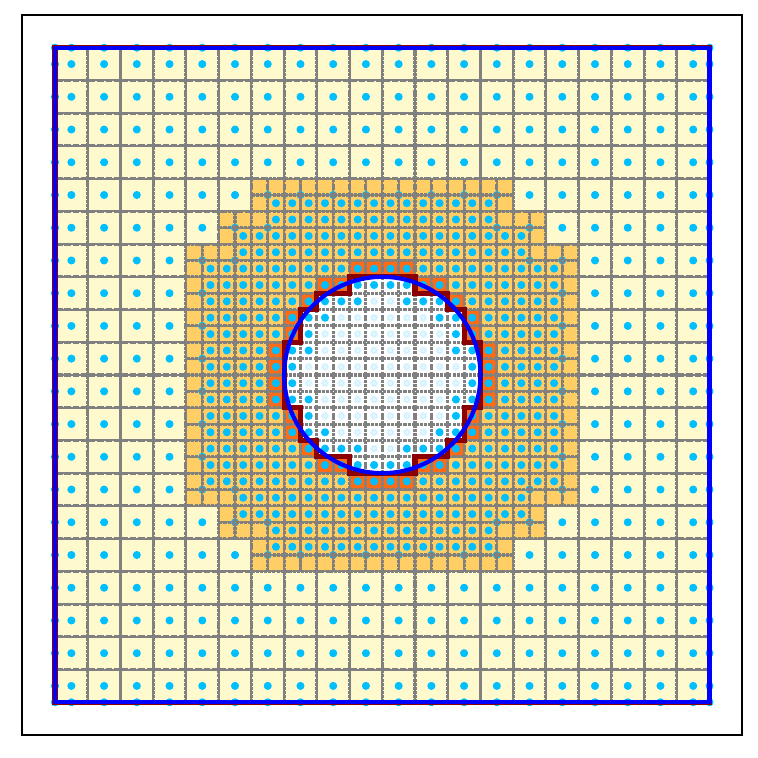}}
        &
        \subcaptionbox{Local $p$-refinement}{%
            \includegraphics[width=0.30\linewidth,
                             trim={0.5cm 0.5cm 0.5cm 0.5cm}, clip]
                             {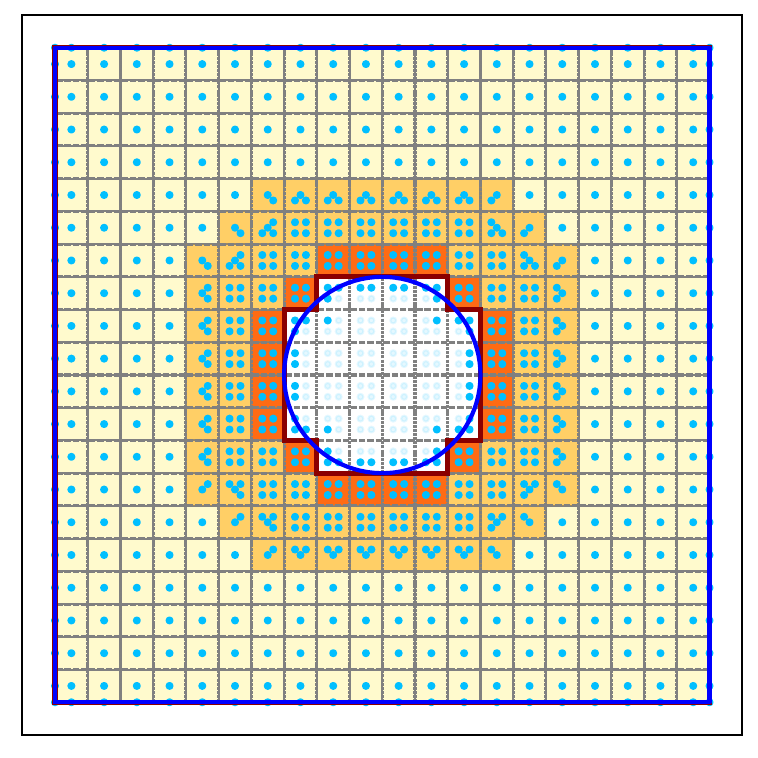}}
        &
        \subcaptionbox{Local $k$-refinement}{%
            \includegraphics[width=0.30\linewidth,
                             trim={0.5cm 0.5cm 0.5cm 0.5cm}, clip]
                             {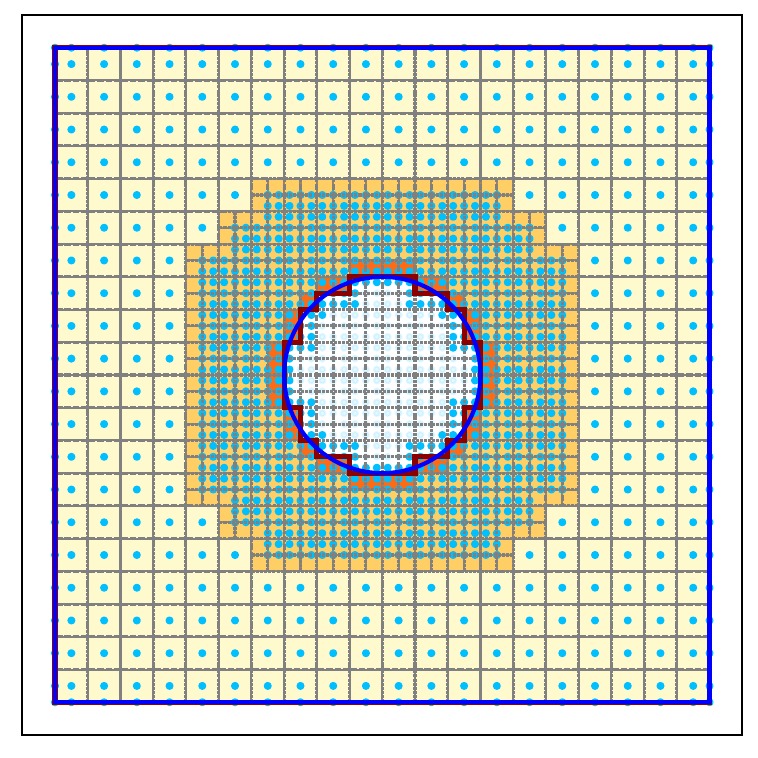}}
    \end{tabular}
    \caption{Surrogate domain for an immersed circle with THB-splines of degree $p=2$ after three different local-refinement strategies.  Each panel shows the active knot spans and the corresponding active control points for (a) $h$-, (b) $p$- and (c) $k$-refinement.}
    \label{fig:sbm-three-refinements}
\end{figure}

{\color{red} Fig.~\ref{fig:sbm-thb} illustrates the active knot spans affected by local refinement for the representative example of an immersed circular hole in a square domain. Fig.~\ref{fig:sbm-three-refinements} shows the resulting active knot spans and control points after one step of $h$-, $p$-, and $k$-refinement, respectively. In all subsequent studies involving local refinement within the SBM, no adaptive refinement strategy is employed. Instead, a single refinement step is applied in a controlled manner to investigate under which conditions it is sufficient to recover the desired convergence behavior.} To preserve integration accuracy, the quadrature rule must be adjusted on each knot span according to the highest polynomial degree of the basis functions supported on it. In particular, the orange‐marked spans in Fig.\ref{fig:sbm-three-refinements} in the case of $p$- and $k$-refinement carry more quadrature points than the rest of the mesh. 

{\color{red}
The following considerations are heuristic and intended to motivate the numerical study. They do not constitute a theoretical analysis of the SBM. Under this perspective, the qualitative effects of the different refinement schemes are expected to behave as follows:

\begin{itemize}
  \item \emph{$h$-refinement} improves the geometric approximation of the surrogate boundary, thereby reducing the projection distance from $\widetilde{\Gamma}$ to $\Gamma$ and decreasing the boundary-consistency term especially for Neumann conditions. In boundary-dominated regimes, this may reduce the overall error and delay the onset of suboptimal convergence behaviour.

  \item \emph{$p$-refinement} leaves the boundary location unchanged but increases the maximal admissible shift order in the Taylor expansion. For Neumann boundary conditions, where the boundary operator involves derivatives of the solution, this local enrichment may reduce the boundary-consistency contribution related to the truncation of the shift expansion, without globally increasing the polynomial degree.

  \item \emph{$k$-refinement} combines local degree elevation with mesh refinement, simultaneously improving the geometric approximation of the surrogate boundary and enriching the Taylor expansion.
\end{itemize}

}

\subsection{Enhanced Shift Operators}
\label{ssec:enhanced_shift}
{\color{red}

In addition to the classical Taylor shift operator defined in \eqref{eq:shift_u_simplified} and \eqref{eq:shift_grad_simplified}, which corresponds to a Taylor expansion of total order $m$ along the direction $\bold{d}$, we introduce a more general tensorial variant.

Let $\bold{d} = (d_1,d_2) \in \mathbb{R}^2$ and let $\boldsymbol{\alpha} = (\alpha_1,\alpha_2) \in \mathbb{N}_0^2$ be a multi-index. For sufficiently smooth fields, we first define the enhanced shift operators with independent orders $m_1$ and $m_2$ for each coordinate direction:
\begin{align} 
    \widetilde{\mathbb{S}}_{\bold{d}}^{(m_1, m_2)}(u(\xTilde)) 
    &:= \sum_{\alpha_1=0}^{m_1} \sum_{\alpha_2=0}^{m_2} \frac{\bold{d}^{\boldsymbol{\alpha}}}{\boldsymbol{\alpha}!} \frac{\partial^{|\boldsymbol{\alpha}|} u(\xTilde)}{\partial \Bx^{\boldsymbol{\alpha}}} 
    \approx u(\Bx), \label{eq:enhanced_shift_u_general} \\
    \widetilde{\mathbb{S}}_{\bold{d}}^{(m_1, m_2)}(\nabla u(\xTilde)) 
    &:= \sum_{\alpha_1=0}^{m_1} \sum_{\alpha_2=0}^{m_2} \frac{\bold{d}^{\boldsymbol{\alpha}}}{\boldsymbol{\alpha}!} \frac{\partial^{|\boldsymbol{\alpha}|} \nabla u(\xTilde)}{\partial \Bx^{\boldsymbol{\alpha}}} 
    \approx \nabla u(\Bx). \label{eq:enhanced_shift_grad_general}
\end{align}
In the following, we consider the case where the expansion order is equal in both directions, $m_1 = m_2 = m$. The operators then read:
\begin{align} 
    \widetilde{\mathbb{S}}_{\bold{d}}^{m}(u(\xTilde)) 
    &:= \sum_{0 \le \alpha_1, \alpha_2 \le m} \frac{\bold{d}^{\boldsymbol{\alpha}}}{\boldsymbol{\alpha}!} \frac{\partial^{|\boldsymbol{\alpha}|} u(\xTilde)}{\partial \Bx^{\boldsymbol{\alpha}}} 
    \approx u(\Bx), \label{eq:enhanced_shift_u_p} \\
    \widetilde{\mathbb{S}}_{\bold{d}}^{m}(\nabla u(\xTilde)) 
    &:= \sum_{0 \le \alpha_1, \alpha_2 \le m} \frac{\bold{d}^{\boldsymbol{\alpha}}}{\boldsymbol{\alpha}!} \frac{\partial^{|\boldsymbol{\alpha}|} \nabla u(\xTilde)}{\partial \Bx^{\boldsymbol{\alpha}}} 
    \approx \nabla u(\Bx). \label{eq:enhanced_shift_grad_p}
\end{align}
In contrast to the classical shift operator, this formulation retains all partial derivatives up to order $m$ in each parametric direction. It should be emphasized that this construction is not a Taylor expansion of total order $m$. While it contains all terms with $|\boldsymbol{\alpha}| \le m$ of the standard operator, it additionally includes mixed derivatives of higher total order, up to $|\boldsymbol{\alpha}| = 2m$, due to the tensor-product index set $0 \le \alpha_1, \alpha_2 \le m$.
In the discrete setting, the field \(u_h\) is represented as
\[
  u_h(\mathbf{x})
  =
  \sum_{A} T_A(\mathbf{x})\, u_A,
\]
where \(T_A\) denotes the spline basis functions and \(u_A\) the
corresponding coefficients. Due to the linearity of the operator, the enhanced shift applied to the discrete field of degree $p$ at the boundary reads
\[
  \widetilde{\mathbb{S}}_{\mathbf{d}}^{p}
  \bigl(u_h(\widetilde{\mathbf{x}})\bigr)
  =
  \sum_A
  \widetilde{\mathbb{S}}_{\mathbf{d}}^{p}
  \bigl(T_A(\widetilde{\mathbf{x}})\bigr)
  u_A.
\]
Hence, it suffices to evaluate the operator at the level of individual
basis functions.

\medskip\noindent
\textbf{Example (\(p=1\)).}
For improved readability, the two parametric directions are denoted
by \(x\) and \(y\) instead of the indices $1$ and $2$.

For a basis function \(T_A\), the classical first-order SBM operator yields
\[
  \mathbb{S}_{\mathbf{d}}^{1}
  \bigl(T_A(\widetilde{\mathbf{x}})\bigr)
  =
  T_A(\widetilde{\mathbf{x}})
  + d_x\,\partial_x T_A(\widetilde{\mathbf{x}})
  + d_y\,\partial_y T_A(\widetilde{\mathbf{x}}).
\]
Bivariate spline basis functions of degree $p=1$ are bilinear on each
element and may contain terms of the form $xy$. Consequently,
the mixed derivative $\partial_x \partial_y T_A$ is in general nonzero,
whereas pure second derivatives vanish identically.

The enhanced operator therefore contains the additional mixed term
\[
  \widetilde{\mathbb{S}}_{\mathbf{d}}^{1}
  \bigl(T_A(\widetilde{\mathbf{x}})\bigr)
  =
  T_A(\widetilde{\mathbf{x}})
  + d_x\,\partial_x T_A(\widetilde{\mathbf{x}})
  + d_y\,\partial_y T_A(\widetilde{\mathbf{x}})
  + d_x d_y\,\partial_x\partial_y T_A(\widetilde{\mathbf{x}}).
\]
Although this includes a second-order mixed derivative,
it does not reproduce the full second-order Taylor expansion,
since the pure second derivatives are zero for bilinear functions. We will use the following operator choice in the remainder based on the empirical study in Sec.~\ref{sec:operator_comparison}.}

\subsection{Final Formulation}
Throughout this work, the parametric and physical domains are taken to coincide. In the context of the SBM, no additional geometric mapping is required, so the identity map suffices. For a discussion of more general (non-coincident) choices in standard IGA settings, see \cite{hughes2005isogeometric}.
Let the physical domain be $\Omega\subset\mathbb{R}^2$ with boundary
$\Gamma=\Gamma_D\cup\Gamma_N$ (Dirichlet and Neumann portions, disjoint).
Following the SBM, a surrogate domain
$\widetilde{\Omega}$ is built from the background mesh, with surrogate
boundary $\widetilde{\Gamma}$ and outward normal
$\widetilde{\mathbf{n}}$.  The closest-point projection
\(\mathbb{M}:\widetilde{\Gamma}\!\to\!\Gamma\) induces the partitions
\[
  \widetilde{\Gamma}_D :=
  \bigl\{\;\widetilde{\mathbf{x}}\in\widetilde{\Gamma}\;|\;
        \mathbb{M}(\widetilde{\mathbf{x}})\in\Gamma_D\bigr\},
  \qquad
  \widetilde{\Gamma}_N :=
  \bigl\{\;\widetilde{\mathbf{x}}\in\widetilde{\Gamma}\;|\;
        \mathbb{M}(\widetilde{\mathbf{x}})\in\Gamma_N\bigr\}.
\]
{\color{red}
Let \(\hat{\mathcal{M}}\) be a hierarchically refined Cartesian background mesh covering \(\widetilde{\Omega}\). We define the discrete space as
\[
V_h(\widetilde{\Omega}) :=
\mathrm{span}\{\, T_A \,\}_{A\in\mathcal{I}}
\subset H^1(\widetilde{\Omega}),
\]
where \(\{T_A\}_{A\in\mathcal{I}}\) denotes the THB-spline basis constructed over \(\hat{\mathcal{M}}\).
In the context of the SBM, only basis functions whose support intersects the surrogate domain $\widetilde{\Omega}$ are retained in \(\mathcal{I}\). Each element \(e \in \hat{\mathcal{M}}\) is associated with a unique hierarchical level \(l\) and the corresponding polynomial degree \(p^l\). Owing to the hierarchical construction, only basis functions from levels less than or equal to \(l\) can be active on \(e\). The restriction of the basis functions to an element remains a tensor-product polynomial, and the maximal local polynomial degree is determined by the level, i.e.,
\[
T_A|_e \in \mathcal{Q}^{p^l}(e).
\]
The same conforming space is employed for both trial and test functions. For a point \(\widetilde{\mathbf{x}}\in\widetilde{\Gamma}\) let
\(\mathbf{d}(\widetilde{\mathbf{x}})=\mathbb{M}(\widetilde{\mathbf{x}})-\widetilde{\mathbf{x}}\).
Assuming that all basis functions whose support intersects the
surrogate boundary $\widetilde{\Gamma}$ have the same polynomial
degree $p$, the Taylor--shift operators
\[
  \mathbb{S}_{\mathbf{d}}^{p}(u_h), \qquad
  \mathbb{S}_{\mathbf{d}}^{p-1}(\nabla u_h)
\] approximate the
pull-backs \(u_h\bigl(\mathbb{M}(\widetilde{\mathbf{x}})\bigr)\) and
\(\nabla u_h\bigl(\mathbb{M}(\widetilde{\mathbf{x}})\bigr)\),
respectively.
The discrete formulation of the Poisson problem for the SBM reads:
\medskip\noindent
Find $u_h \in V_h(\widetilde{\Omega})$ such that
\[
a_h(u_h, v_h) = \ell_h(v_h)
\quad \text{for all } v_h \in V_h(\widetilde{\Omega}),
\]
where
\begin{subequations}\label{eq:sbm-nitsche}
\begin{align}
  a_h(u_h,v_h) &= (\nabla u_h,\nabla v_h)_{\widetilde{\Omega}}
  - \langle \nabla u_h\!\cdot\!\widetilde{\mathbf{n}},\,v_h\rangle_{\widetilde{\Gamma}_D}
  - \theta\,\langle \mathbb{S}_{\mathbf{d}}^{p}(u_h),\nabla v_h\!\cdot\!\widetilde{\mathbf{n}}\rangle_{\widetilde{\Gamma}_D}
  + \frac{\alpha}{h}\,\langle\mathbb{S}_{\mathbf{d}}^{p}(u_h),v_h\rangle_{\widetilde{\Gamma}_D} \notag\\
  &\quad - \langle \nabla u_h\!\cdot\!\widetilde{\mathbf{n}},\,v_h\rangle_{\widetilde{\Gamma}_N}
  + \langle (\mathbb{S}_{\mathbf{d}}^{p-1}(\nabla u_h)\!\cdot\!\mathbf{n})(\widetilde{\mathbf{n}}\!\cdot\!\mathbf{n}),\,v_h\rangle_{\widetilde{\Gamma}_N}
\label{eq:sbm-nitsche-b} \\[4pt]
  \ell_h(v_h) &= (f,v_h)_{\widetilde{\Omega}}
          - \theta\,\langle u_D,\nabla v_h\!\cdot\!\widetilde{\mathbf{n}}\rangle_{\widetilde{\Gamma}_D}
  + \frac{\alpha}{h}\,\langle u_D,v_h\rangle_{\widetilde{\Gamma}_D}
          +\langle t_N\,(\widetilde{\mathbf{n}}\!\cdot\!\mathbf{n}),\,v_h\rangle_{\widetilde{\Gamma}_N}.
\label{eq:sbm-nitsche-c}
\end{align}
\end{subequations}
Here, $(\cdot,\cdot)_{\widetilde{\Omega}}$ denotes the $L^2$ inner product over $\widetilde{\Omega}$, while $\langle \cdot , \cdot \rangle_{\widetilde{\Gamma}_D}$ and $\langle \cdot , \cdot \rangle_{\widetilde{\Gamma}_N}$ denote the $L^2$ inner product on the surrogate boundary $\widetilde{\Gamma}_D$ and $\widetilde{\Gamma}_N$, respectively. The parameter $h$ is a representative local mesh size. 
The vector $\mathbf{n}$ is the outward unit normal on the true boundary, evaluated at $\mathbb{M}(\widetilde{\mathbf{x}})$. The bilinear form is consistent with the original Poisson problem in the sense that the exact solution of the continuous problem satisfies the discrete formulation when evaluated on the surrogate boundary. The linear functional $\ell_h(\cdot)$ contains the body force $f$ on $\widetilde{\Omega}$, the prescribed Neumann flux $t_N$ on $\Gamma_N$, and the Nitsche terms associated with the Dirichlet data $u_D$ on $\Gamma_D$. Setting $\theta = -1$ and $\alpha = 0$ in \eqref{eq:sbm-nitsche} yields the penalty-free formulation, which is used for all subsequent analysis in this paper.
For a detailed derivation and analysis, see
\cite{main2018shifted,antonelli2024shifted}.}

\clearpage

\section{Numerical examples} \label{sec:numerical_examples}
This section gathers a series of numerical experiments designed to investigate the use of THB-splines within the SBM. All tests employ the penalty-free formulation. {\color{red} After verifying the correctness of the local refinement implementation, we compare the classical and enhanced shift operators. All numerical studies involving local refinement are intentionally restricted to a single refinement step to isolate its impact on convergence and to maintain a controlled comparison between refinement strategies. In particular, when local $p$- or $k$-refinement is applied, all basis functions whose support intersects the surrogate boundary are uniformly elevated to the same polynomial degree, such that the above assumption of a homogeneous boundary degree remains satisfied.}

{\color{red}In all simulations where a body-fitted reference solution is considered, this reference is obtained through a \emph{pseudo body-fitted} approach, as in~\cite{antonelli2024shifted}. In this setting, the exact boundary conditions are imposed directly on the surrogate boundary $\widetilde{\Gamma}$, without applying any shift. This eliminates the boundary-consistency error associated with the SBM while retaining the same background mesh. This choice allows for a fair comparison with the SBM formulation while avoiding the need for a separate geometrically exact mesh.

For the convergence analysis, a sequence of meshes is generated by reparameterizing the domain: at each step, the number of uniformly sized elements is increased by one in each coordinate direction while keeping the geometry fixed. This procedure corresponds to rebuilding an equidistant knot vector with higher resolution. Local refinement at the boundary is applied subsequently. The resulting mesh defines the background discretization for the SBM. Since the surrogate boundary is derived from the background mesh, its position relative to the true boundary changes slightly between successive reconstructions. This enables a controlled assessment of the convergence behavior of the SBM.
}


Convergence is assessed with respect to the manufactured exact
solution defined on the surrogate computational domain
$\widetilde{\Omega}$.
The error between the numerical solution $u_h$ and the exact solution $u$,
measured in the $L^2(\widetilde{\Omega})$ and
$H^1(\widetilde{\Omega})$ norms, is defined as
\begin{subequations}\label{eq:error_norms}
\begin{align}
    \|u - u_h\|_{L^2(\widetilde{\Omega})}
        &= \Bigl(\int_{\widetilde{\Omega}} (u - u_h)^2
                 \,\mathrm d\tilde\Omega\Bigr)^{1/2},\\[4pt]
    \|u - u_h\|_{H^1(\widetilde{\Omega})}
        &= \Bigl(
           \|u - u_h\|^2_{L^2(\widetilde{\Omega})}
           + \|\nabla u - \nabla u_h\|^2_{L^2(\widetilde{\Omega})}
           \Bigr)^{1/2}.
\end{align}
\end{subequations}


Here,
\[
\|\nabla v\|_{L^2(\widetilde{\Omega})}
:=
\Bigl(
\sum_{i=1}^{d}
\|\partial_{x_i} v\|_{L^2(\widetilde{\Omega})}^2
\Bigr)^{1/2},
\]
where $d$ denotes the spatial dimension. In this paper, we report relative errors, obtained by dividing
the quantities defined in Eq.~\eqref{eq:error_norms}
by the corresponding $L^2(\widetilde{\Omega})$
or $H^1(\widetilde{\Omega})$ norm of the exact solution $u$. All integrals in Eq.~\eqref{eq:error_norms} are evaluated by $(p+1)$ Gauss quadrature points in each direction, where $p$ is the degree of the element.

We begin by examining the isolated effects of localized $h$-, $p$-, and $k$-refinements in an untrimmed body-fitted setting, establishing a baseline for each refinement mode. Next, we quantify the stability improvement obtained with the enhanced
shift operator. We then transition to immersed-boundary
configurations, namely:
\begin{itemize}
    \item an immersed hole (Fig.~\ref{fig:sbm-thb}), with Dirichlet boundary conditions,
    \item an immersed hole (Fig.~\ref{fig:sbm-thb}), with Neumann boundary conditions,
    \item and an immersed object with a central hole (Fig.~\ref{fig:sbm-overview}), subject to external Dirichlet and internal Neumann conditions.
\end{itemize}
These cases allow us to assess how local refinement influences convergence when the surrogate boundary deviates from the true geometry. Finally, we present additional experiments on more complex immersed geometries to demonstrate the stability and robustness of the THB-SBM framework under increased geometric complexity.
{\color{blue}
All numerical experiments were performed using an in-house Python implementation based on the formulation described in this work.
}

\subsection{Verification of local refinement}
\begin{figure}[h]
    \centering
    \includegraphics[width=0.6\textwidth]{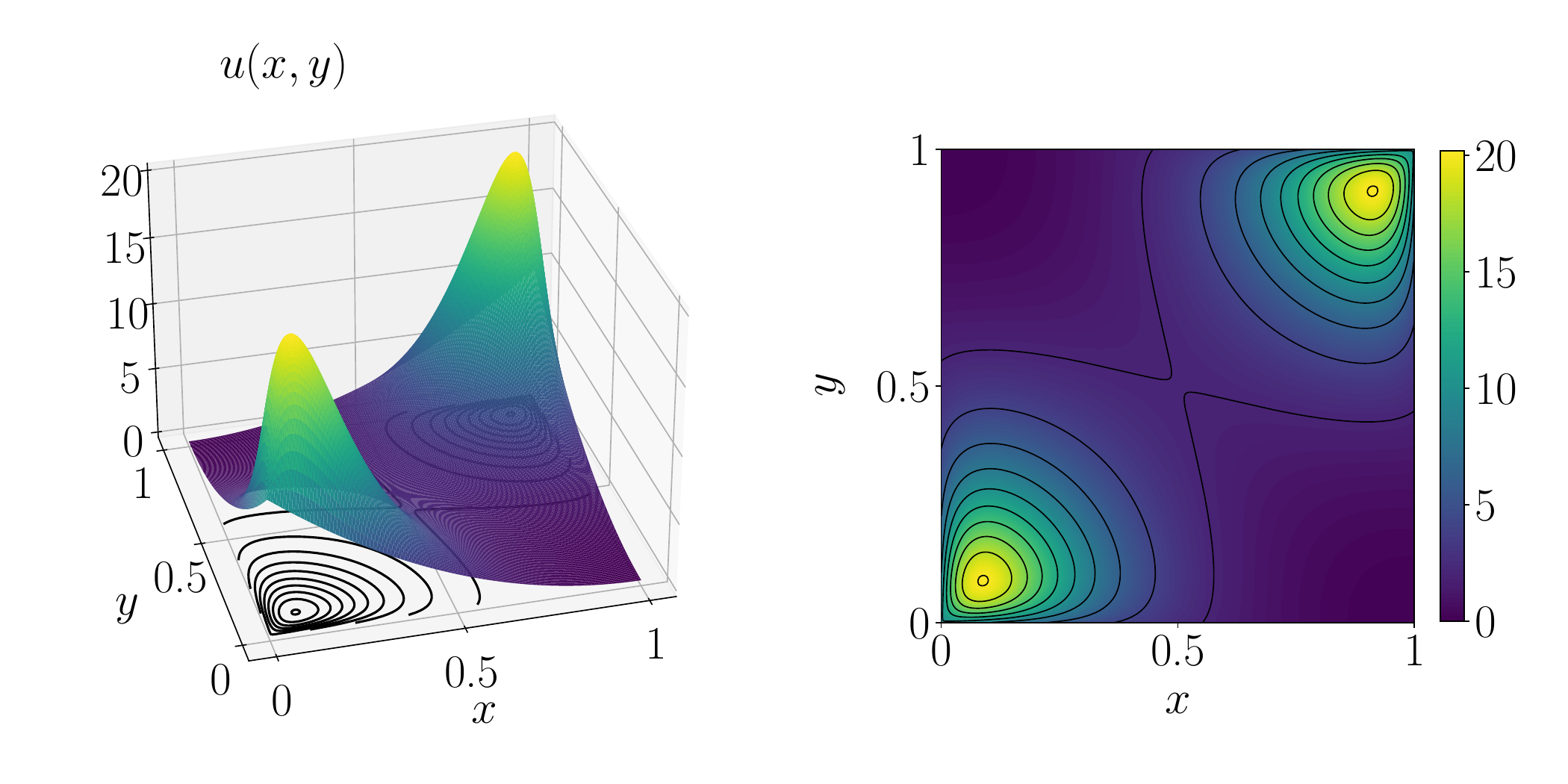}
    \caption{Manufactured solution 
    $u(x, y) = x y\,(1 - x)(1 - y) \left[\sinh\left(10 (1 - x)(1 - y)\right) + \sinh\left(10 x y\right)\right] + 10 \left[(1 - x)^2 (1 - y)^2 + x^2 y^2\right]$.}
    \label{fig:sdsf}
\end{figure}

\begin{figure}[h]
    \centering

    \begin{subfigure}[t]{0.3\textwidth}
        \centering
        \includegraphics[width=\textwidth]{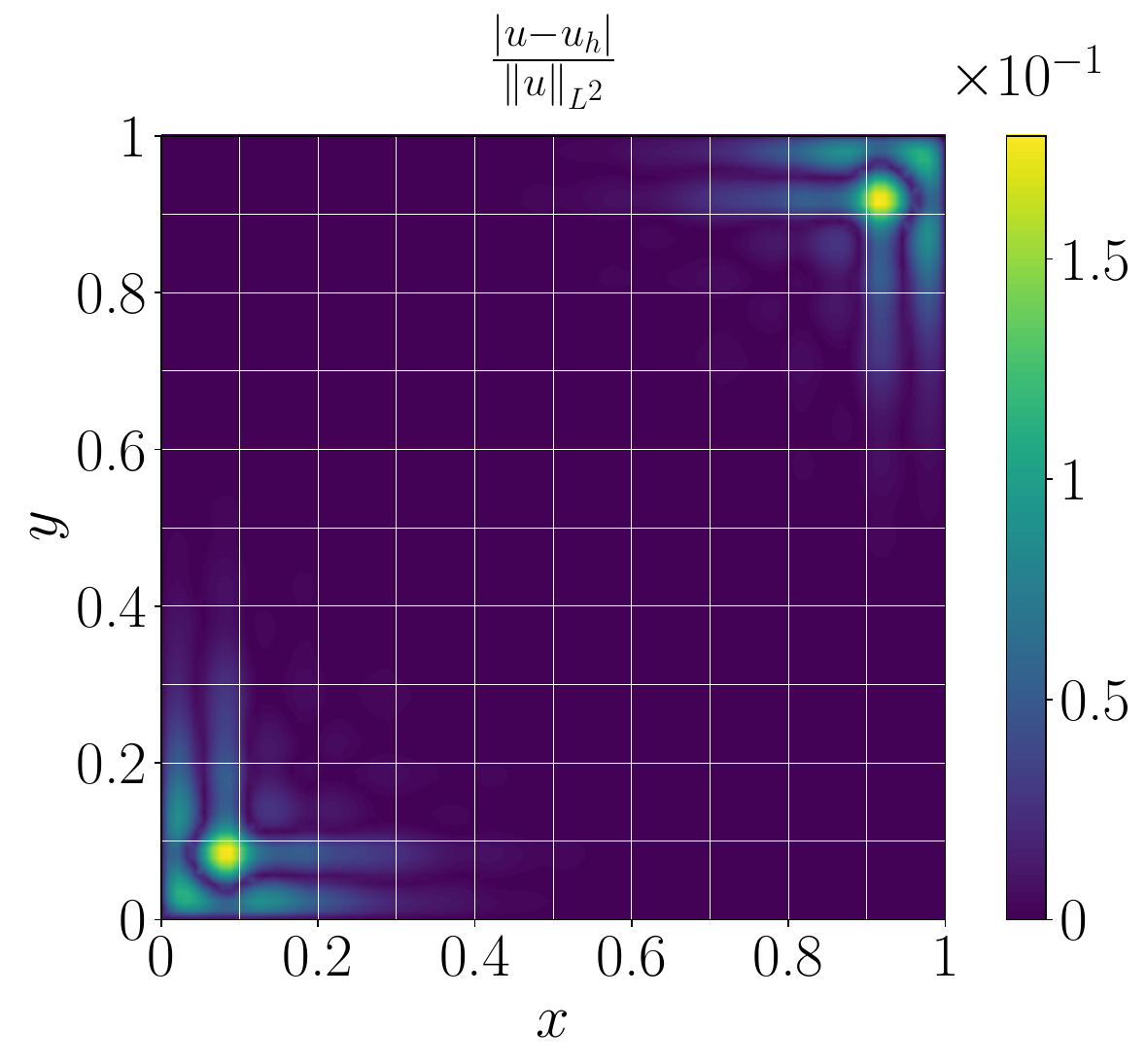}
        \caption{no refinement,\\ 144 DOFs, 100 Elements}
        \label{fig:loc_ref_a}
    \end{subfigure}
    \hfill
    \begin{subfigure}[t]{0.3\textwidth}
        \centering
        \includegraphics[width=\textwidth]{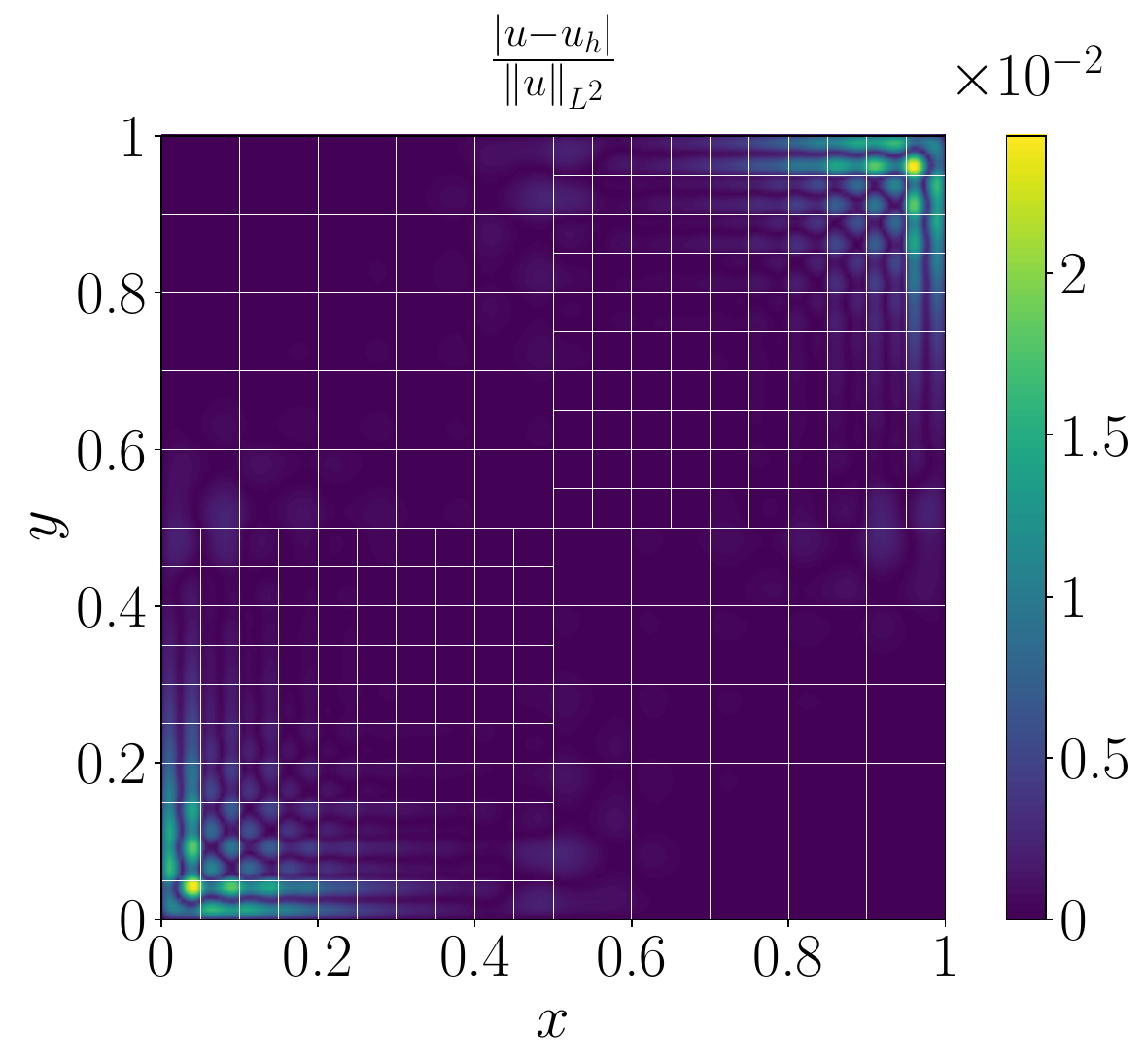}
        \caption{$h$-refinement (1 step),\\ 294 DOFs, 250 Elements}
        \label{fig:loc_ref_b}
    \end{subfigure}
    \hfill
    \begin{subfigure}[t]{0.3\textwidth}
        \centering
        \includegraphics[width=\textwidth]{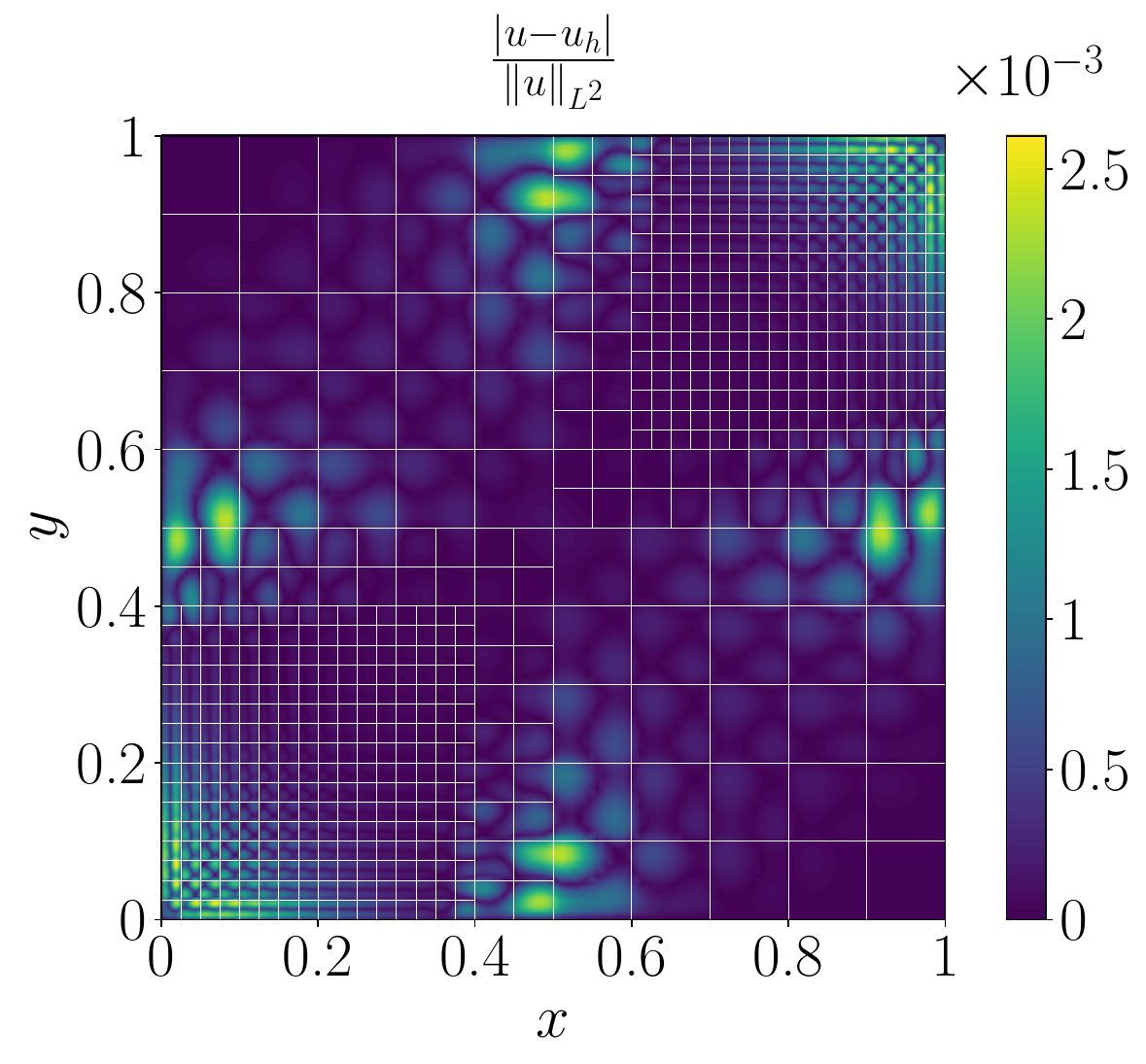}
        \caption{$h$-refinement (2 steps),\\ 678 DOFs, 634 Elements}
        \label{fig:loc_ref_c}
    \end{subfigure}

    \vspace{1em}

    \begin{subfigure}[t]{0.3\textwidth}
        \centering
        \includegraphics[width=\textwidth]{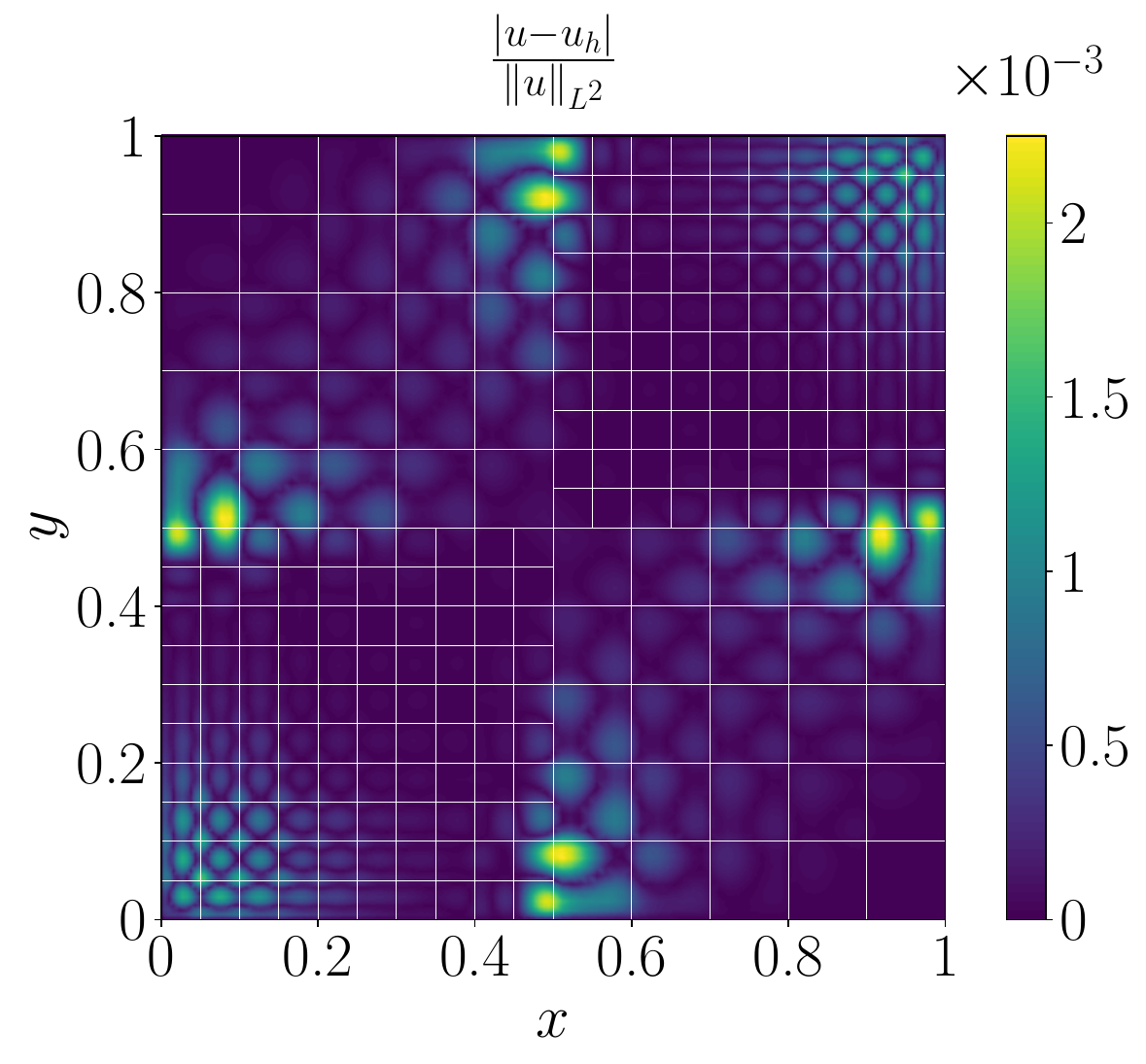}
        \caption{$k$-refinement (1 step),\\ 544 DOFs, 250 Elements}
        \label{fig:loc_ref_d}
    \end{subfigure}
    \hfill
    \begin{subfigure}[t]{0.3\textwidth}
        \centering
        \includegraphics[width=\textwidth]{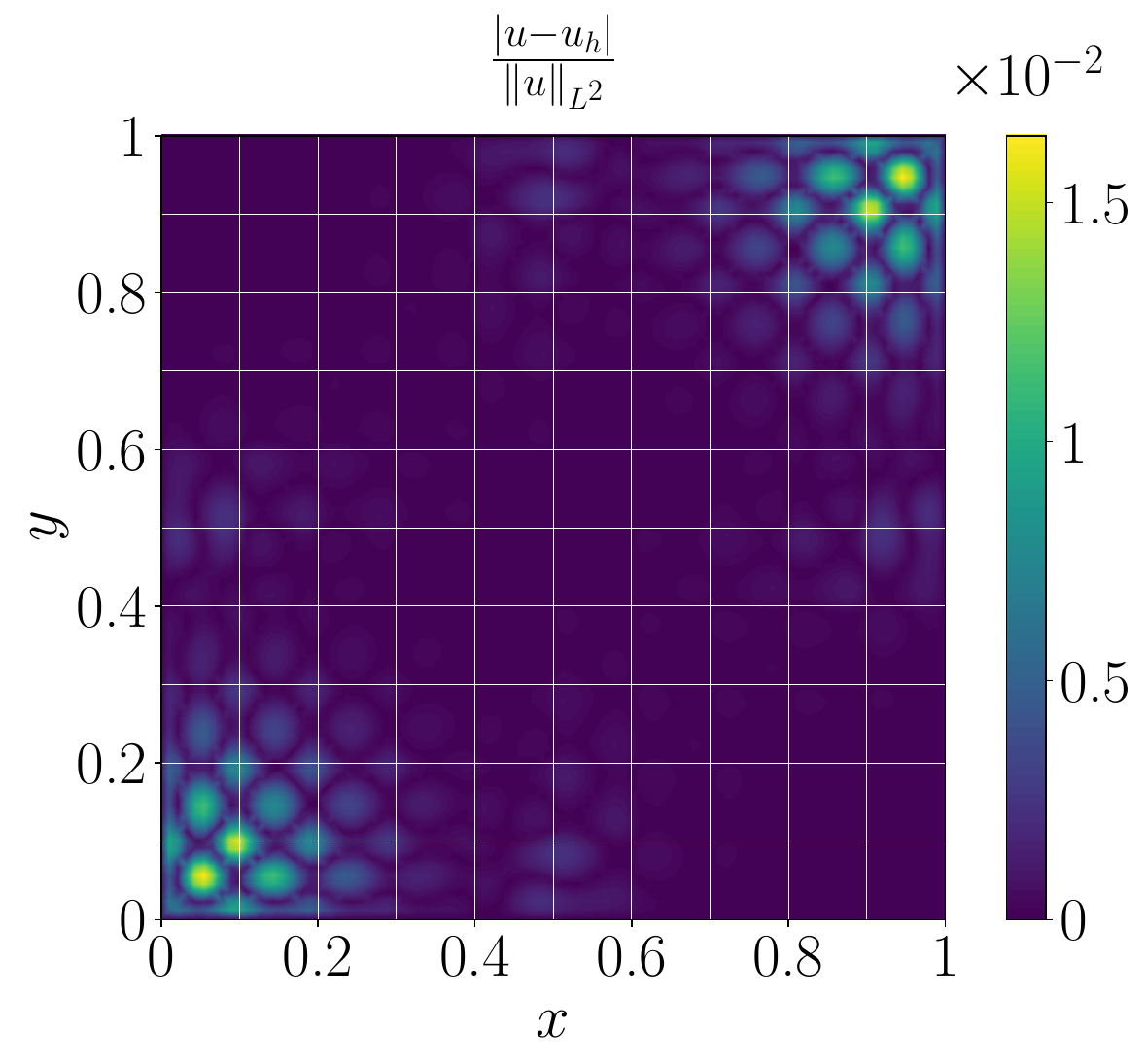}
        \caption{$p$-refinement (1 step),\\ 294 DOFs, 100 Elements}
        \label{fig:loc_ref_e}
    \end{subfigure}
    \hfill
    \begin{subfigure}[t]{0.3\textwidth}
        \centering
        \includegraphics[width=\textwidth]{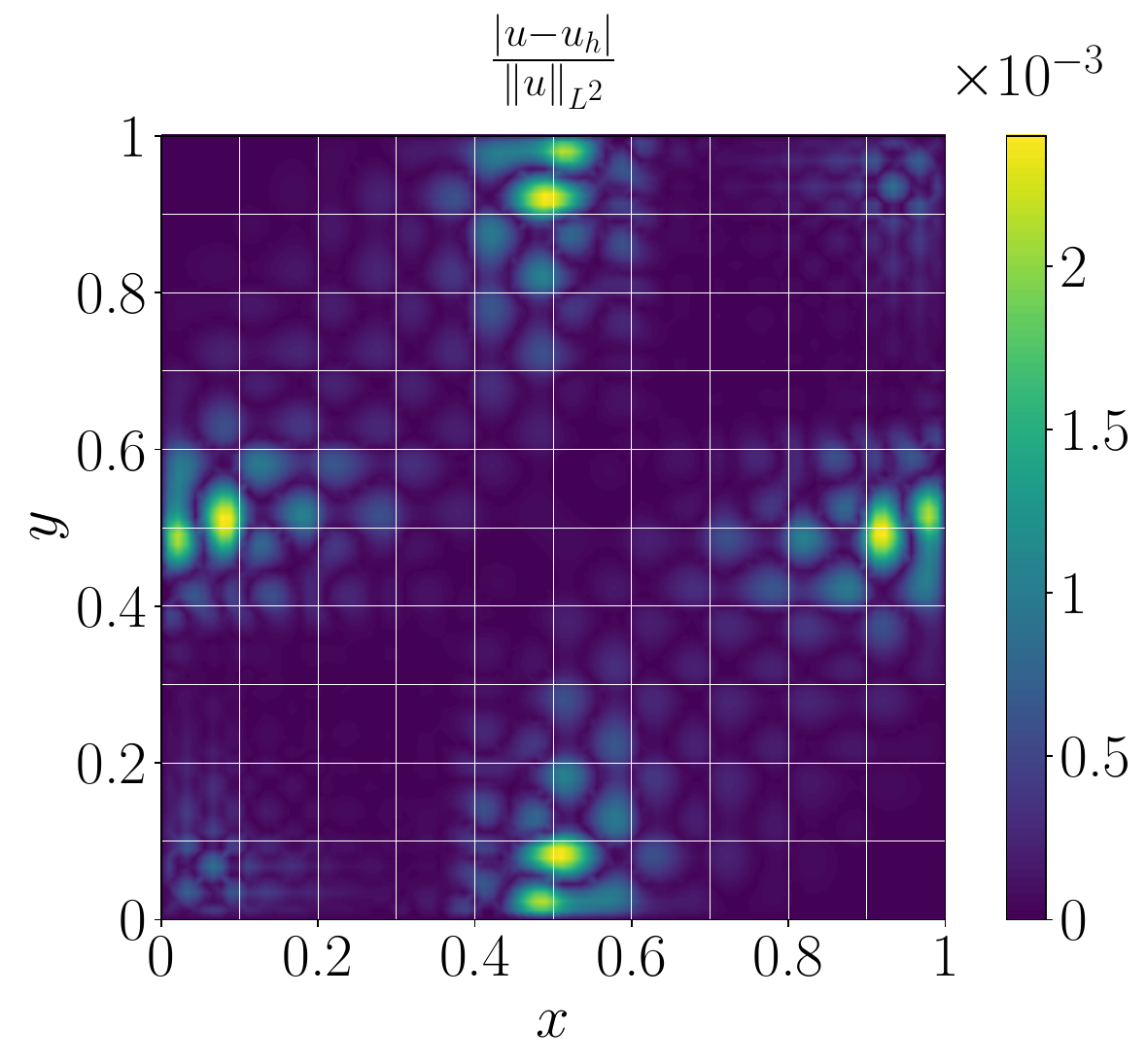}
        \caption{$p$-refinement (2 steps),\\ 454 DOFs, 100 Elements}
        \label{fig:loc_ref_f}
    \end{subfigure}

    \caption{Error distribution for various local refinement strategies.}
    \label{fig:loc_ref_validate}
\end{figure}

{\color{red}
Before employing local refinement in the subsequent studies, we first verify the correct behavior of the refinement implementation. To this end, we analyze the pointwise error of a symmetric manufactured solution with non-homogeneous boundary conditions on an initial $10 \times 10$ mesh with polynomial degree $p=2$. Fig.~\ref{fig:sdsf} shows both a surface and contour representation of the solution $u(x,y)$.

Fig.~\ref{fig:loc_ref_validate} presents the normalized pointwise error for the unrefined mesh and for locally refined meshes obtained using $h$-, $p$-, and $k$-refinement strategies. Refinement is applied in regions of elevated error, first within the domain $\Omega^1 = [0, 0.5]^2 \cup [0.5, 1]^2$ and subsequently within $\Omega^2 = [0, 0.4]^2 \cup [0.6, 1]^2$. Pure $p$-refinement increases the polynomial order while leaving the mesh topology unchanged.

After one $k$-refinement step, as well as after two successive $h$- or $p$-refinement steps, the error within the refined regions becomes smaller than in the remaining coarse parts of the domain. This confirms that the local refinement procedures act in the intended manner. The remaining error consequently concentrates in regions that are not refined. More rigorous convergence studies are presented in the following sections.

This example is not intended as a comprehensive study of the different refinement strategies, but rather as an implementation verification demonstrating that targeted refinement reduces the error in the expected regions. Since the manufactured solution is smooth, no singular behavior is present and the example is not designed to assess optimal refinement strategies.

}

\captionsetup[subfigure]{justification=centering}

\subsection{Effect of the Enhanced Shift Operator}
\label{sec:operator_comparison}

We compare the classical shift operator \(\mathbb{S}_{\mathbf{d}}^{p}\)
with the mixed‐derivative variant
\(\widetilde{\mathbb{S}}_{\mathbf{d}}^{p}\) introduced in
Sec.~\ref{ssec:enhanced_shift}.  
The test configuration is the ring geometry of
Fig.~\ref{fig:sbm-overview}: an annulus with outer radius
\(r_{\mathrm{ext}} = 0.47\) and inner radius \(r = 0.10\), centered at
\((0.5,0.5)\) and embedded in the square
\([0,1]^2\).
Throughout the numerical examples we employ the manufactured solution
\begin{equation}\label{analytical_sol}
  u(x,y) = x\,y\,\cosh(x)\,\sin(y),
\end{equation}
from which the forcing term and boundary conditions are derived.
Previous SBM studies~\cite{main2018shifted,atallah2021shifted,%
antonelli2024shifted}
chose to refine by halving the mesh size
\(\bigl(h_{\text{new}} = h_{\text{old}}/2\bigr)\),
thereby ensuring that the surrogate boundary did not move farther from the true boundary.
{\color{red}
To examine the method’s sensitivity to small changes in mesh size, and to
verify whether more degrees of freedom always increase accuracy, we adopt a step-by-step strategy: at each iteration, the number of equally sized elements is increased by one in each coordinate direction (e.g., $10\times10 \rightarrow 11\times11 \rightarrow 12\times12$...). This corresponds to rebuilding an equidistant knot vector with increased resolution on a fixed geometry. Consequently, the surrogate boundary may move either closer to or farther from the true boundary at each step, providing a fine-grained assessment of the convergence behaviour for the two shift operators.}

\begin{figure}[h]
    \begin{subfigure}[b]{0.45\textwidth}
        \centering
        \includegraphics[width=\linewidth,
                         trim={0.0cm 0.0cm 0.0cm 0.0cm}, clip]
                         {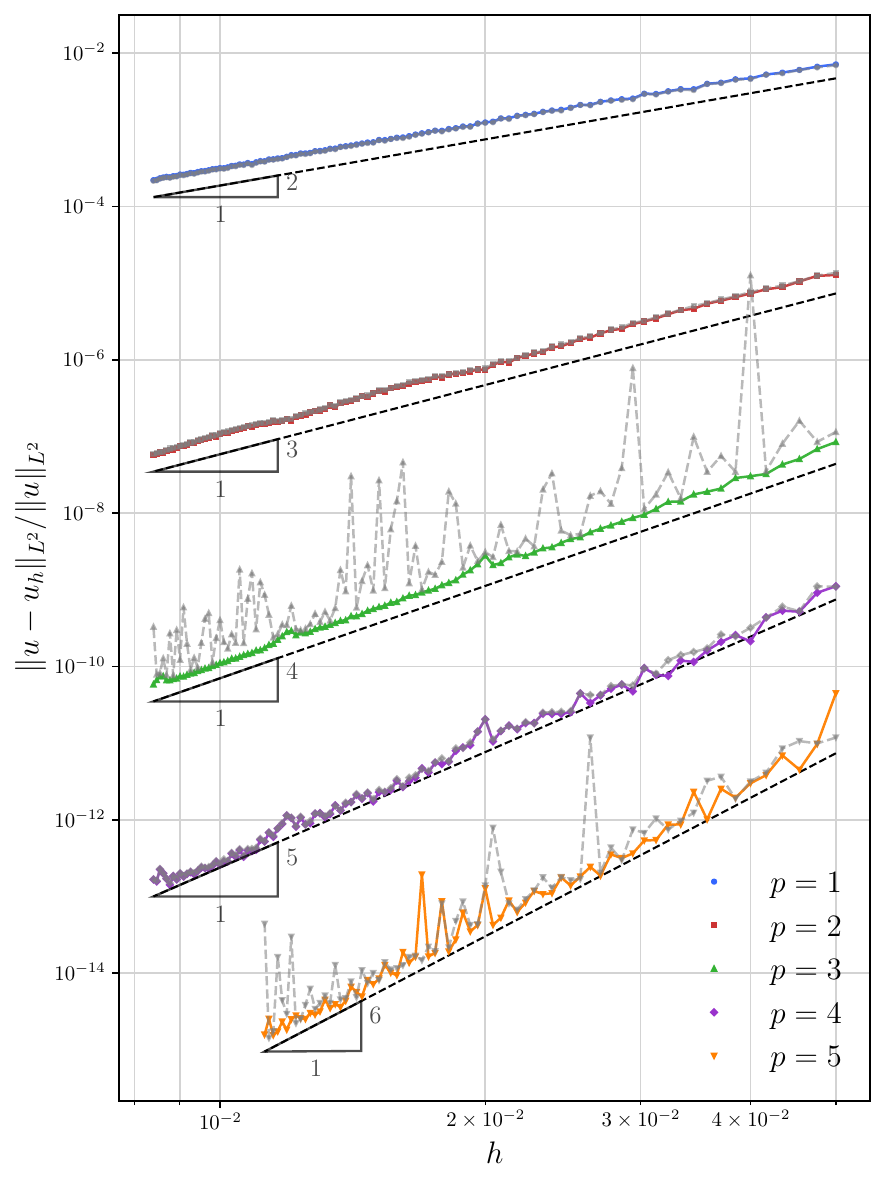}
        \caption{$L^2$ norm}
    \end{subfigure}
    \hspace{0.05\textwidth}
    \begin{subfigure}[b]{0.45\textwidth}
        \centering
        \includegraphics[width=\linewidth,
                         trim={0.0cm 0.0cm 0.0cm 0.0cm}, clip]
                         {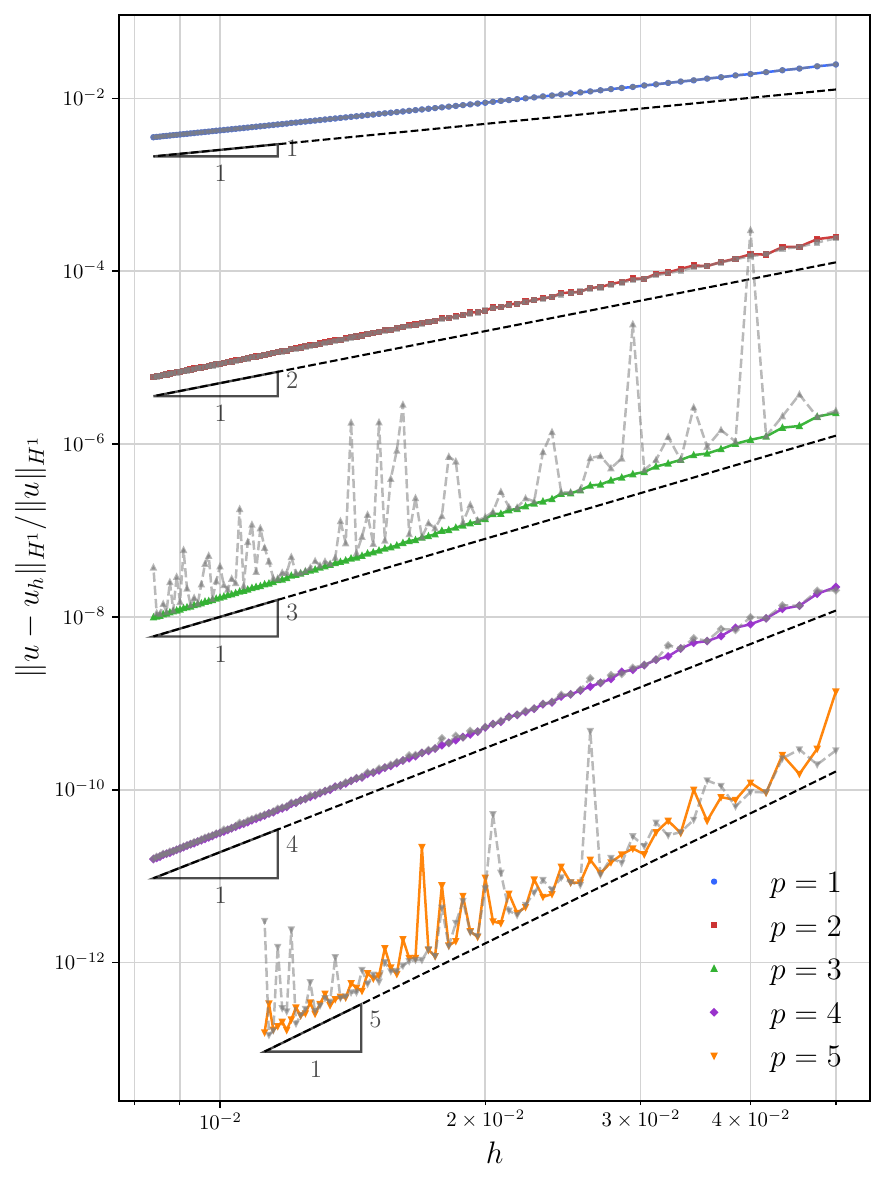}
        \caption{$H^1$ norm}
    \end{subfigure}
    \caption{ Effect of the enhanced shift operator
    \(\widetilde{\mathbb{S}}_{\mathbf{d}}^{p}\) on the convergence for the
    manufactured solution in Eq.~\ref{analytical_sol} with Dirichlet boundary conditions
    applied to the geometry shown in Fig.~\ref{fig:sbm-overview}.  
    Grey dashed lines depict the convergence obtained with the standard
    operator \(\mathbb{S}_{\mathbf{d}}^{p}\), while dashed black lines
    indicate the reference slope corresponding to the optimal rate.}
    \label{fig:shift operator Dirichlet}
\end{figure}

Fig.~\ref{fig:shift operator Dirichlet} shows the resulting errors for
Dirichlet boundary conditions on both boundaries.  
For \(p=1\) and \(p=2\) the two operators give almost identical, smoothly
convergent curves in both norms.  
For \(p=3\) the standard operator exhibits pronounced oscillations upon
refinement, whereas the enhanced operator removes almost all of them.  
For \(p=4\) and \(p=5\) the benefit is less obvious, although the
enhanced operator appears less sensitive to the occasional “lucky” or
“unlucky’’ placement of the surrogate boundary; the remaining peaks at
\(p=5\) are already close to machine precision.

\begin{figure}[h]
    \begin{subfigure}[b]{0.45\textwidth}
        \centering
        \includegraphics[width=\linewidth,
                         trim={0.0cm 0.0cm 0.0cm 0.0cm}, clip]
                         {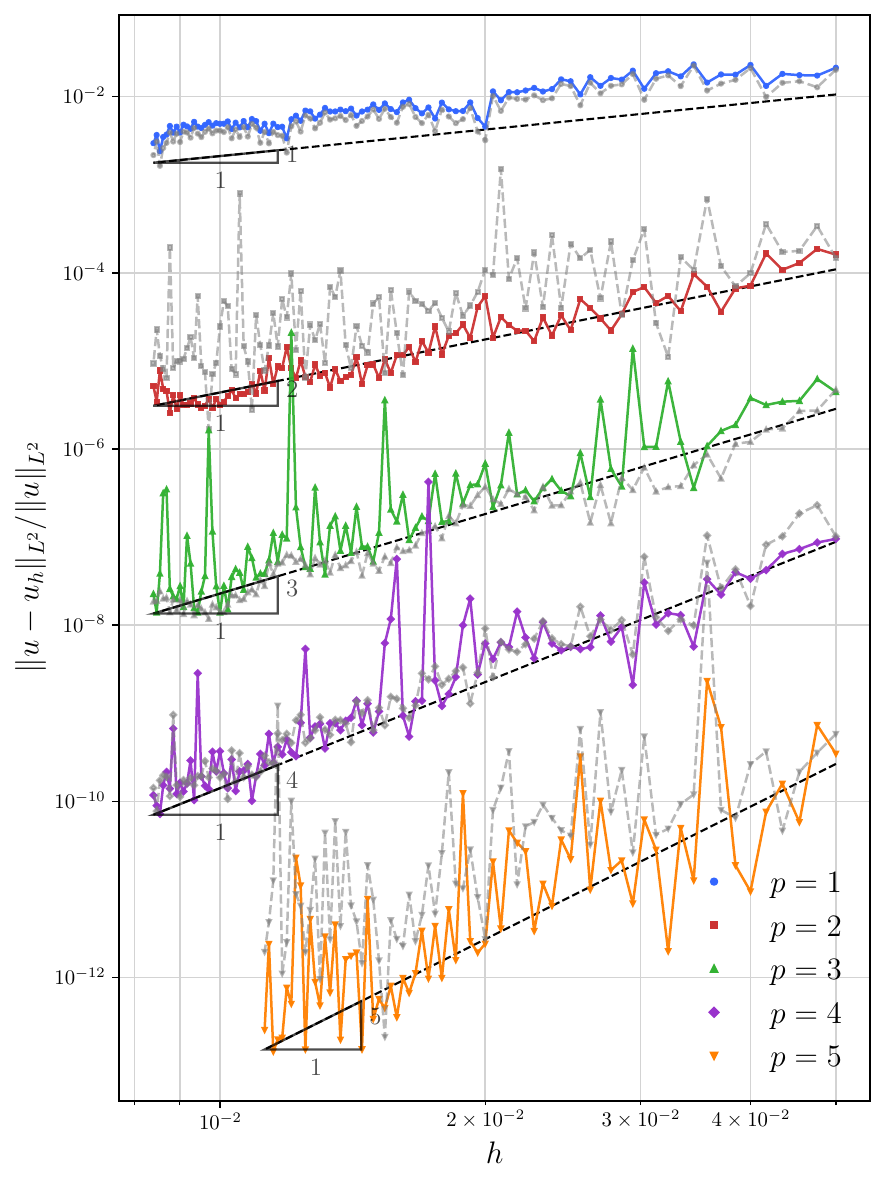}
        \caption{$L^2$ norm}
    \end{subfigure}
    \hspace{0.05\textwidth}
    \begin{subfigure}[b]{0.45\textwidth}
        \centering
        \includegraphics[width=\linewidth,
                         trim={0.0cm 0.0cm 0.0cm 0.0cm}, clip]
                         {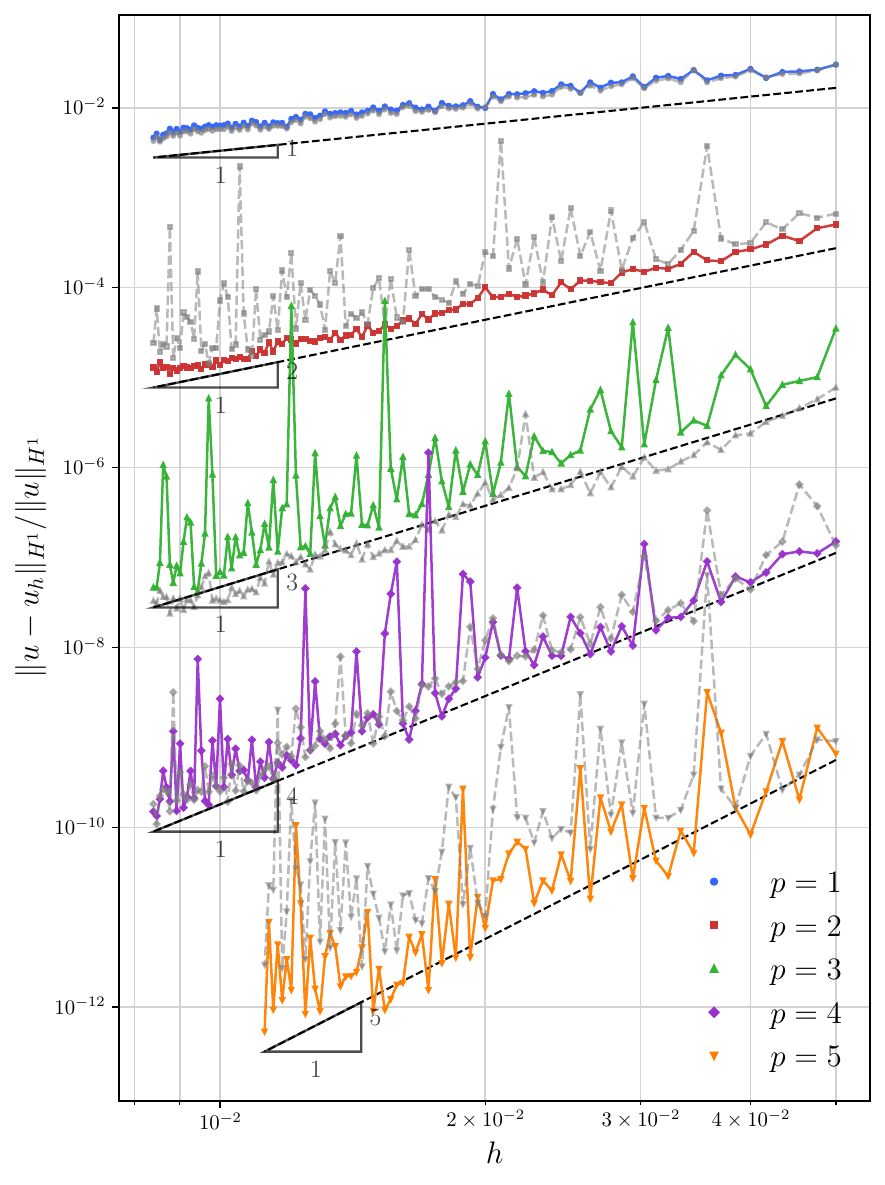}
        \caption{$H^1$ norm}
    \end{subfigure}
     \caption{Effect of the enhanced shift operator
\(\widetilde{\mathbb{S}}_{\mathbf{d}}^{p}\) on convergence for the manufactured
solution in Eq.~\ref{analytical_sol}, with Neumann conditions on the
exterior boundary and Dirichlet conditions on the interior boundary of
the geometry shown in Fig.~\ref{fig:sbm-overview}.  
Grey dashed lines depict the convergence obtained with the standard
operator \(\mathbb{S}_{\mathbf{d}}^{p}\), while dashed black lines
represent the reference slopes.}
\label{fig:shift operator Neumann}
\end{figure}

Fig.~\ref{fig:shift operator Neumann} presents the corresponding study
with Neumann conditions on the exterior boundary and Dirichlet
conditions on the interior one. The Dirichlet data are always treated with the enhanced operator, as it
proved to never be worse than the standard version. Here the advantage is less clear-cut: for \(p=1\) the two curves overlap
almost exactly; for \(p=2\) the enhanced operator yields smoother
convergence, but for \(p=3\) the trend reverses, with the standard
operator producing the smoother lines.
For \(p=4\) the enhanced operator introduces slightly more oscillations,
while for \(p=5\) it introduces slightly fewer.

We repeated these tests with different manufactured solutions and
geometries and obtained the same qualitative behaviour.  
For Dirichlet boundary conditions the enhanced operator is consistently
beneficial, particularly for \(p=3\), where the classical operator
reliably shows stronger oscillations.  
For Neumann data it improves the behaviour for \(p=2\). For higher degrees
it can sometimes introduce additional oscillations. 
We finally remark that the convergence in the $H^1$ norm for Neumann boundary conditions has the optimal convergence rate.
{\color{red} These observations are purely empirical and do not constitute a theoretical robustness result. 
They suggest that the enhanced shift operator may reduce step-by-step convergence oscillations in certain SBM regimes.
}
Accordingly, the remainder of this paper adopts the following rule:

\begin{itemize}
  \item For Neumann boundary conditions with boundary‐supporting basis
        functions of degree \(p\le 2\), \emph{and for all Dirichlet
        conditions regardless of degree}, we use the enhanced operator
        \(\widetilde{\mathbb{S}}_{\mathbf{d}}^{p}\).
  \item For shifted Neumann conditions with boundary basis functions of
        degree \(p>2\), we retain the classical operator
        \(\mathbb{S}_{\mathbf{d}}^{p}\).
\end{itemize}

\subsection{Local refinement with Dirichlet boundary conditions} 
Let the domain be as described in Fig.~\ref{fig:sbm-thb}, specifically, a square $[0,1]^2$ with a circular hole of radius $r=0.15$ centered at $(0.5,0.5)$. Dirichlet boundary conditions are imposed on the exterior boundary by a standard penalty-free Nitsche approach, while a surrogate boundary is constructed on the interior boundary to impose the Dirichlet data via the SBM. {\color{red}For all convergence studies presented henceforth, the characteristic mesh size is defined as the maximum edge length among all quadrilateral elements. As a reference configuration, we employ a body-fitted setup in which the background mesh used for the SBM is retained, but the boundary conditions are imposed directly on the surrogate boundary rather than being shifted to the true boundary. This eliminates the Taylor-expansion consistency error and allows the influence of the boundary-shift operator to be isolated.

}

\begin{figure}[htbp]
    \centering

    \begin{subfigure}[t]{0.3\textwidth}
        \centering
        \includegraphics[width=\textwidth]{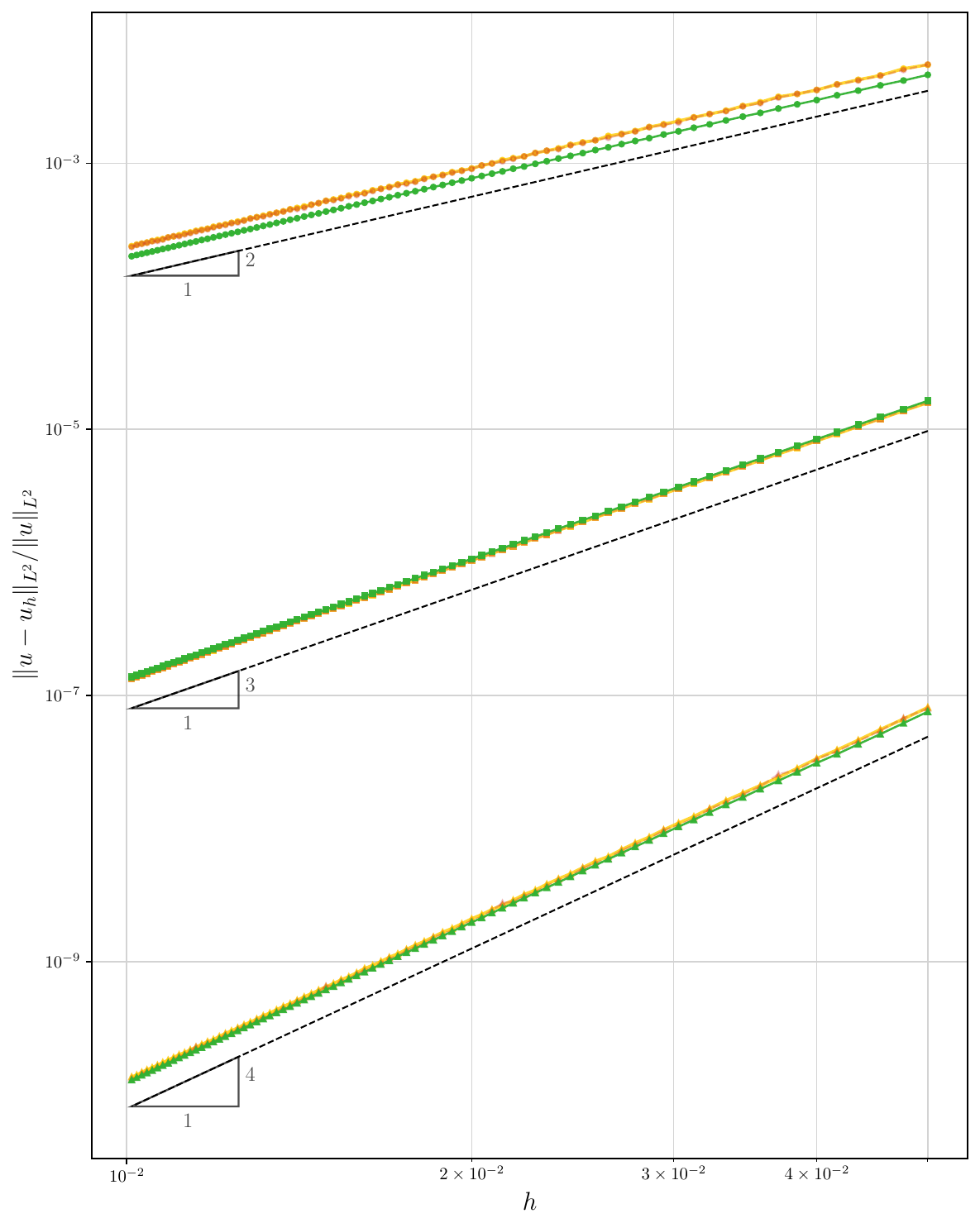}
        \caption{$h$-refinement}
        \label{fig:ex11_a}
    \end{subfigure}
    \hfill
    \begin{subfigure}[t]{0.3\textwidth}
        \centering
        \includegraphics[width=\textwidth]{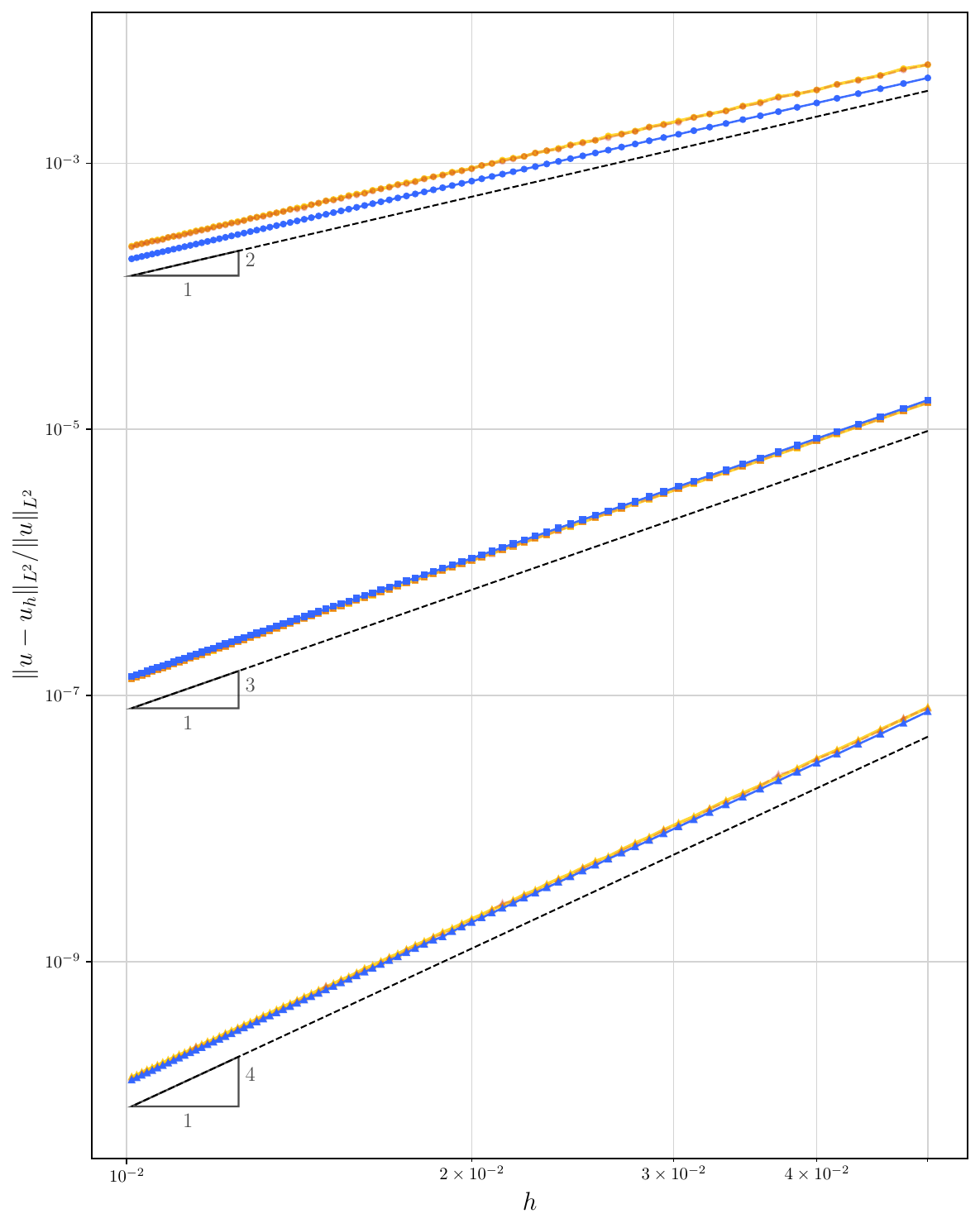}
        \caption{$p$-refinement}
        \label{fig:ex11_b}
    \end{subfigure}
    \hfill
    \begin{subfigure}[t]{0.3\textwidth}
        \centering
        \includegraphics[width=\textwidth]{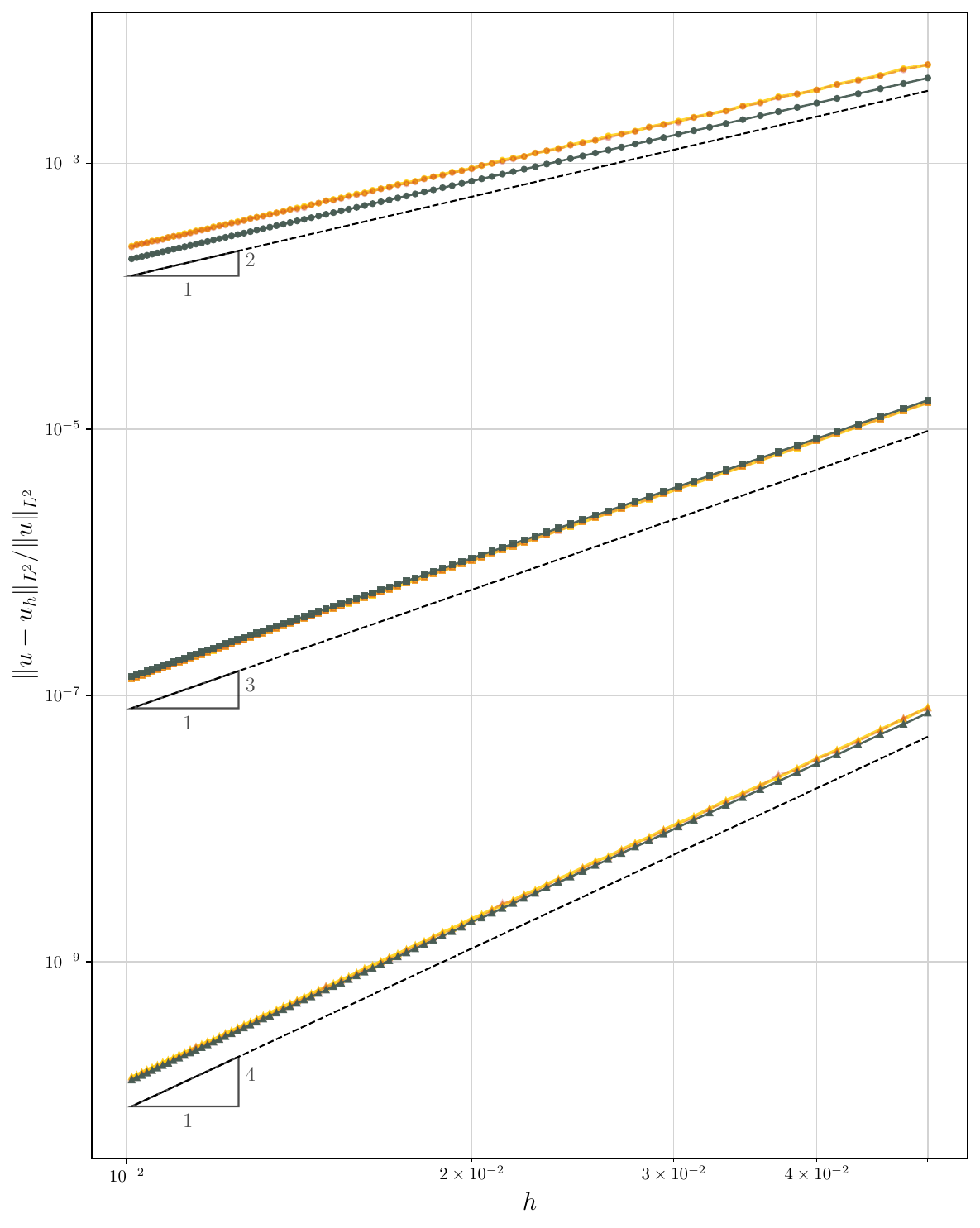}
        \caption{$k$-refinement}
        \label{fig:ex11_c}
    \end{subfigure}

    \vspace{1em}

    \begin{subfigure}[t]{0.3\textwidth}
        \centering
        \includegraphics[width=\textwidth]{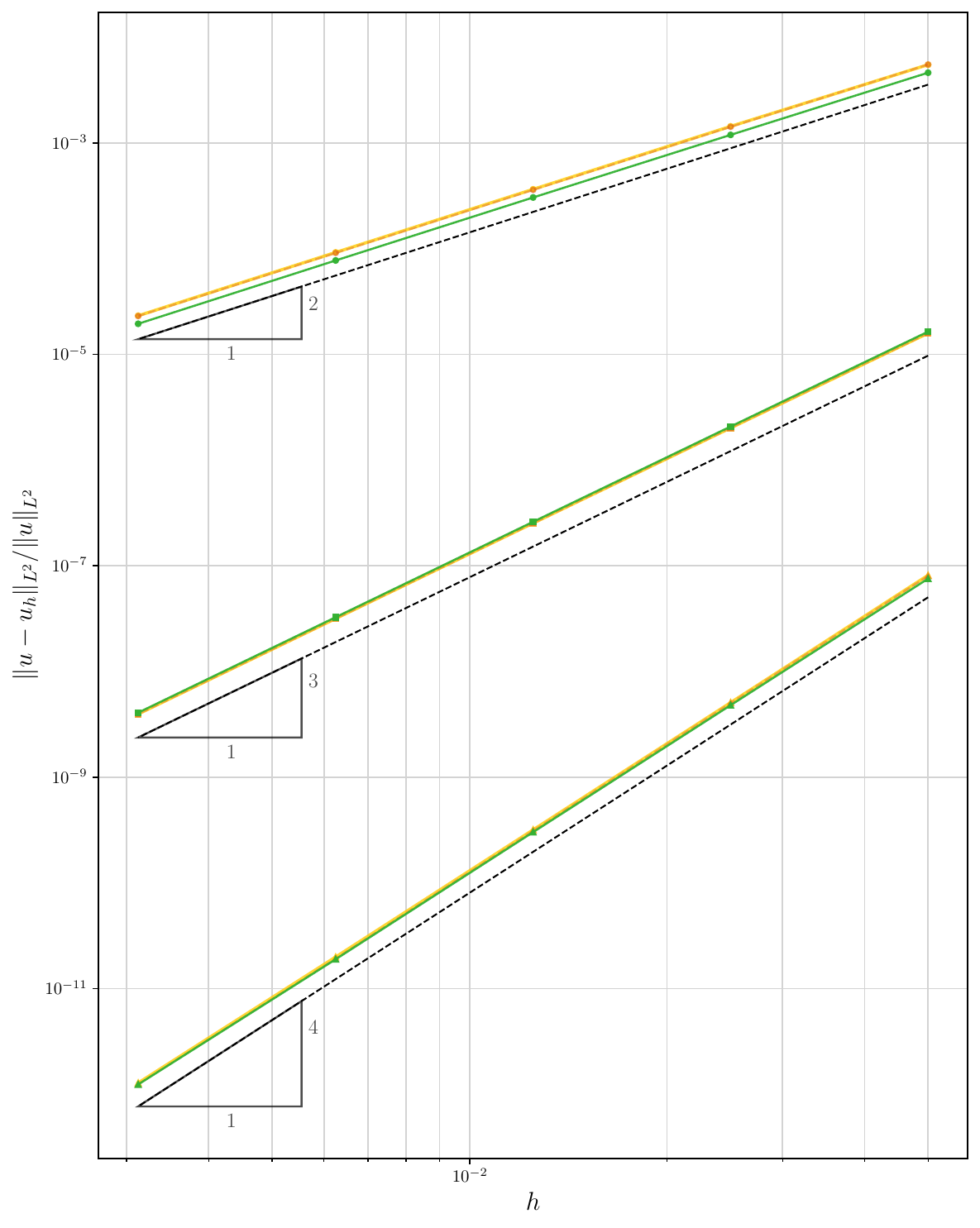}
        \caption{$h$-refinement}
        \label{fig:ex11_d}
    \end{subfigure}
    \hfill
    \begin{subfigure}[t]{0.3\textwidth}
        \centering
        \includegraphics[width=\textwidth]{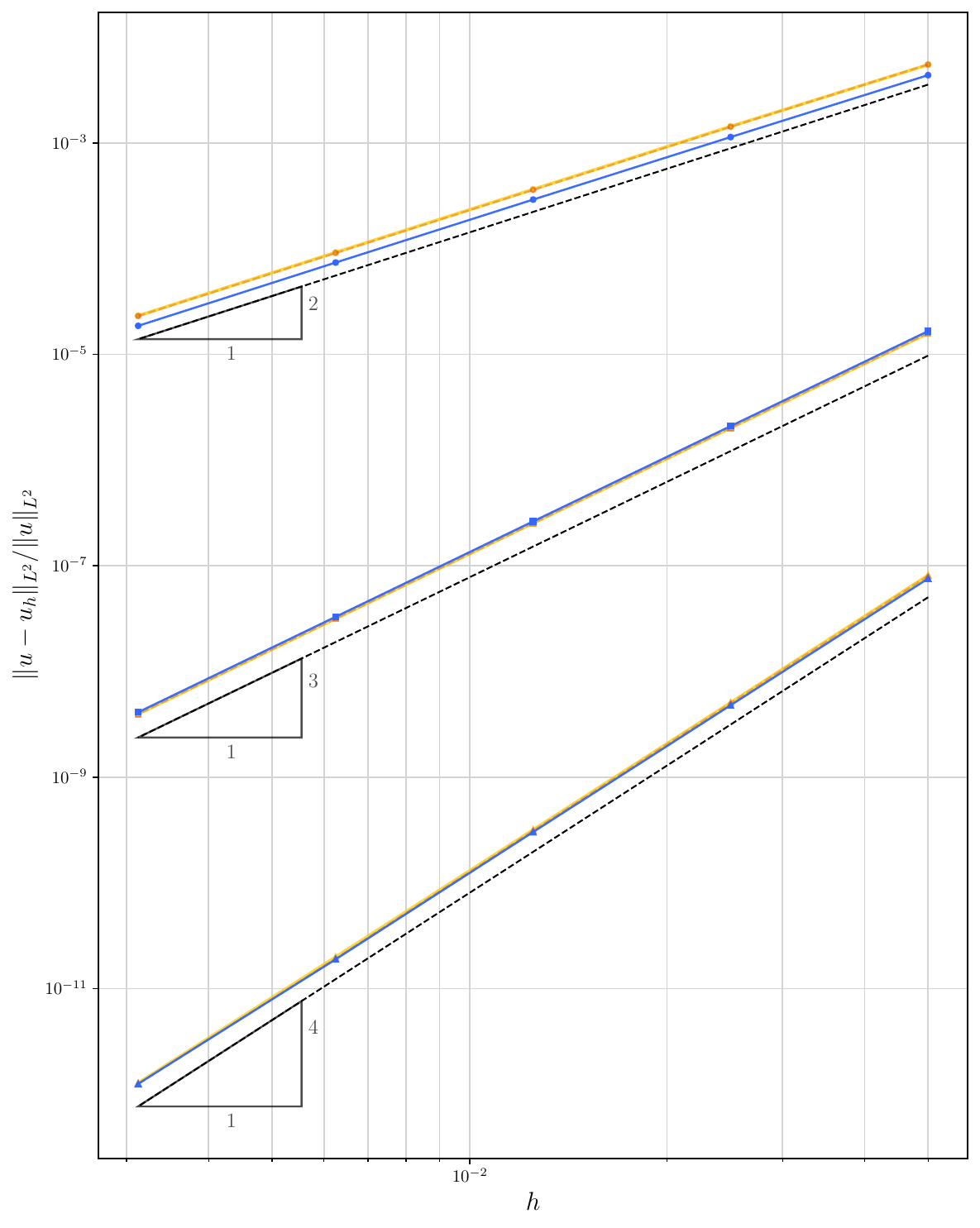}
        \caption{$p$-refinement}
        \label{fig:ex11_e}
    \end{subfigure}
    \hfill
    \begin{subfigure}[t]{0.3\textwidth}
        \centering
        \includegraphics[width=\textwidth]{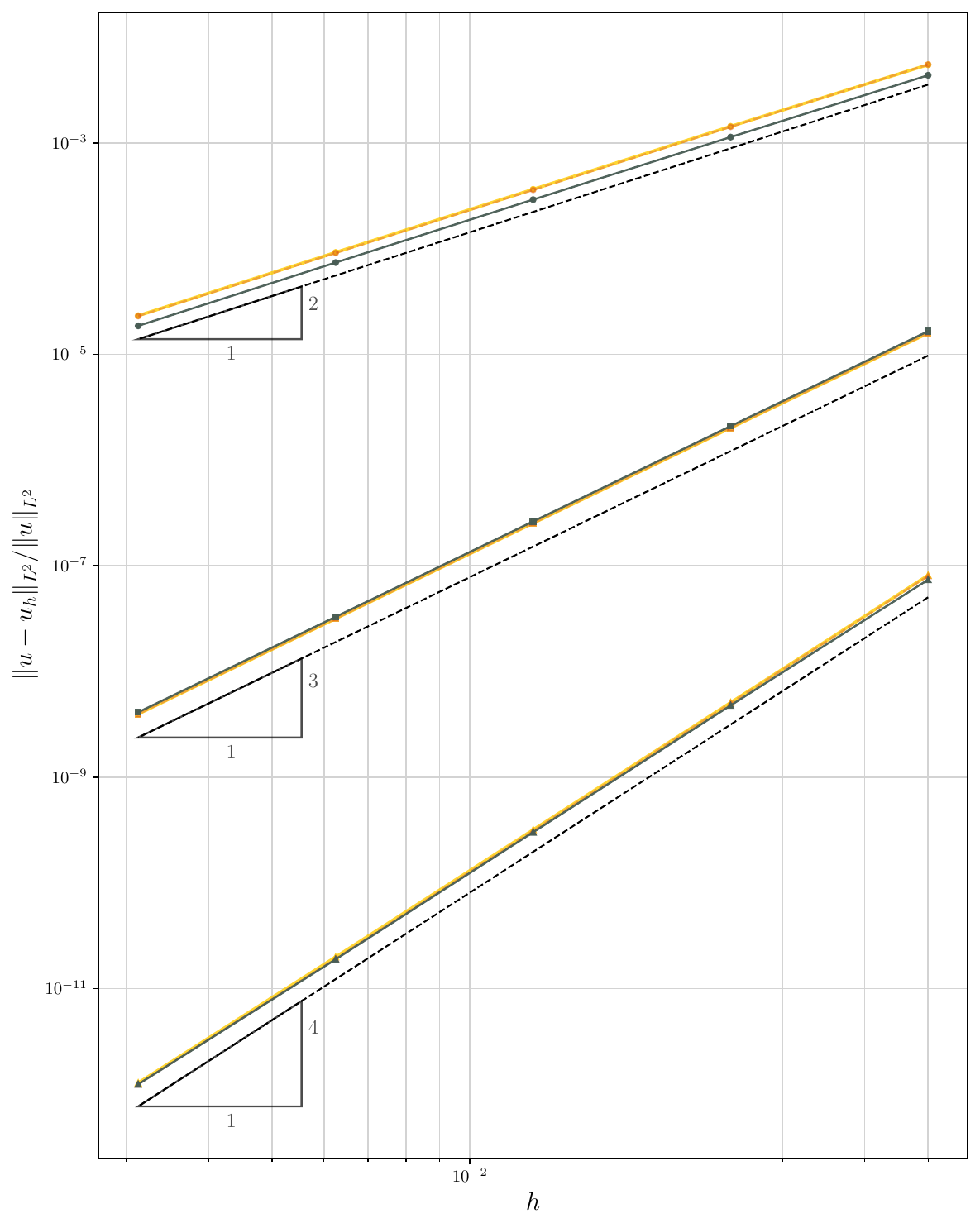}
        \caption{$k$-refinement}
        \label{fig:ex11_f}
    \end{subfigure}

    \caption{Convergence study for the geometry in Fig.~\ref{fig:sbm-thb} with different local refinement strategies and only Dirichlet boundary conditions. 
    First row: step-by-step convergence from $20\times20$ to $100\times100$ knot spans; second row: standard halved convergence from $20\times20$ to $320\times320$ knot spans.
    Legend: yellow: body-fitted; dashed red: non-refined SBM; dashed black: reference line for the optimal convergence; green: $h$-refinement; blue: $p$-refinement; black: $k$-refinement.  Circles: $p = 1$; squares: $p=2$; triangles $p=3$.}
    \label{fig:ex11}
\end{figure}

Fig.~\ref{fig:ex11} compares the convergence obtained with local $h$-, $p$-, and $k$-refinements {\color{red}after a single local refinement step applied to the surrogate boundary}. {\color{red}As expected for Dirichlet boundary conditions,} none of the refinement strategies improves the convergence slope, which is already {\color{red}of order $p+1$ in $L^2$} and thus optimal with the standard SBM; the additional DOFs introduced only slightly improve the convergence constant. None of the local refinement strategies introduces oscillations in the convergence, and the step-by-step curves are essentially straight lines, as in the halving study.

\begin{figure}[htbp]
    \centering

    \begin{subfigure}[t]{0.3\textwidth}
        \centering
        \includegraphics[width=\textwidth]{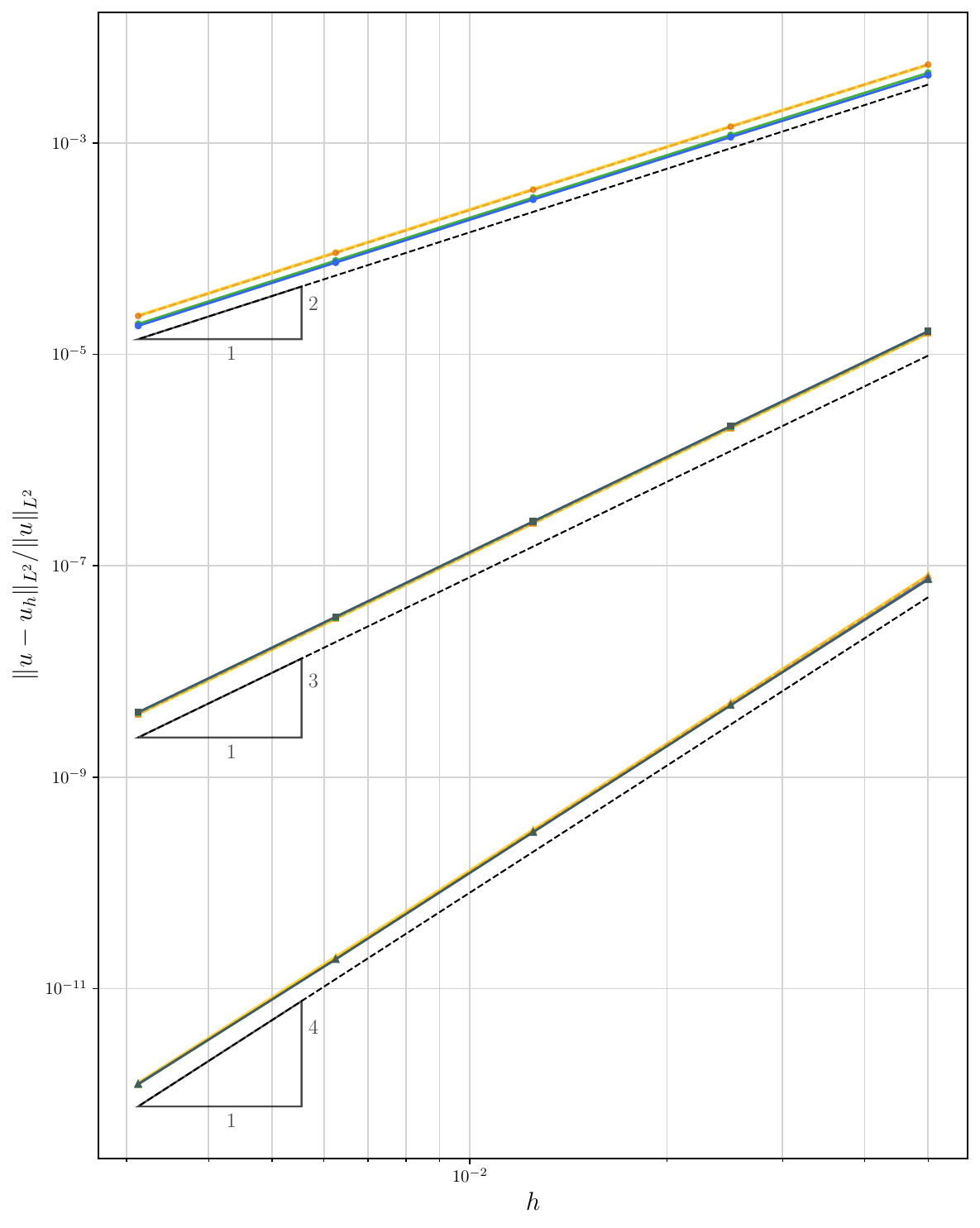}
        \caption{Convergence comparison in the $L^2$ norm.}
        \label{fig:ex11_comp_a}
    \end{subfigure}
    \hfill
    \begin{subfigure}[t]{0.3\textwidth}
        \centering
        \includegraphics[width=\textwidth]{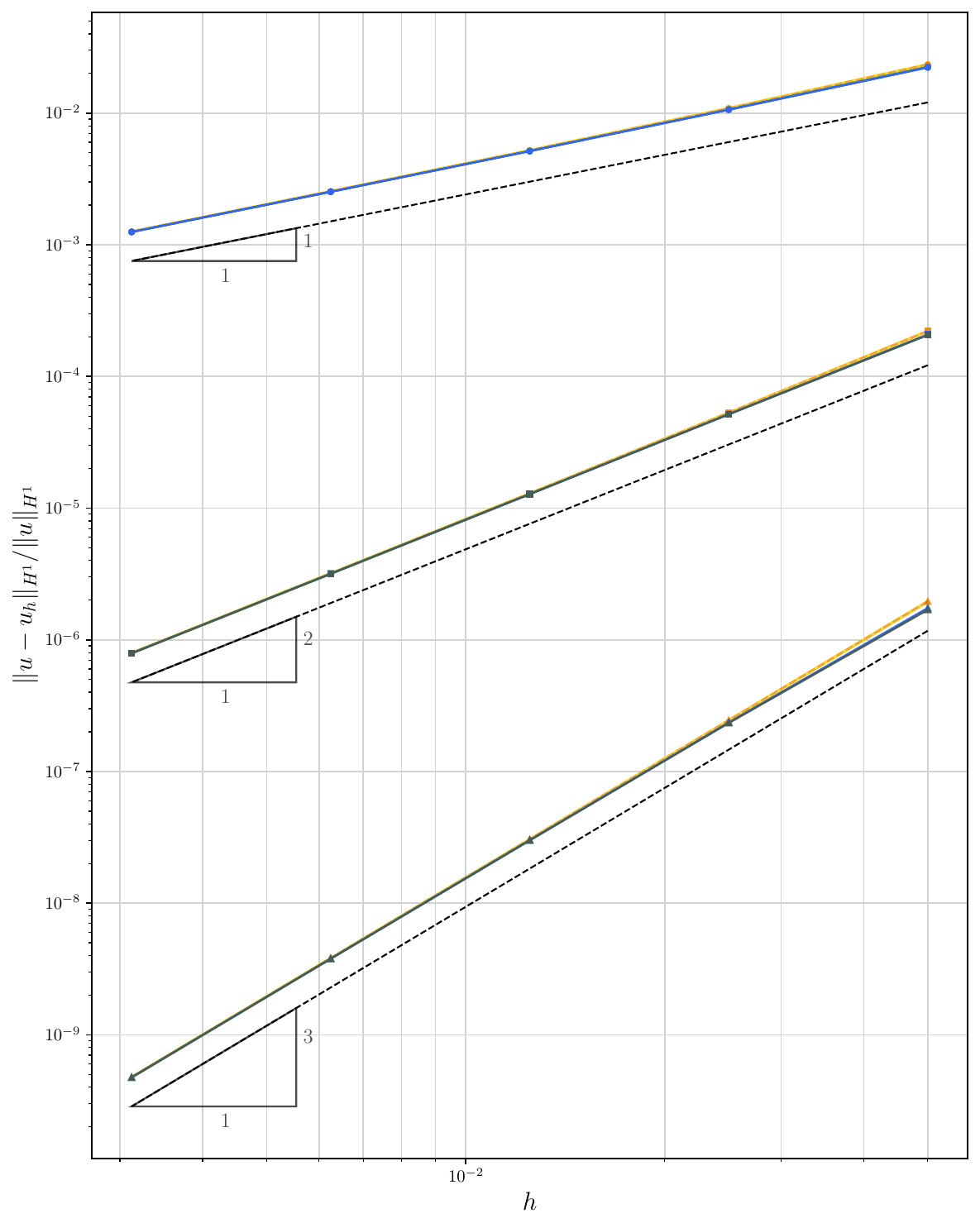}
        \caption{Convergence comparison in the $H^1$ norm.}
        \label{fig:ex11_comp_b}
    \end{subfigure}
    \hfill
    \begin{subfigure}[t]{0.3\textwidth}
        \centering
        \includegraphics[width=\textwidth]{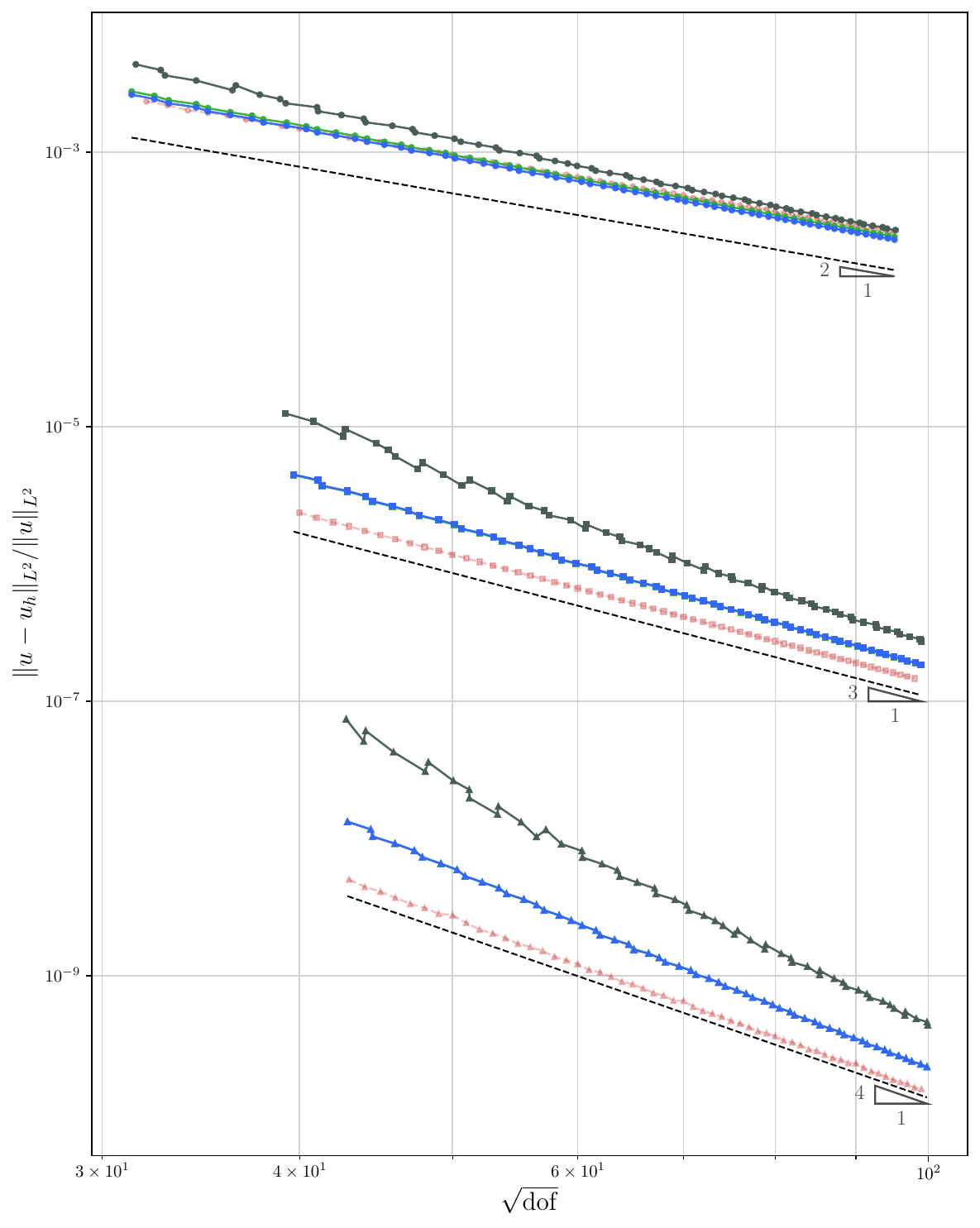}
        \caption{Convergence comparison in the $L^2$ norm against the DOFs.}
        \label{fig:ex11_comp_c}
    \end{subfigure}

    \caption{Head-to-head convergence comparison, from $20\times20$ to $320\times320$ knot spans, for the $h$-, $p$-, and $k$-refinement schemes and the non-refined SBM, applied to the geometry presented in Fig.~\ref{fig:sbm-thb} with only Dirichlet BCs. The color and marker legends are the same as in Fig.~\ref{fig:ex11}.}
    \label{fig:ex_11_comparison}
\end{figure}

Finally, Fig.~\ref{fig:ex_11_comparison} compares the standard SBM with the three refinement strategies. No clear difference can be observed in terms of the $L^2$ or $H^1$ convergence, but Fig.~\ref{fig:ex11_comp_c} clearly highlights the higher cost in terms of DOFs for the locally refined cases. 
{\color{red}These results indicate that, for immersed Dirichlet boundaries, additional local refinement does not provide a meaningful advantage over the standard SBM.}
{\color{red}In this configuration, the boundary-consistency error is not dominant and the global error is already comparable to the interior discretization error. Consequently, refining only in the vicinity of the boundary does not reduce the overall error norm.
}
Consequently, {\color{red}no additional refinement is applied to immersed Dirichlet boundaries in the following examples.}



\subsection{Local refinement with Neumann boundary conditions}
\begin{figure}[htbp]
    \centering

    \begin{subfigure}[t]{0.3\textwidth}
        \centering
        \includegraphics[width=\textwidth]{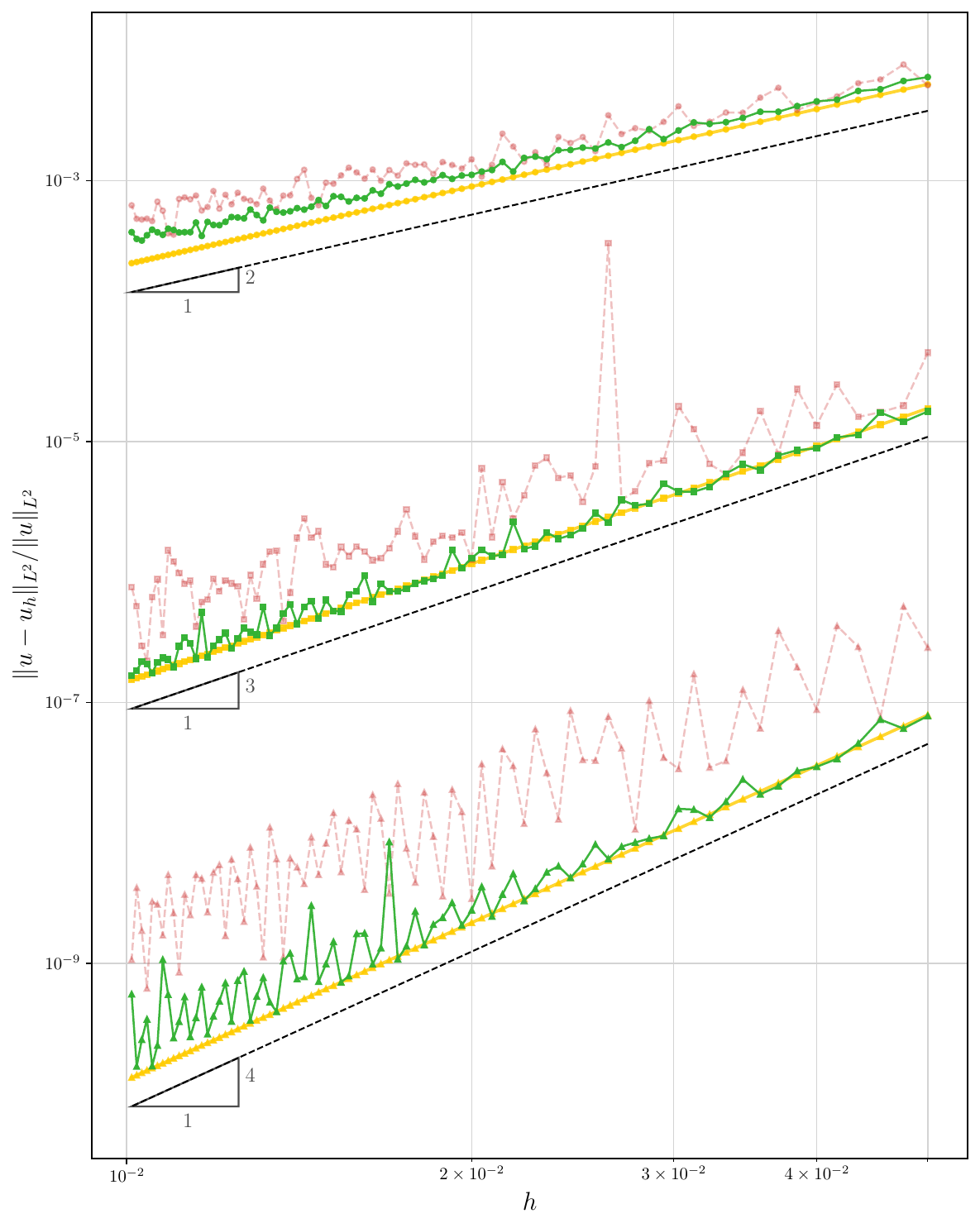}
        \caption{$h$-refinement}
        \label{fig:ex12_a}
    \end{subfigure}
    \hfill
    \begin{subfigure}[t]{0.3\textwidth}
        \centering
        \includegraphics[width=\textwidth]{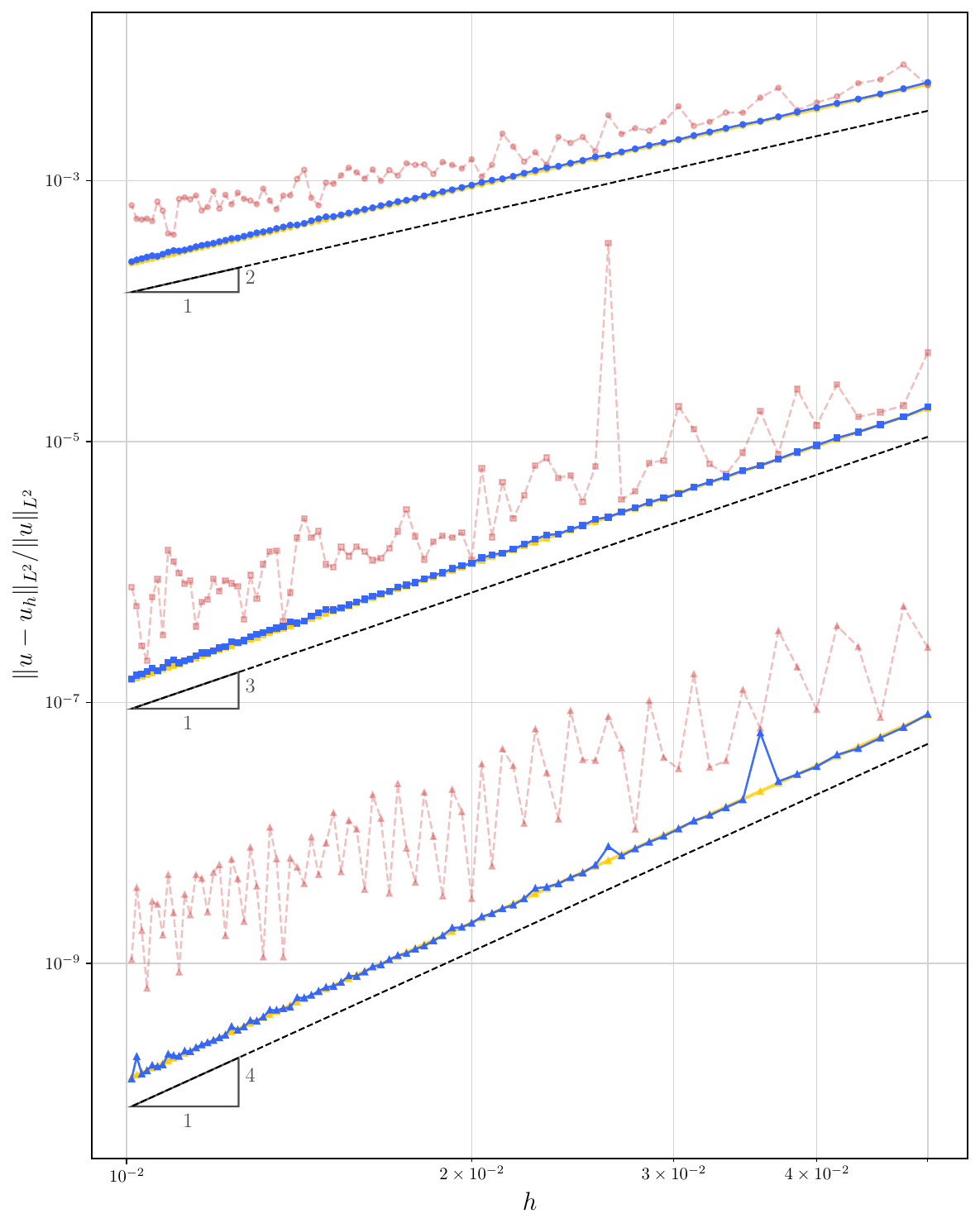}
        \caption{$p$-refinement}
        \label{fig:ex12_b}
    \end{subfigure}
    \hfill
    \begin{subfigure}[t]{0.3\textwidth}
        \centering
        \includegraphics[width=\textwidth]{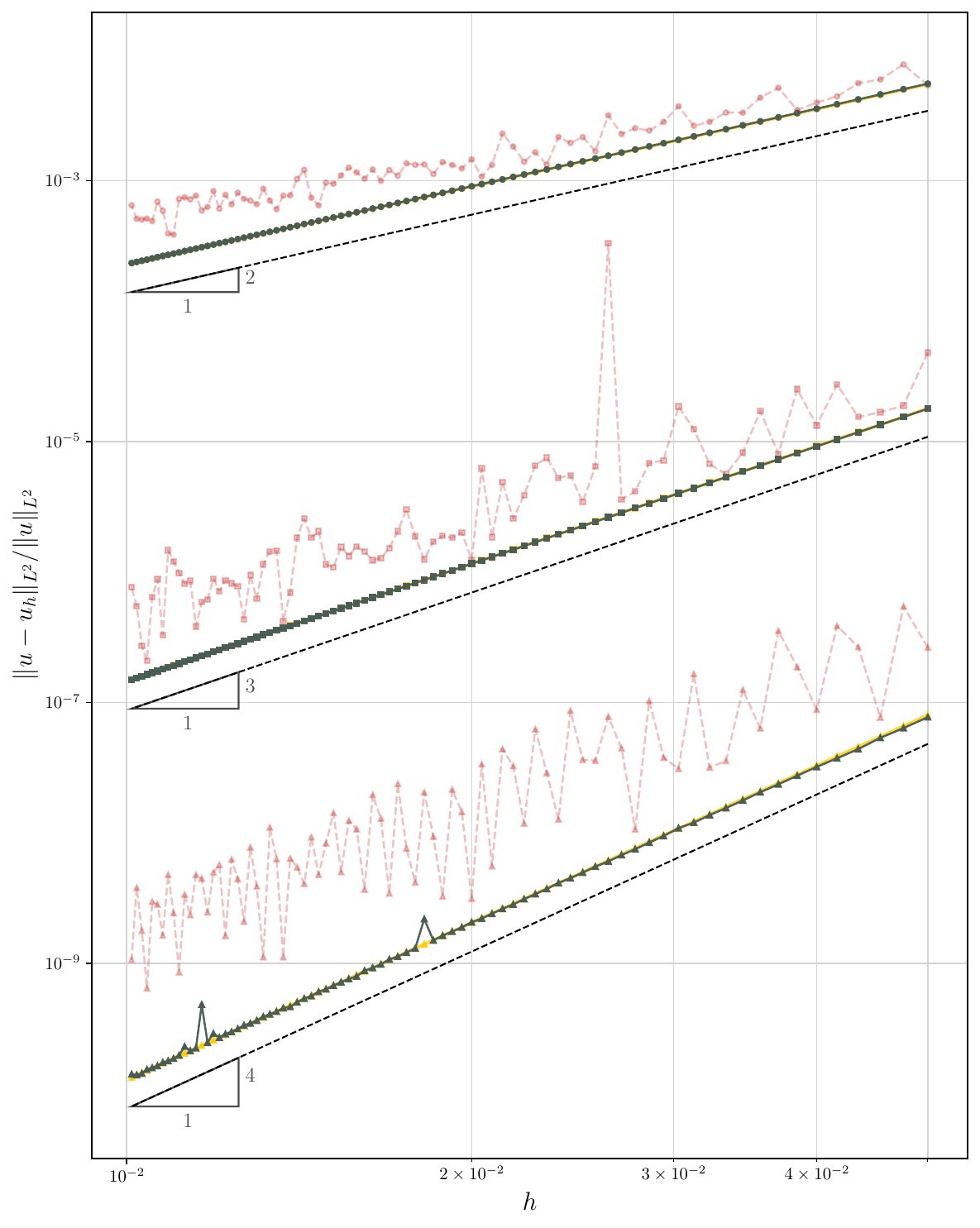}
        \caption{$k$-refinement}
        \label{fig:ex12_c}
    \end{subfigure}

    \vspace{1em}

    \begin{subfigure}[t]{0.3\textwidth}
        \centering
        \includegraphics[width=\textwidth]{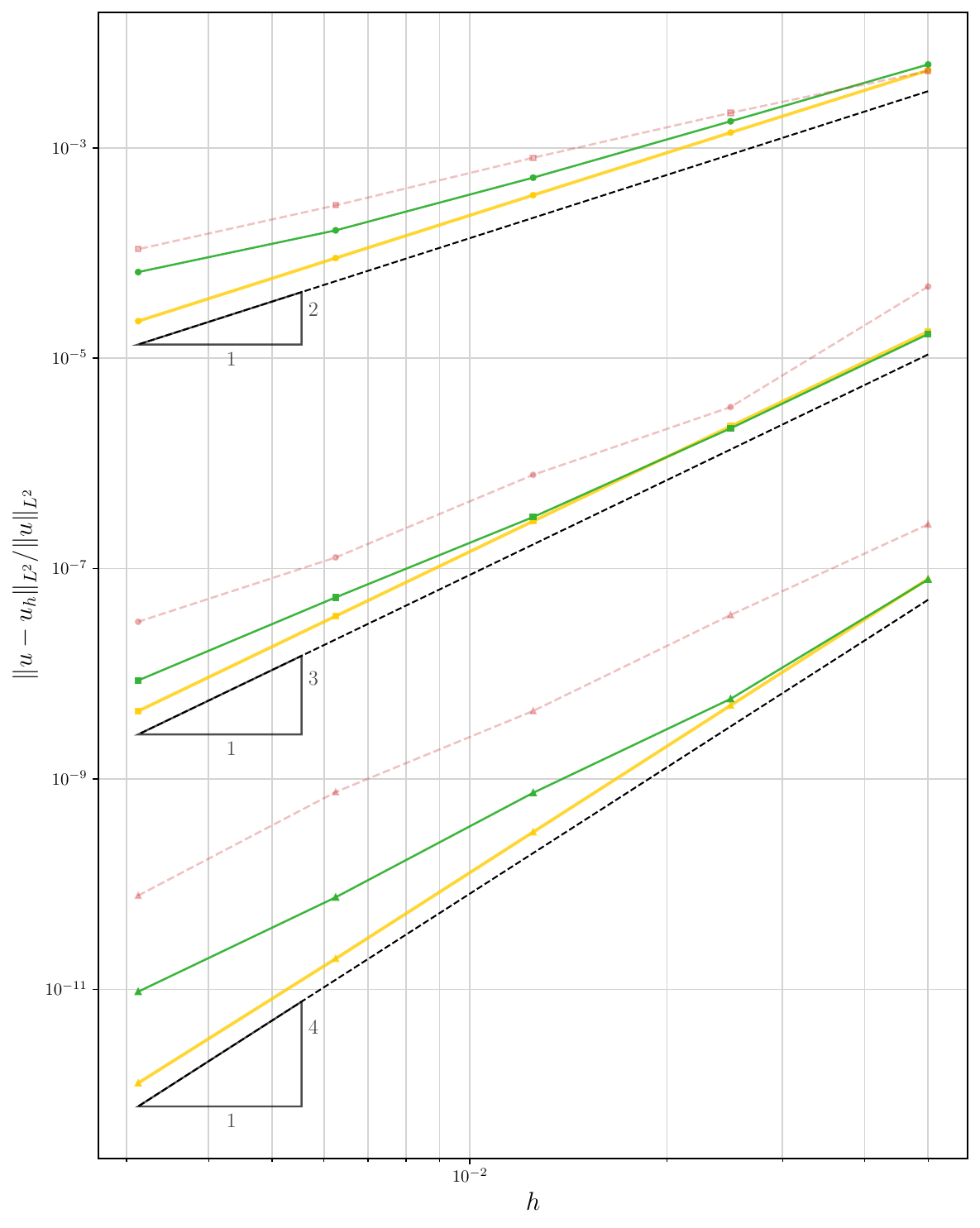}
        \caption{$h$-refinement}
        \label{fig:ex12_d}
    \end{subfigure}
    \hfill
    \begin{subfigure}[t]{0.3\textwidth}
        \centering
        \includegraphics[width=\textwidth]{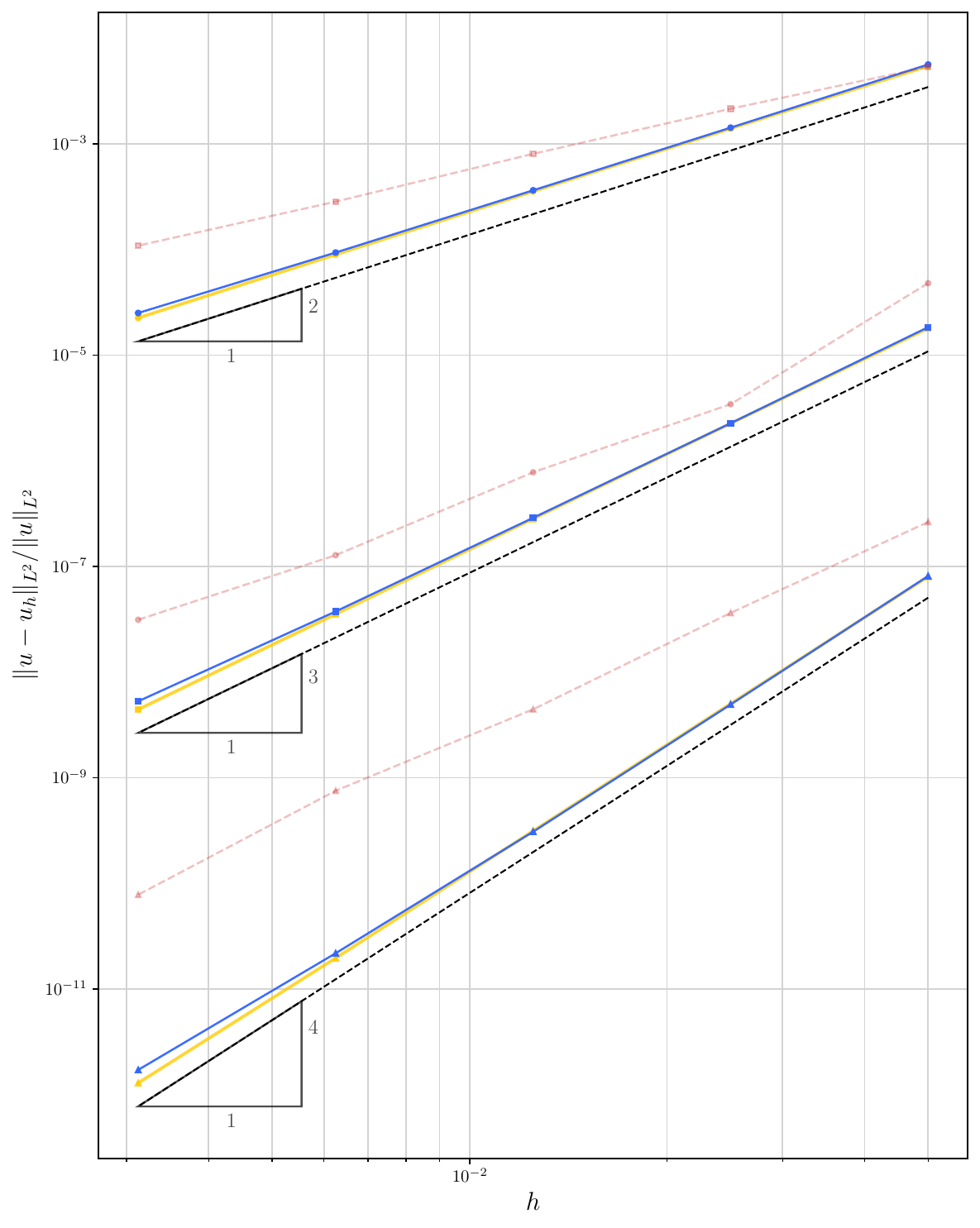}
        \caption{$p$-refinement}
        \label{fig:ex12_e}
    \end{subfigure}
    \hfill
    \begin{subfigure}[t]{0.3\textwidth}
        \centering
        \includegraphics[width=\textwidth]{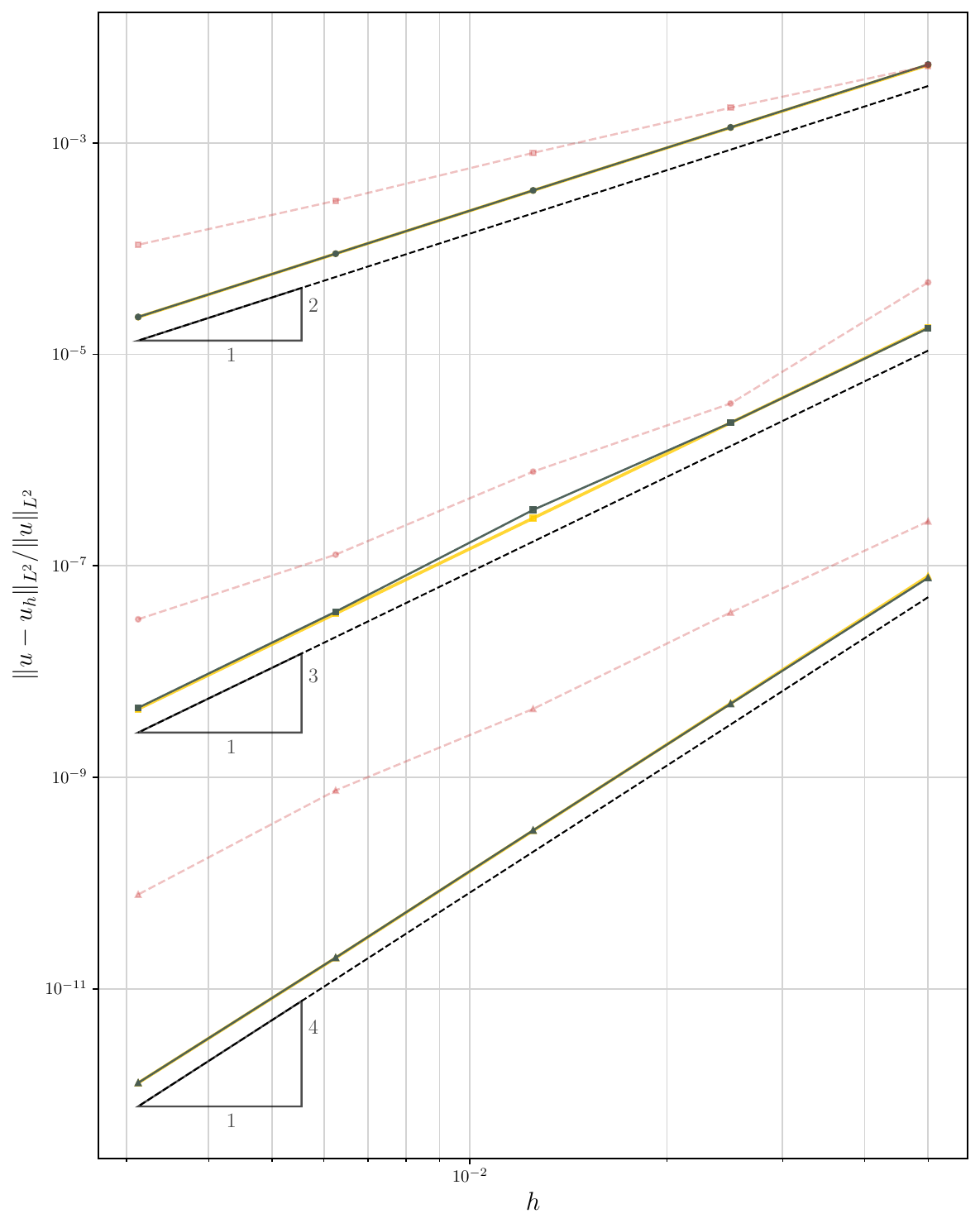}
        \caption{$k$-refinement}
        \label{fig:ex12_f}
    \end{subfigure}

    \caption{Convergence study for the geometry \ref{fig:sbm-thb} with different local refinement strategies and immersed Neumann BCs. First row: step-by-step convergence from $20\times20$ to $100\times100$ knot spans; second row: standard halved convergence from $20\times20$ to $320\times320$ knot spans. 
    Legend: yellow: body-fitted; dashed red: non-refined SBM; dashed black: reference line for the optimal convergence; green: $h$-refinement; blue: $p$-refinement; black: $k$-refinement. Circles: $p = 1$; squares: $p=2$; triangles $p=3$.}
    \label{fig:ex12}
\end{figure}

Let us now consider the same geometry, but with Neumann conditions applied on the immersed circular boundary of Fig.~\ref{fig:sbm-thb} instead of Dirichlet conditions.
The convergence results for the three refinement schemes are shown in Fig.~\ref{fig:ex12}.
The color legend is the same as in the previous example. In this case we clearly observe the oscillatory and suboptimal convergence for the standard SBM (red lines).
{\color{red} In this configuration, the error is predominantly governed by the SBM boundary approximation, so that local refinement directly targets the dominant error source.
}
All local refinement strategies substantially improve the accuracy, with some differences.
The $h$-refinement (green lines) recovers the accuracy of the body-fitted just on coarse meshes and still exhibits minor oscillations as shown in Fig.~\ref{fig:ex12_a}.
From the halving study in Fig.~\ref{fig:ex12_d}, carried out up to the finest meshes considered, we clearly see the convergence curve detach from the optimal body-fitted reference.

The $p$-refinement (blue lines) in Fig.~\ref{fig:ex12_b} yields better results: it remains very close to the reference line for the optimal convergence even on very fine meshes, showing noticeable oscillations only for cubic splines.
Still, from the halving study for $p=3$ in Fig.~\ref{fig:ex12_e}, one can observe an incipient detachment from the optimal behaviour at the last refinement step, indicating that $p$-refinement does not eliminate the suboptimal convergence, but rather postpones the departure to very fine meshes.

Finally, the $k$-refinement in Fig.~\ref{fig:ex12_c} exhibits the best overall performance among the three techniques: only mild oscillations appear in the step-by-step study, and no detachment from the optimal rate is visible even for the finest meshes in the halving study in Fig.~\ref{fig:ex12_f}.

\begin{figure}[htbp]
    \centering

    \begin{subfigure}[t]{0.3\textwidth}
        \centering
        \includegraphics[width=\textwidth]{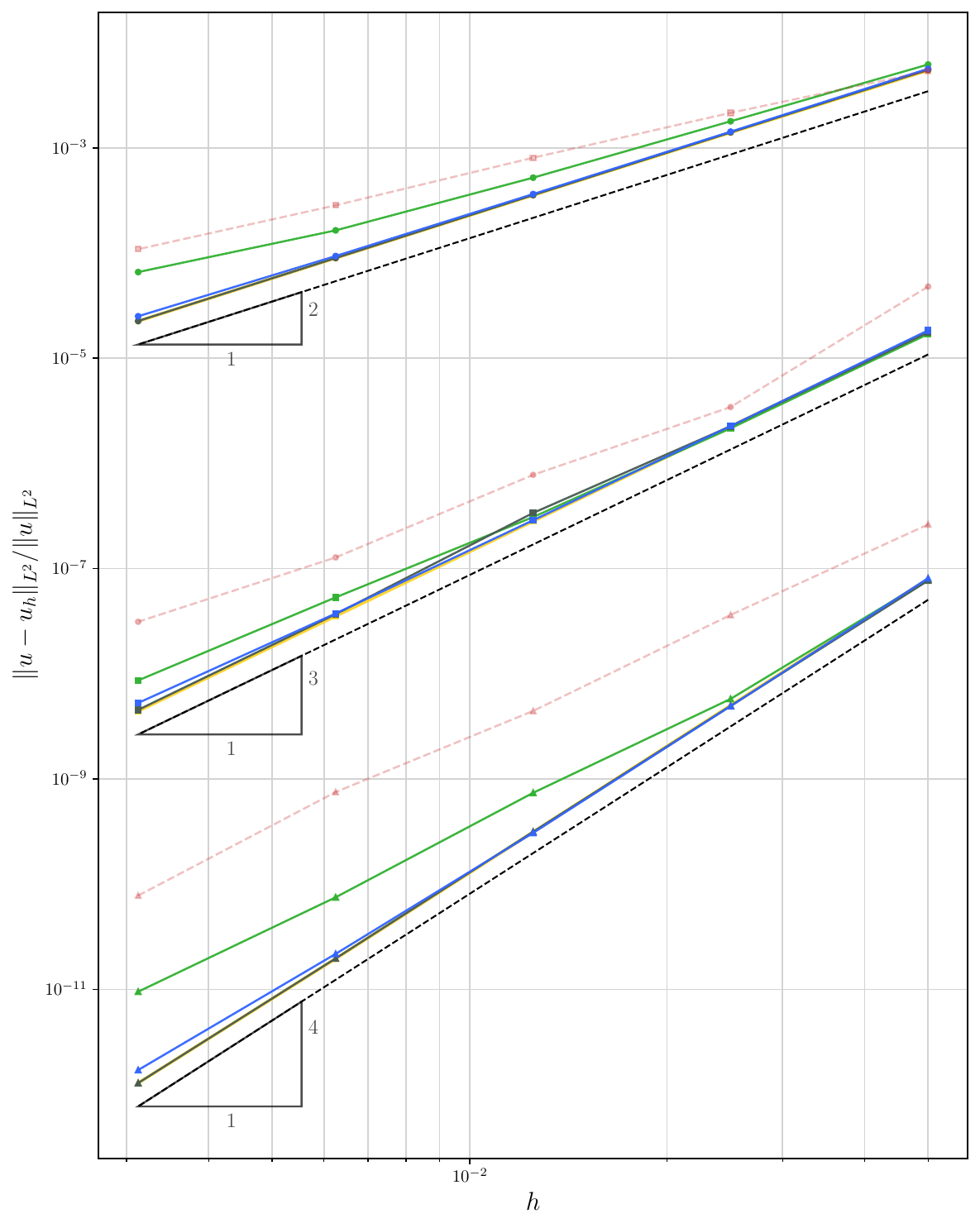}
        \caption{Average comparison $L^2$ norm}
        \label{fig:ex12_comp_a}
    \end{subfigure}
    \hfill
    \begin{subfigure}[t]{0.3\textwidth}
        \centering
        \includegraphics[width=\textwidth]{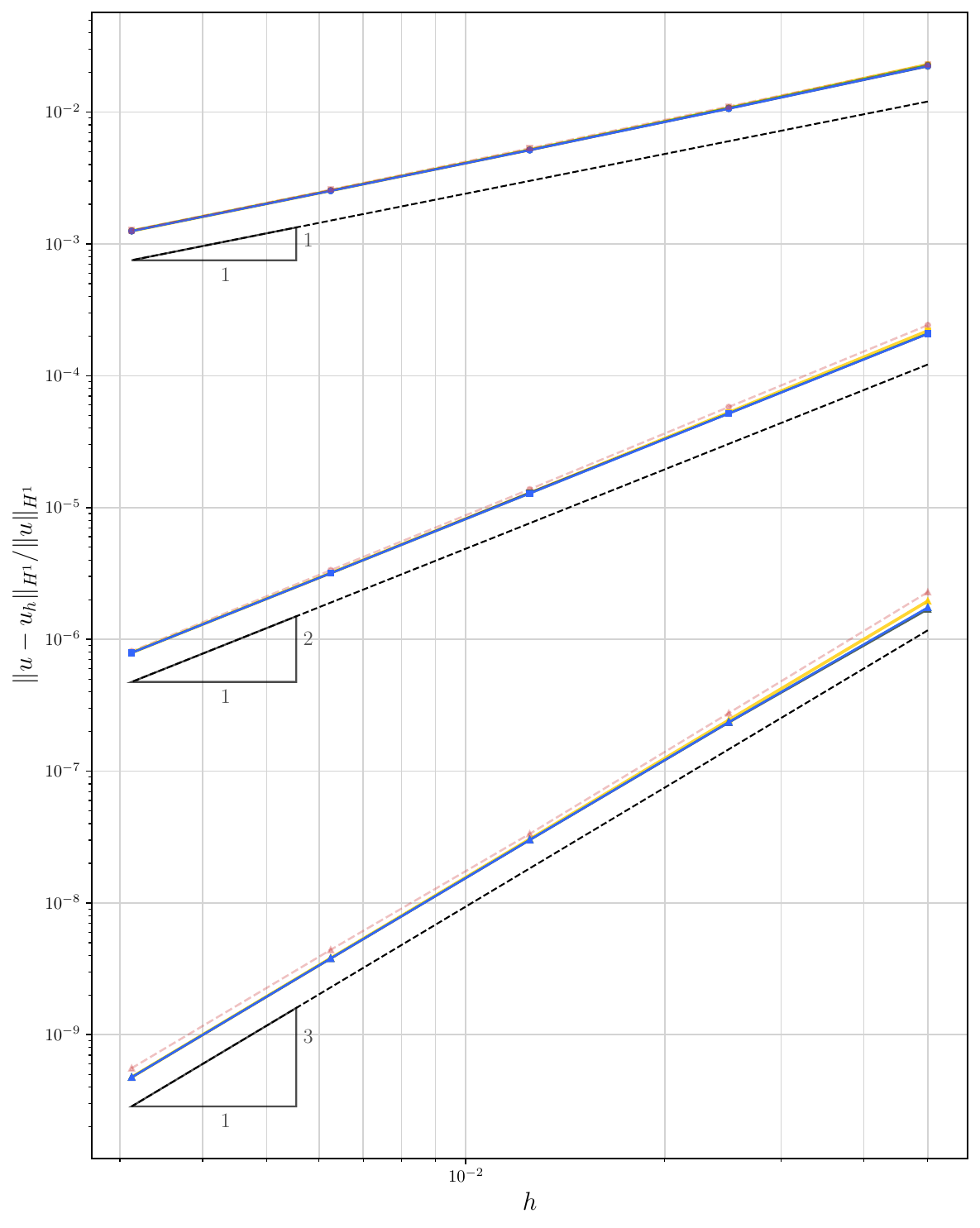}
        \caption{Average comparison $H^1$ norm}
        \label{fig:ex12_comp_b}
    \end{subfigure}
    \hfill
    \begin{subfigure}[t]{0.3\textwidth}
        \centering
        \includegraphics[width=\textwidth]{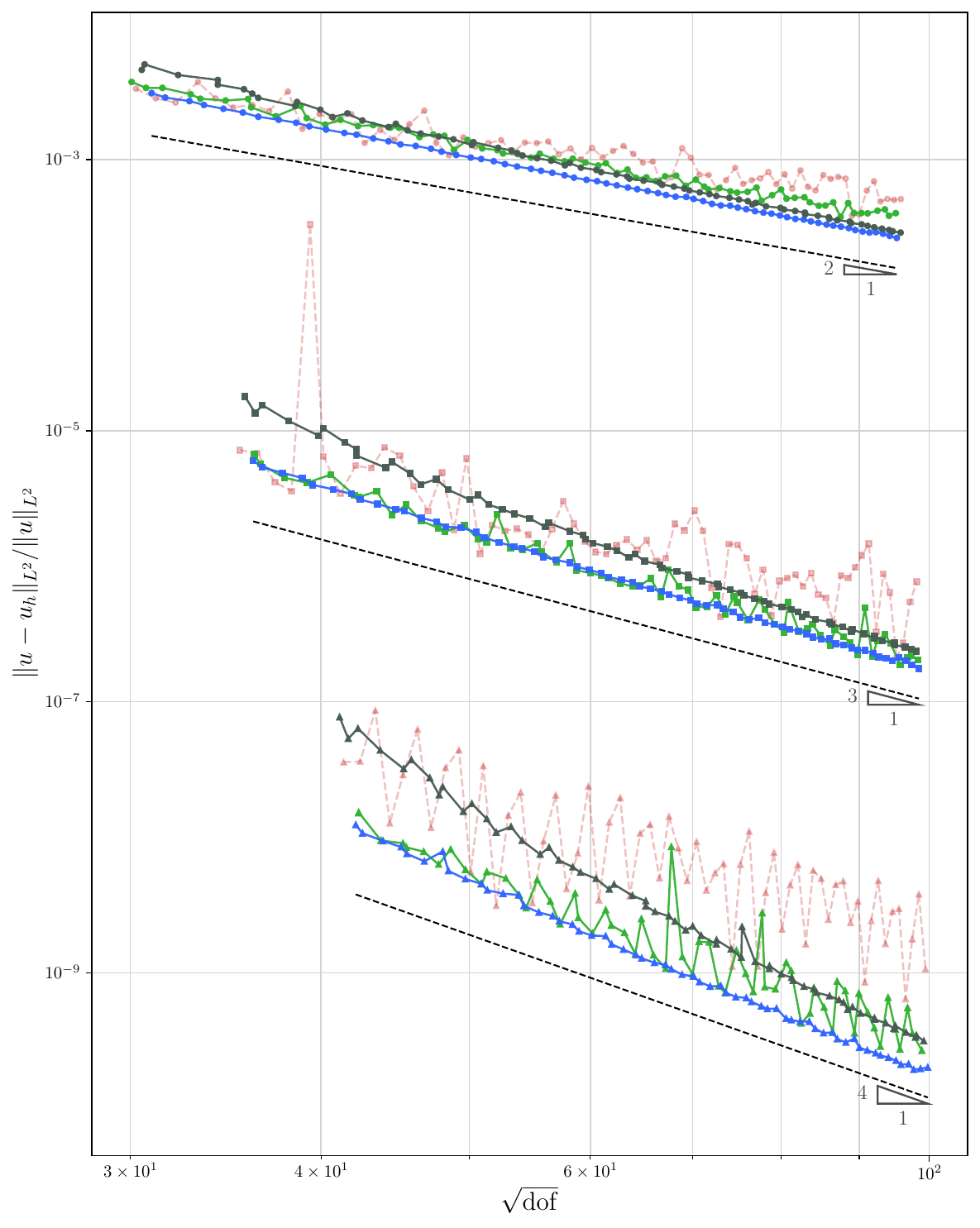}
        \caption{DOFs comparison}
        \label{fig:ex12_comp_c}
    \end{subfigure}

    \caption{Head-to-head convergence comparison, from $20\times20$ to $320\times320$ knot spans, for the $h$-, $p$-, and $k$-refinement schemes and the non-refined SBM, applied to the geometry presented in Fig.~\ref{fig:sbm-thb} with immersed Neumann BCs. The color and marker legends are the same as in Fig.~\ref{fig:ex12}.}
    \label{fig:ex12_comparison}
\end{figure}

Fig.~\ref{fig:ex12_comparison} presents a direct comparison between the standard SBM and the locally refined approaches for this scenario.
As expected, the $H^1$ convergence is already optimal for the standard SBM, so none of the refinements improves it.
In the $L^2$ norm, for a given mesh size, $p$- and $k$-refinements are the most accurate, with $k$-refinement slightly outperforming $p$-refinement.
Considering the number of degrees of freedom required to reach a target accuracy, $p$-refinement appears most efficient, followed by $h$-refinement.
$k$-refinement seems to be advantageous only on coarse meshes; for finer meshes it tends to require the same or more DOFs than the standard SBM to achieve comparable accuracy.
{\color{red} This behaviour should not be interpreted as a contradiction of standard IGA continuity arguments, since the comparison here does not isolate continuity at fixed degree, but rather reflects the interaction between local boundary refinement and interior discretization errors.
}

\subsection{Local Refinement with Mixed Dirichlet–Neumann Boundary Conditions}

\begin{figure}[htbp]
    \centering

    \begin{subfigure}[t]{0.3\textwidth}
        \centering
        \includegraphics[width=\textwidth]{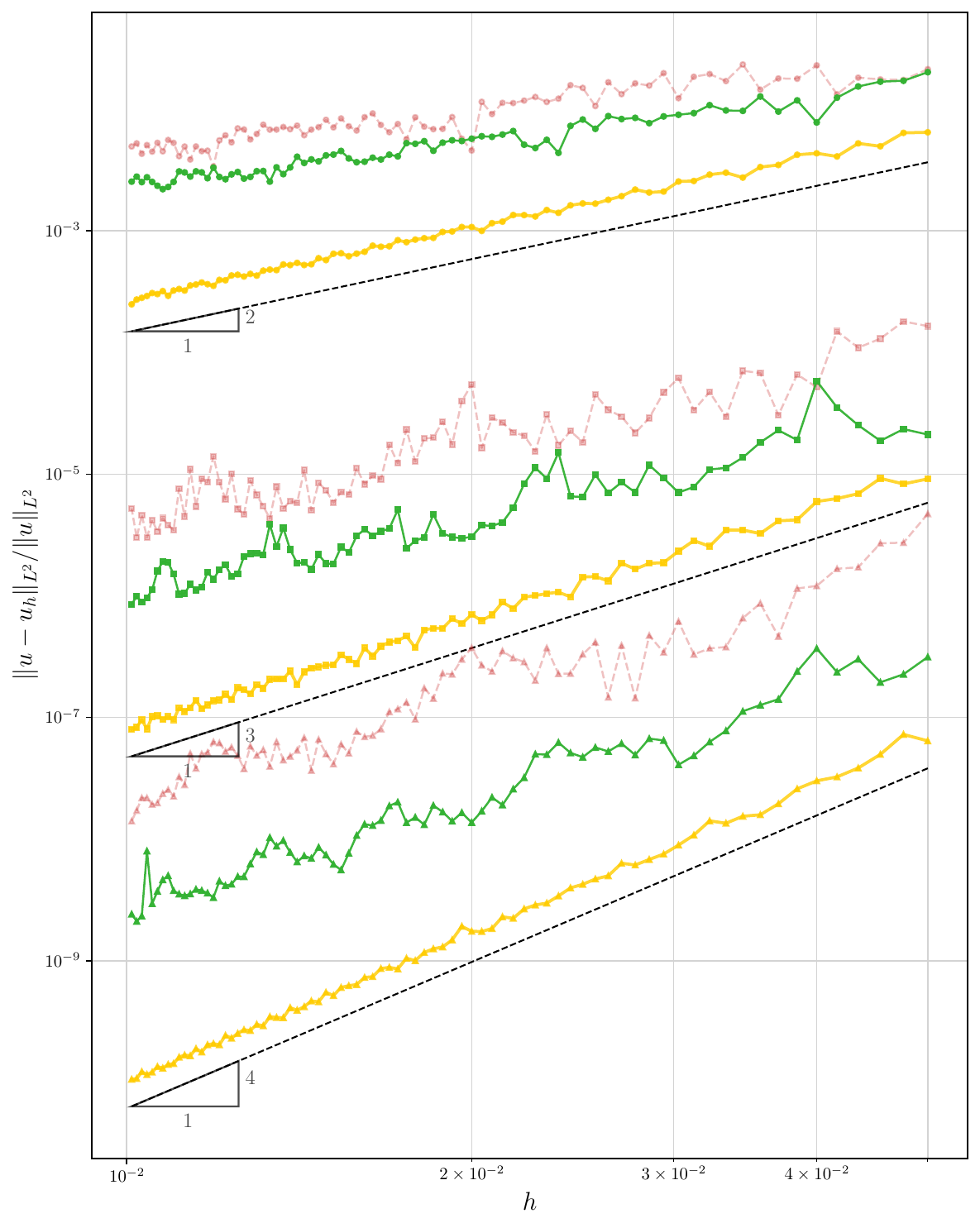}
        \caption{$h$-refinement}
        \label{fig:ex13_a}
    \end{subfigure}
    \hfill
    \begin{subfigure}[t]{0.3\textwidth}
        \centering
        \includegraphics[width=\textwidth]{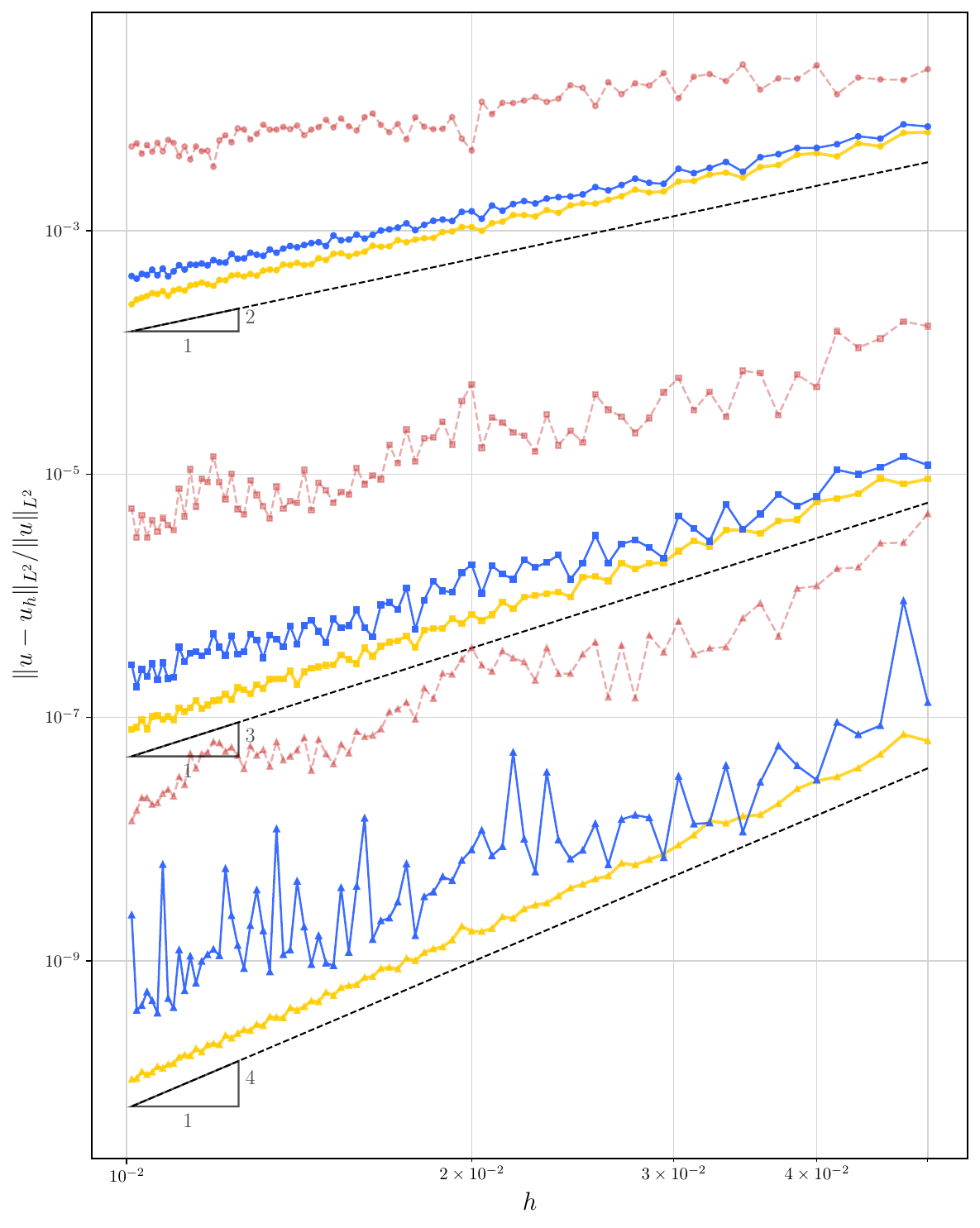}
        \caption{$p$-refinement}
        \label{fig:ex13_b}
    \end{subfigure}
    \hfill
    \begin{subfigure}[t]{0.3\textwidth}
        \centering
        \includegraphics[width=\textwidth]{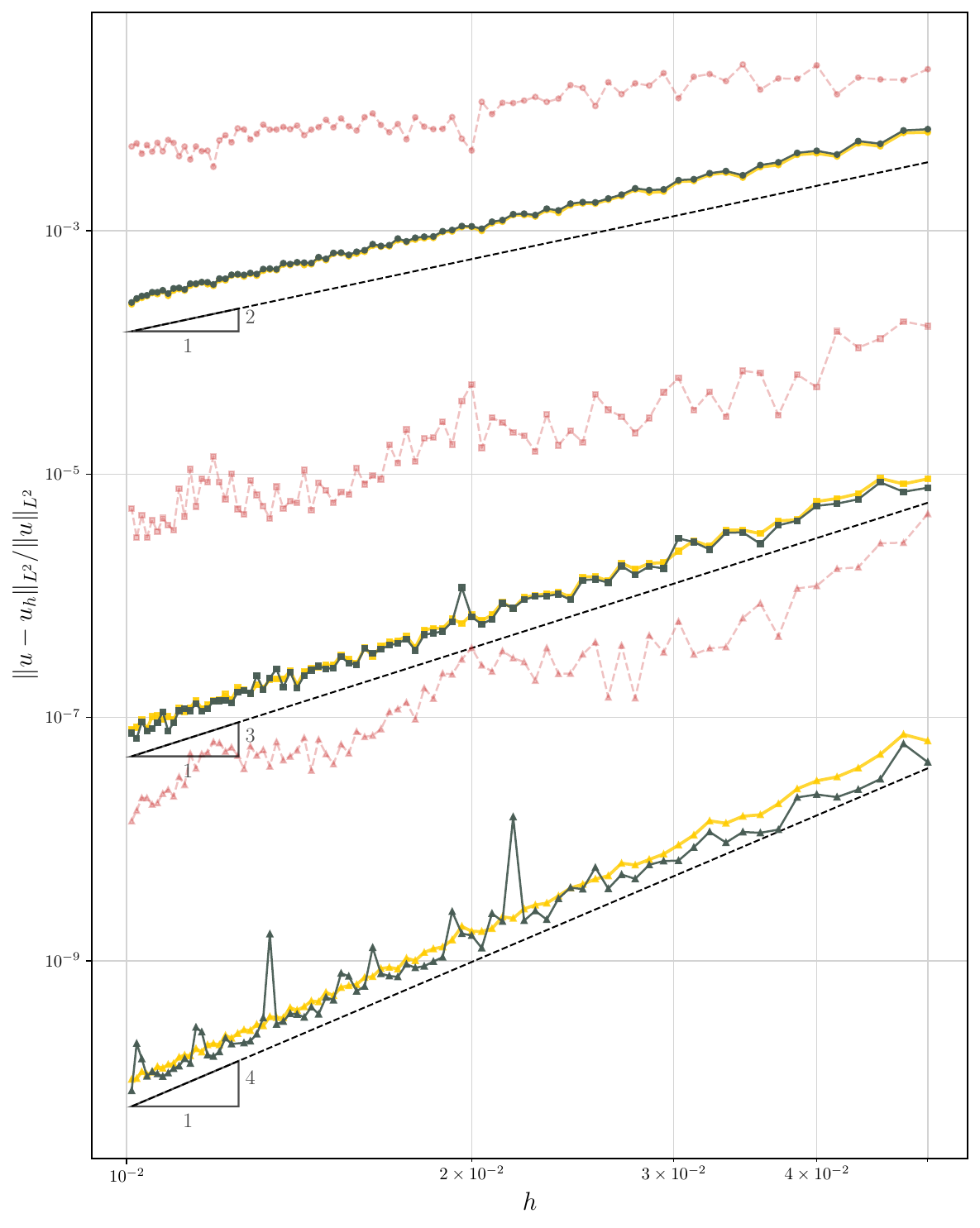}
        \caption{$k$-refinement}
        \label{fig:ex13_c}
    \end{subfigure}

    \vspace{1em}

    \begin{subfigure}[t]{0.3\textwidth}
        \centering
        \includegraphics[width=\textwidth]{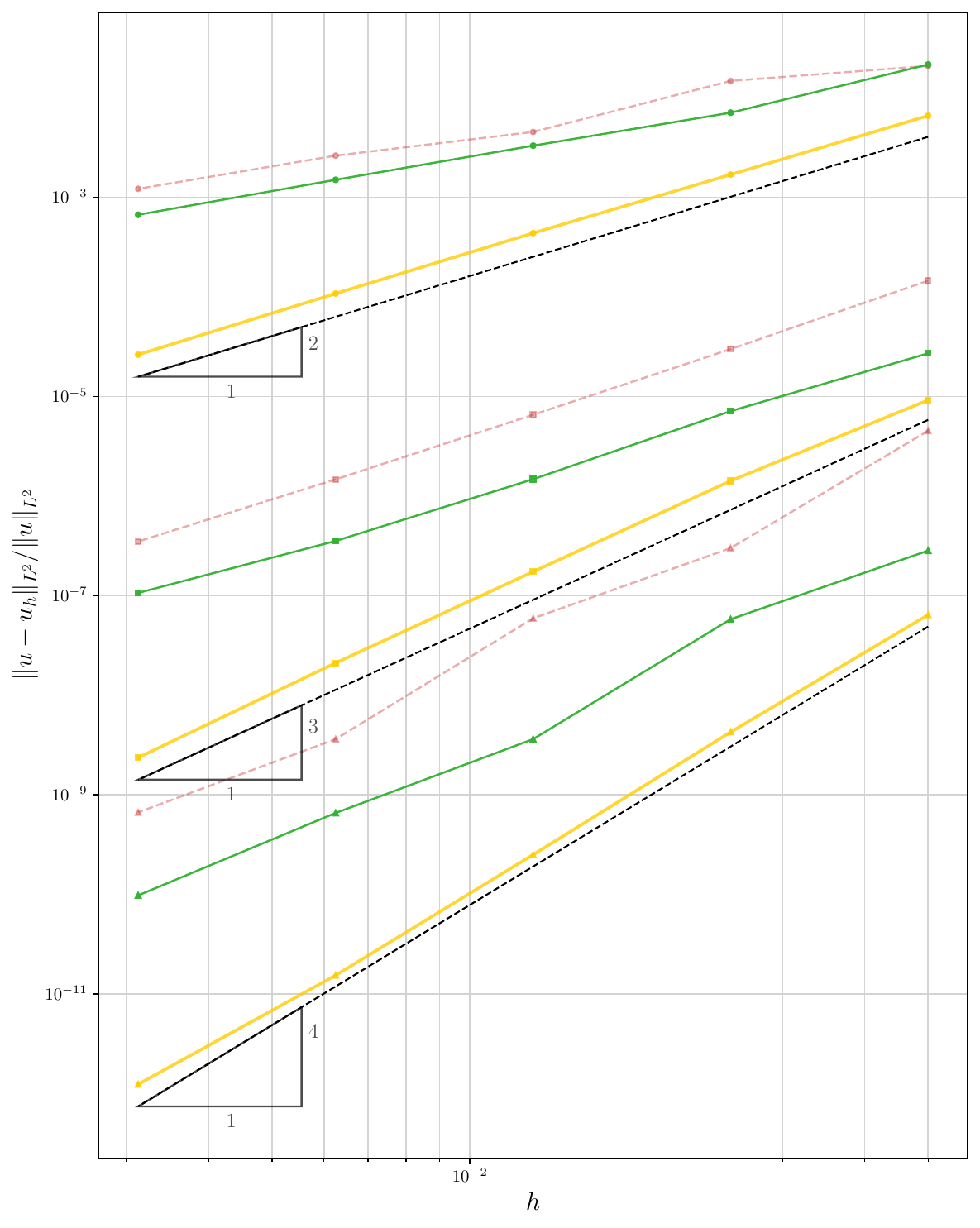}
        \caption{$h$-refinement}
        \label{fig:ex13_d}
    \end{subfigure}
    \hfill
    \begin{subfigure}[t]{0.3\textwidth}
        \centering
        \includegraphics[width=\textwidth]{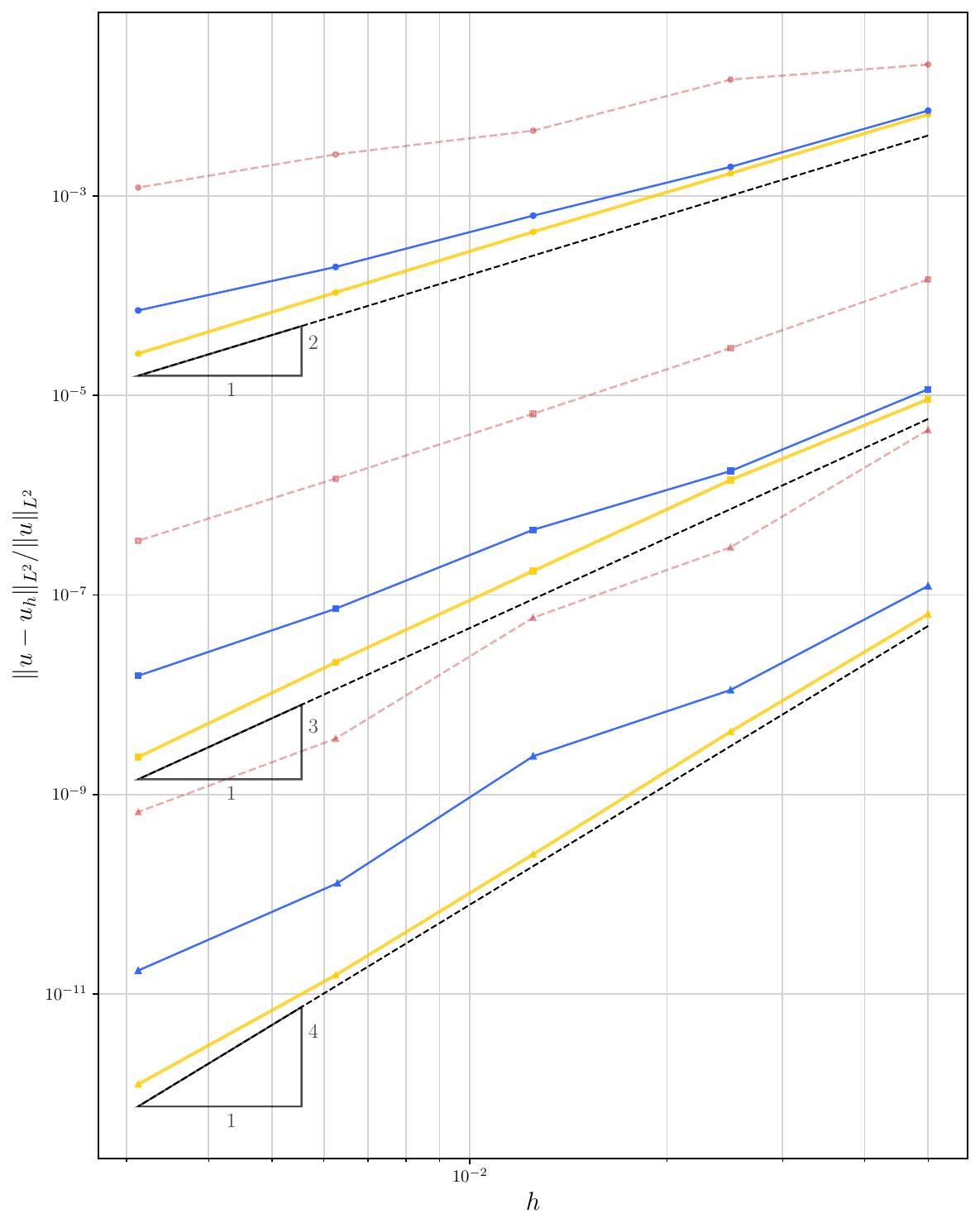}
        \caption{$p$-refinement}
        \label{fig:ex13_e}
    \end{subfigure}
    \hfill
    \begin{subfigure}[t]{0.3\textwidth}
        \centering
        \includegraphics[width=\textwidth]{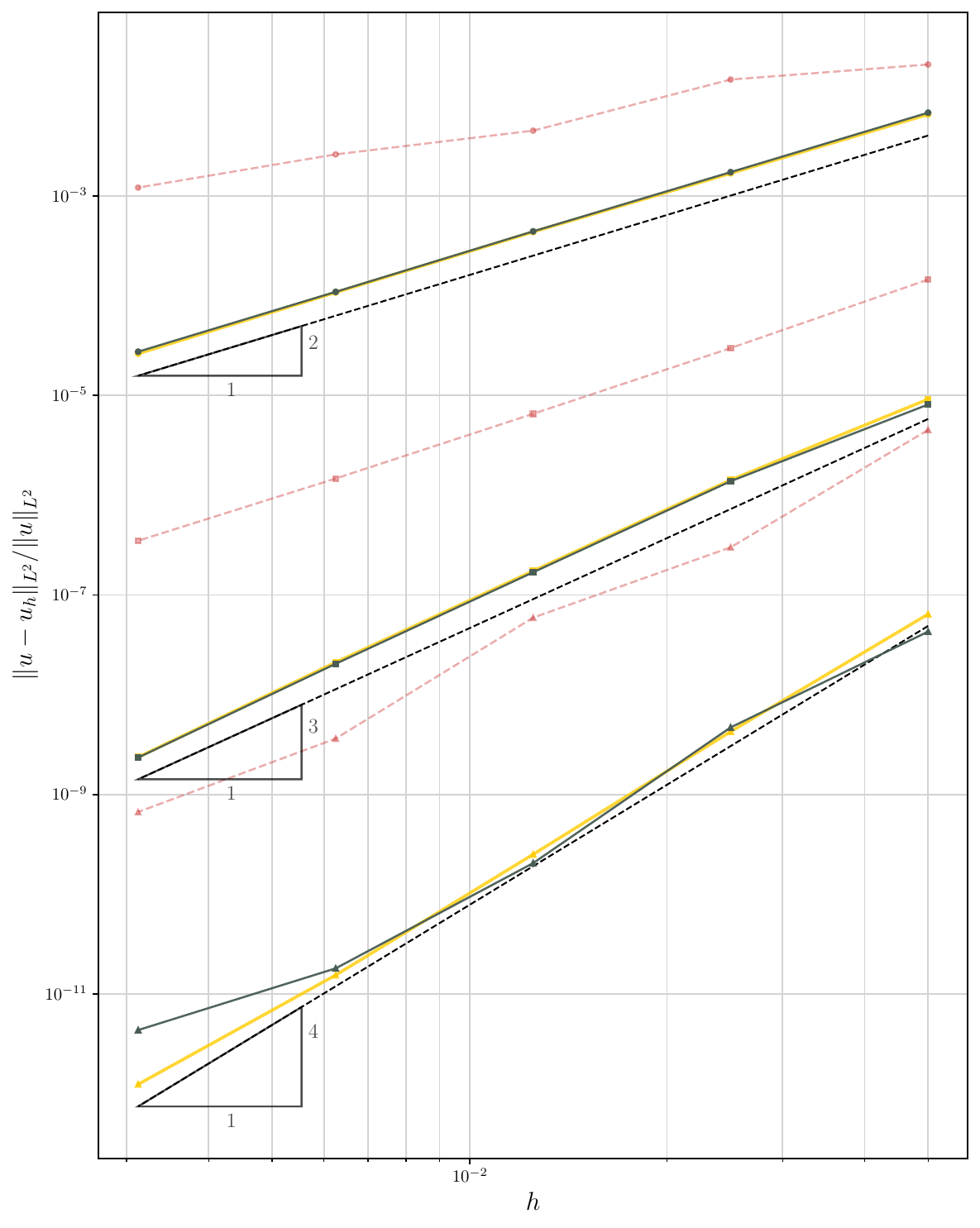}
        \caption{$k$-refinement}
        \label{fig:ex13_f}
    \end{subfigure}

    \caption{Convergence study for the geometry \ref{fig:sbm-overview} with immersed exterior Neumann boundary and immersed interior Dirichlet boundary, different local refinement strategies for the Neumann boundary. First row: step-by-step convergence from $20\times20$ to $100\times100$ knot spans; second row: standard halved convergence from $20\times20$ to $320\times320$ knot spans. Legend: yellow: body-fitted; dashed red: non-refined SBM; dashed black: reference line for the optimal convergence; green: $h$-refinement; blue: $p$-refinement; black: $k$-refinement. Circles: $p = 1$; squares: $p=2$; triangles $p=3$.}
    \label{fig:ex13}
\end{figure}

We now consider again the geometry of Fig.~\ref{fig:sbm-overview}: a circle of
radius \(r_{\mathrm{ext}} = 0.47\) centered at \((0.5,0.5)\) that
contains an interior circular hole of radius \(r = 0.10\), also centered
at \((0.5,0.5)\). Dirichlet boundary conditions are imposed on the inner hole, while Neumann conditions are prescribed on the outer boundary.
Fig.~\ref{fig:ex13} compares the convergence obtained with local
\(h\)-, \(p\)-, and \(k\)-refinements. All plots use the same color legend as in the previous examples.
{\color{red} All refinement schemes lead to a systematic reduction of the error. 
For \(h\)-refinement (Figs.~\ref{fig:ex13_a} and \ref{fig:ex13_d}), the convergence slope remains suboptimal, although the error constant decreases with increasing spline degree. This effect is more pronounced than in the previous example, since Neumann contributions dominate more strongly in the present setting. As discussed in Secs.~\ref{SBM} and \ref{THB for SBM}, local \(h\)-refinement near the boundary primarily reduces the projection distance \(\|\mathbf{d}\|\) without directly modifying the interior mesh size \(h\). A heuristic interpretation of the error behavior for Neumann boundary
conditions can be obtained from a
simplified model consistent with the Taylor-based shift expansion,
\[
E(h,\|\mathbf{d}\|)
\approx
C_{\mathrm{int}}\, h^{p+1}
+
C_{\mathrm{bnd}}\, \|\mathbf{d}\|^{p},
\]
where the first term represents the interior approximation error and
the second term reflects the boundary-consistency contribution associated
with the shift distance.
Recovering the optimal asymptotic slope requires a coupling between $h$ and 
$\|\mathbf{d}\|$ such that $\|\mathbf{d}\|^{p} = \mathcal{O}(h^{p+1})$. 
In the present study, only a single local refinement step at the boundary is 
employed by design in order to isolate the influence of the refinement strategy. 
Since the reduction of $\|\mathbf{d}\|$ is not explicitly coupled to the global 
mesh size, the boundary-consistency term becomes dominant for sufficiently small $h$. As a consequence, the convergence curve may eventually deviate from 
the optimal interior-dominated slope.}

The \(p\)-refinement in Figs.~\ref{fig:ex13_b} and \ref{fig:ex13_e} yields a larger accuracy gain, bringing the
convergence curve very close to the body-fitted reference, although at
fine meshes it still falls short of the optimal slope.  
By contrast, a single \(k\)-refinement in Figs.~\ref{fig:ex13_c} and \ref{fig:ex13_f}  maintains the optimal rate up to highly refined meshes, losing it just at the last refinement step of the halved convergence.

\begin{figure}[htbp]
    \centering

    \begin{subfigure}[t]{0.3\textwidth}
        \centering
        \includegraphics[width=\textwidth]{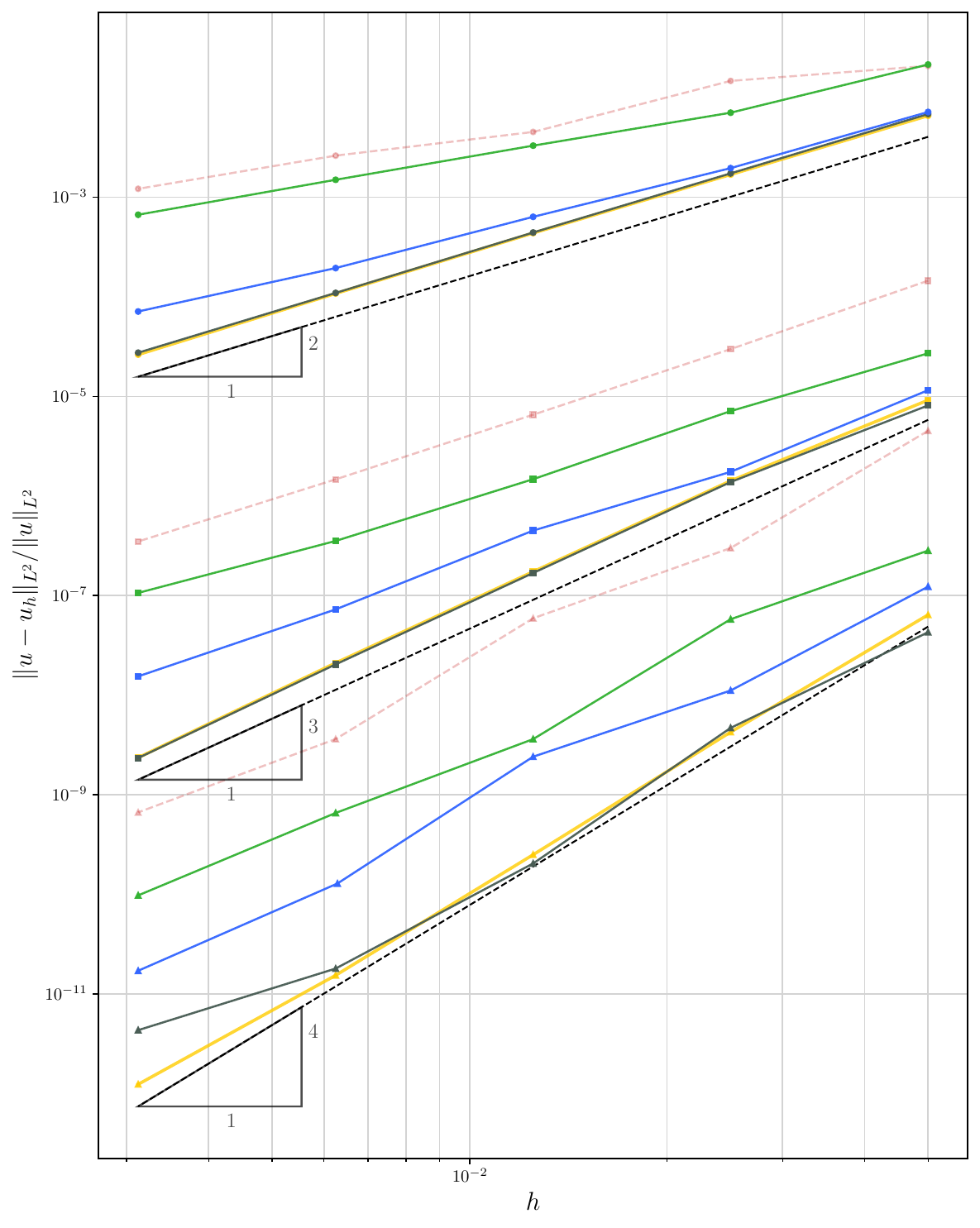}
        \caption{Average comparison $L^2$ norm}
        \label{fig:ex13_comp_a}
    \end{subfigure}
    \hfill
    \begin{subfigure}[t]{0.3\textwidth}
        \centering
        \includegraphics[width=\textwidth]{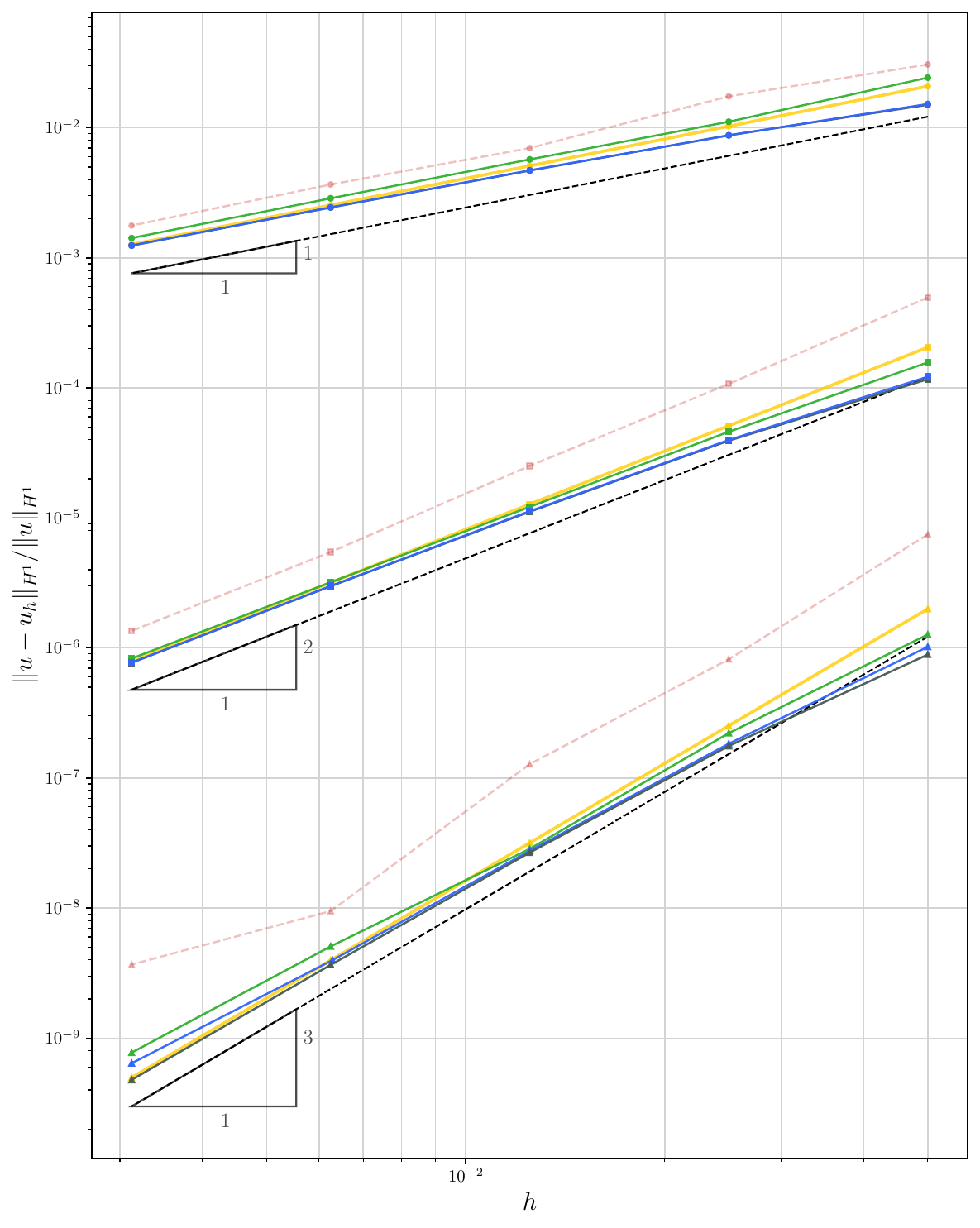}
        \caption{Average comparison $H^1$ norm}
        \label{fig:ex13_comp_b}
    \end{subfigure}
    \hfill
    \begin{subfigure}[t]{0.3\textwidth}
        \centering
        \includegraphics[width=\textwidth]{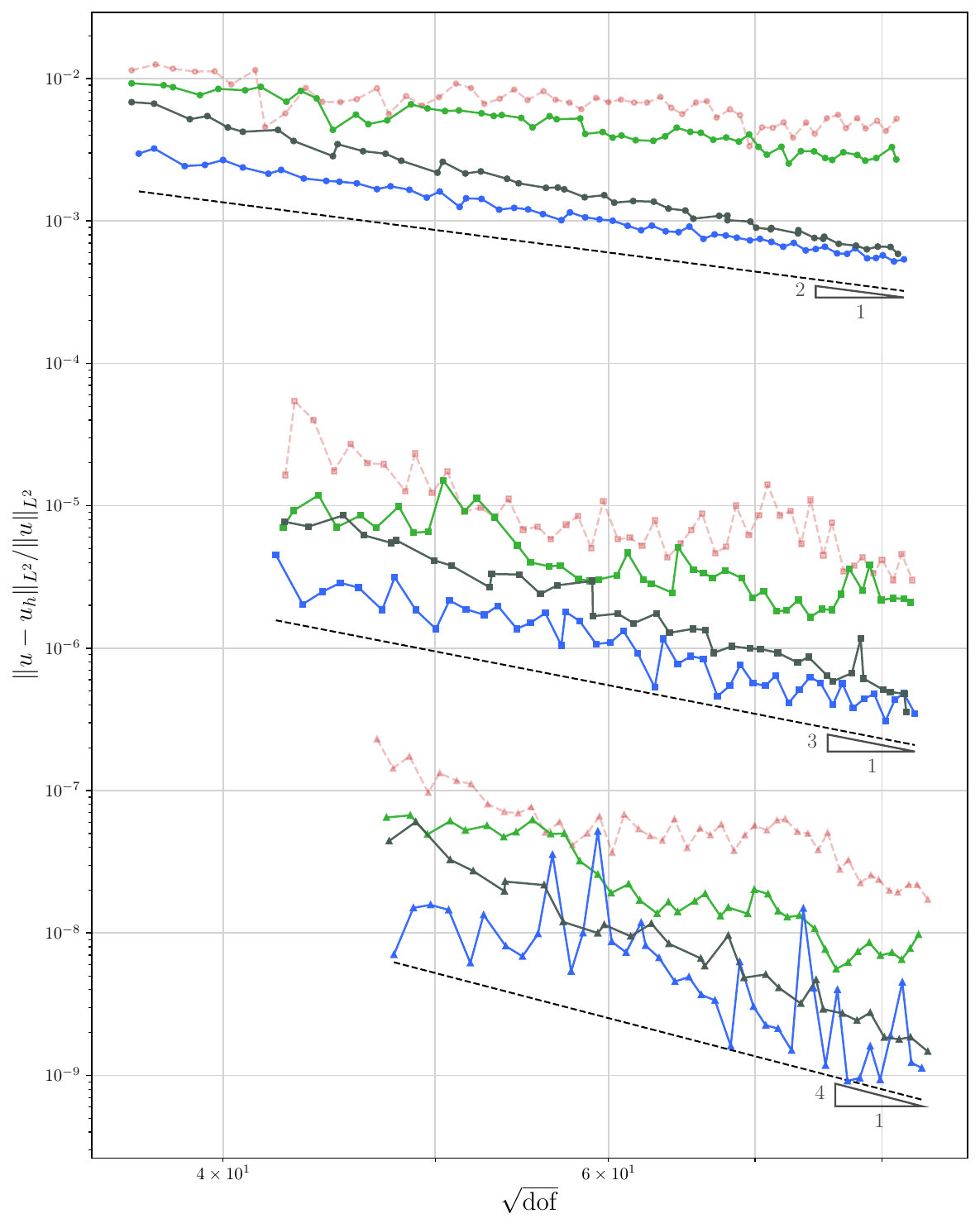}
        \caption{DOFs comparison}
        \label{fig:ex13_comp_c}
    \end{subfigure}

    \caption{Head-to-head convergence comparison, from $20\times20$ to $320\times320$ knot spans, for the $h$-, $p$-, and $k$-refinement schemes and the non-refined SBM, applied to the geometry presented in Fig.~\ref{fig:sbm-overview} with immersed exterior Neumann boundary and immersed interior Dirichlet boundary. The color and marker legends are the same as in Fig.~\ref{fig:ex13}.}
    \label{fig:ex13_comparison}
\end{figure}

Fig.~\ref{fig:ex13_comparison} summarises the three refinement
strategies.  The \(H^{1}\) error is already optimal without refinement,
and local refinement chiefly lowers the error constant, with little
difference among the three approaches. In terms of accuracy-to-cost ratio, \(p\)-refinement achieves a given error with the
fewest degrees of freedom, followed by \(k\)- and then \(h\)-refinement as can be seen in Fig.~\ref{fig:ex13_comp_c}.

We note here that in some simulations with local refinement the errors turn out to be smaller than in the unrefined body-fitted case. The reason is that the locally refined meshes contain more degrees of freedom than the body-fitted baseline, which allows a more accurate approximation of the solution.

\subsection{Applications to complex geometries with Neumann boundary conditions}

\begin{figure}[H]
  \centering
  \begin{subfigure}[t]{0.488\textwidth}
    \includegraphics[width=\linewidth]{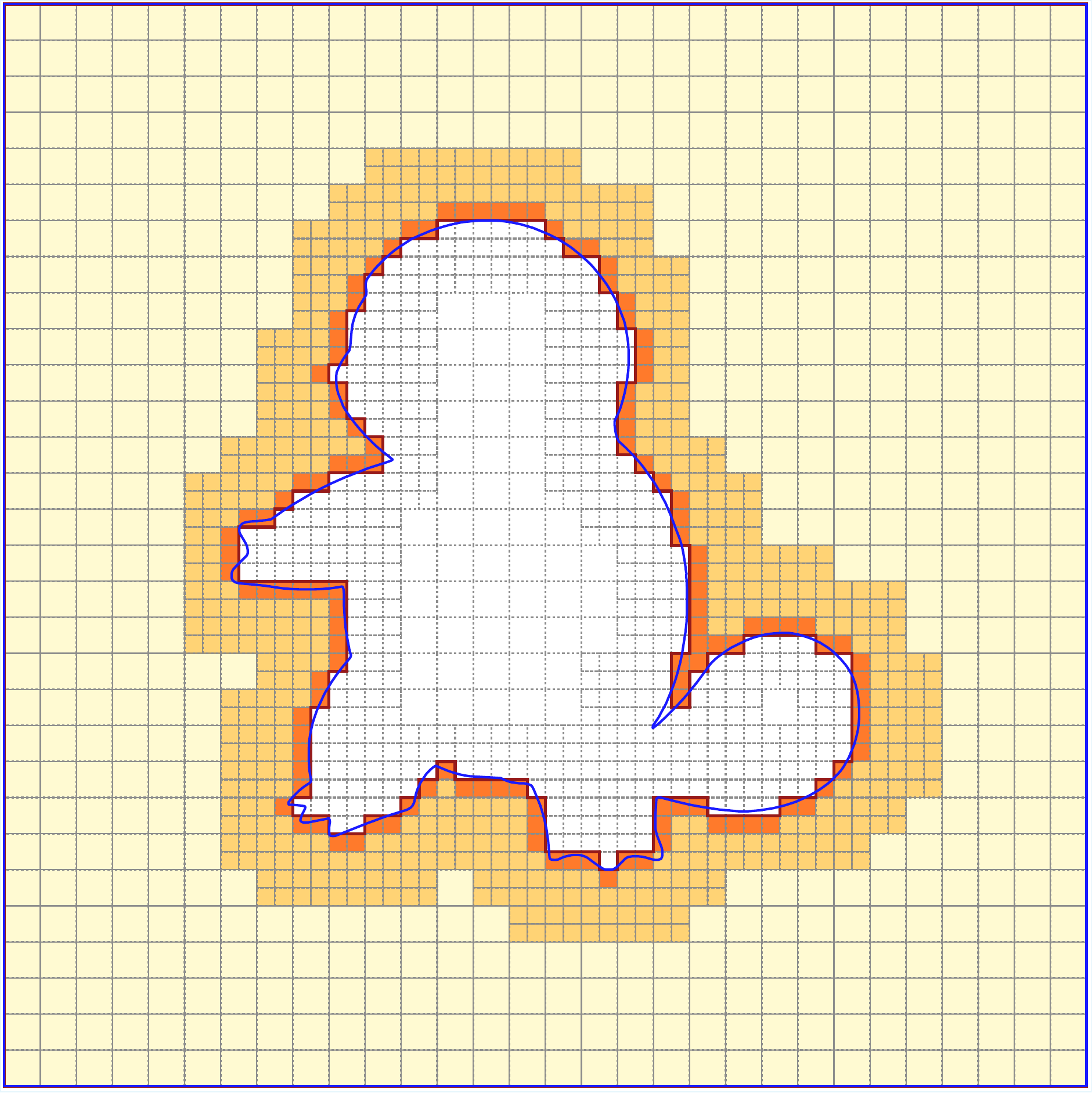}
    \caption{Squirtle}
    \label{fig:squirtle}
  \end{subfigure}
  \hspace{0.01\textwidth} 
  \begin{subfigure}[t]{0.49\textwidth}
    \includegraphics[width=\linewidth]{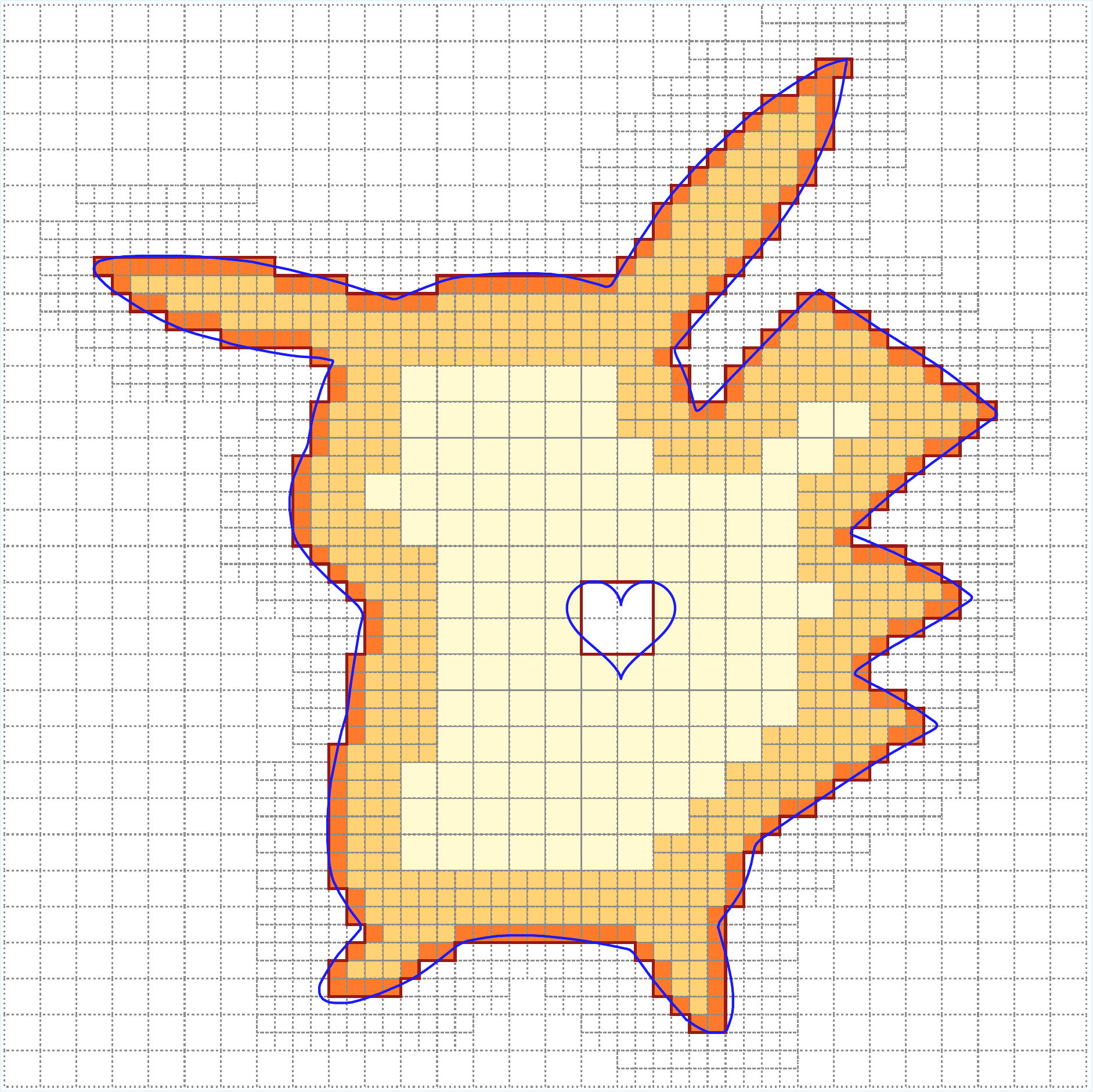}
    \caption{Pikachu}
    \label{fig:pikachu}
  \end{subfigure}
\caption{Example visualizations of the geometries with $30\times 30$ elements, polynomial degree $p=1$, and local refinement applied along the Neumann boundaries.}
  \label{fig:squirtle_pikachu}
\end{figure}

Now we turn to the more complex geometries of Fig.~\ref{fig:squirtle_pikachu}.
In Fig.~\ref{fig:squirtle}, a unit square with Dirichlet conditions on the outer boundary is considered, where the shape of a Squirtle is removed and a Neumann condition is prescribed along its boundary.
In Fig.~\ref{fig:pikachu}, the outer boundary carries Neumann conditions, while the shape of a heart is cut out and assigned Dirichlet conditions.

\begin{figure}[H]
  \centering
  \begin{subfigure}[t]{0.45\textwidth}
    \includegraphics[width=\linewidth]{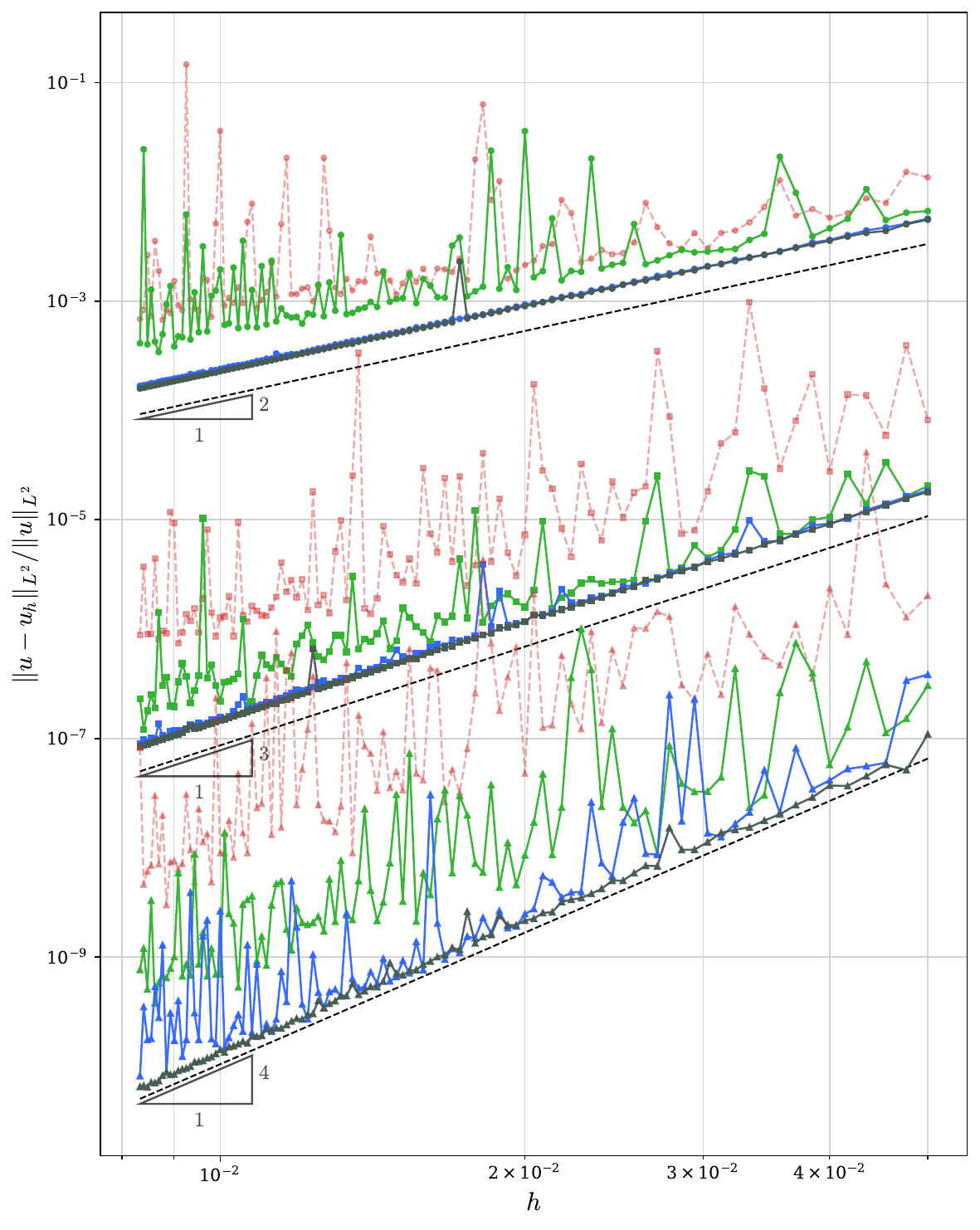}
    \caption{$L^2$ convergence study vs. characteristic mesh size.}
  \end{subfigure}
  \hspace{0.05\textwidth}
  \begin{subfigure}[t]{0.45\textwidth}
    \includegraphics[width=\linewidth]{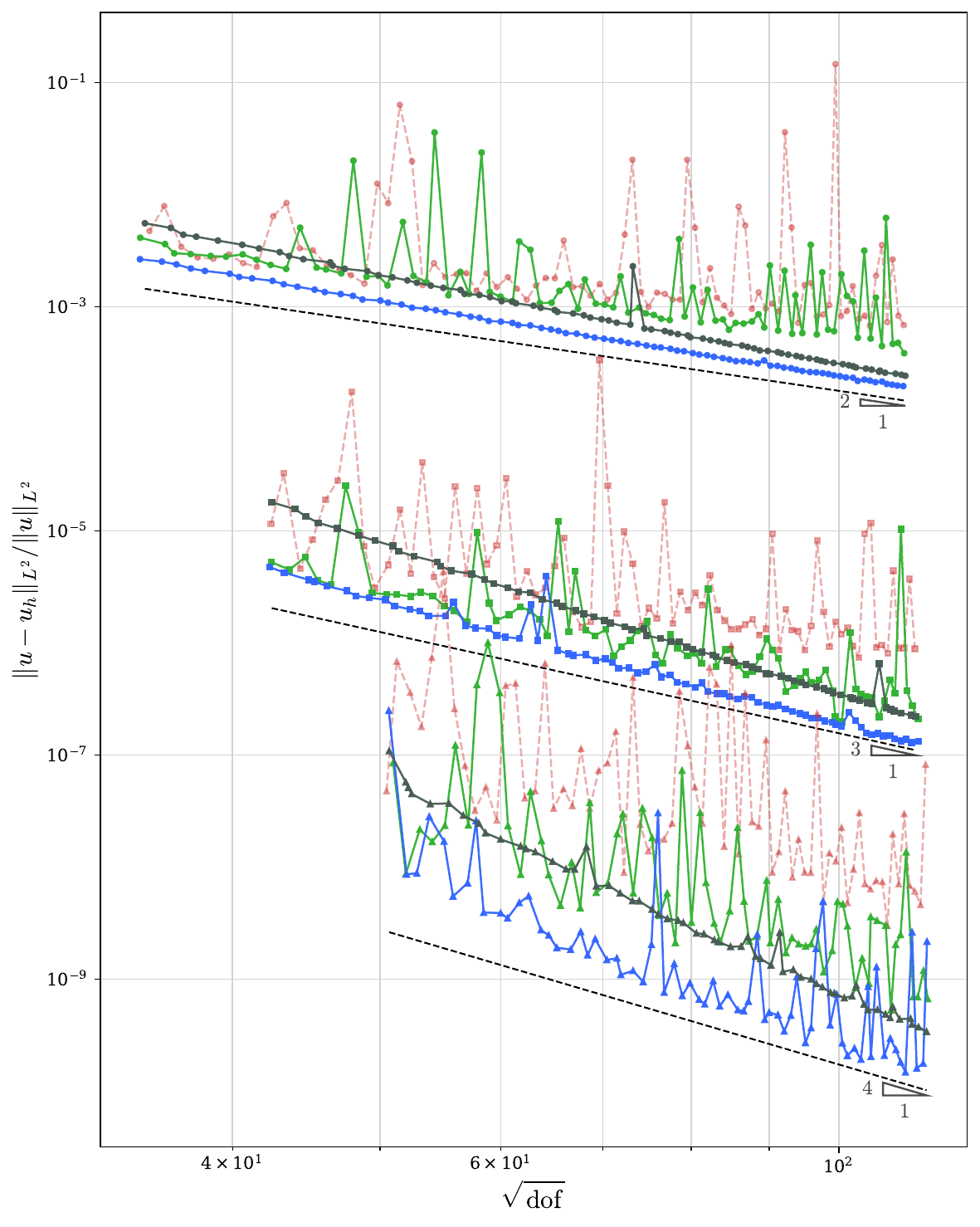}
    \caption{$L^2$ convergence study vs. DOFs.}
  \end{subfigure}
\caption{ Step-by-step convergence study for the geometry in Fig.~\ref{fig:squirtle}, using different local refinement strategies, internal Dirichlet boundary conditions, and external Neumann boundary conditions. The study spans from an initial mesh of $20\times20$ up to $100\times100$ knot spans.
     Legend: dashed red: non-refined SBM; dashed black: reference line for optimal convergence; green: $h$-refinement; blue: $p$-refinement; black: $k$-refinement. Circles: $p = 1$; squares: $p = 2$; triangles: $p = 3$.}

  \label{fig:squirtle convergence}
\end{figure}

\begin{figure}[H]
  \centering
  \begin{subfigure}[t]{0.45\textwidth}
    \includegraphics[width=\linewidth]{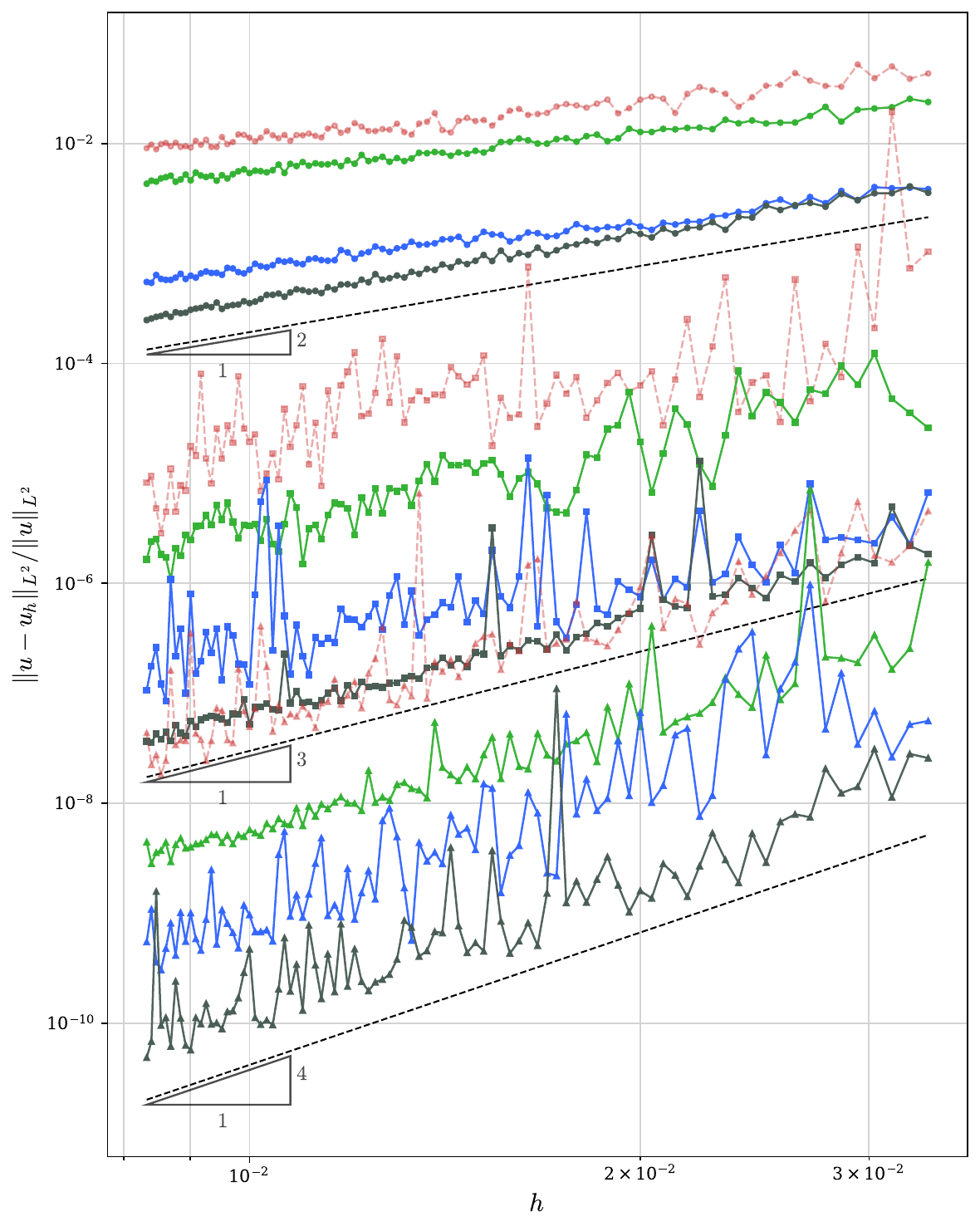}
    \caption{$L^2$ convergence study vs. characteristic mesh size.}
  \end{subfigure}
  \hspace{0.05\textwidth}
  \begin{subfigure}[t]{0.45\textwidth}
    \includegraphics[width=\linewidth]{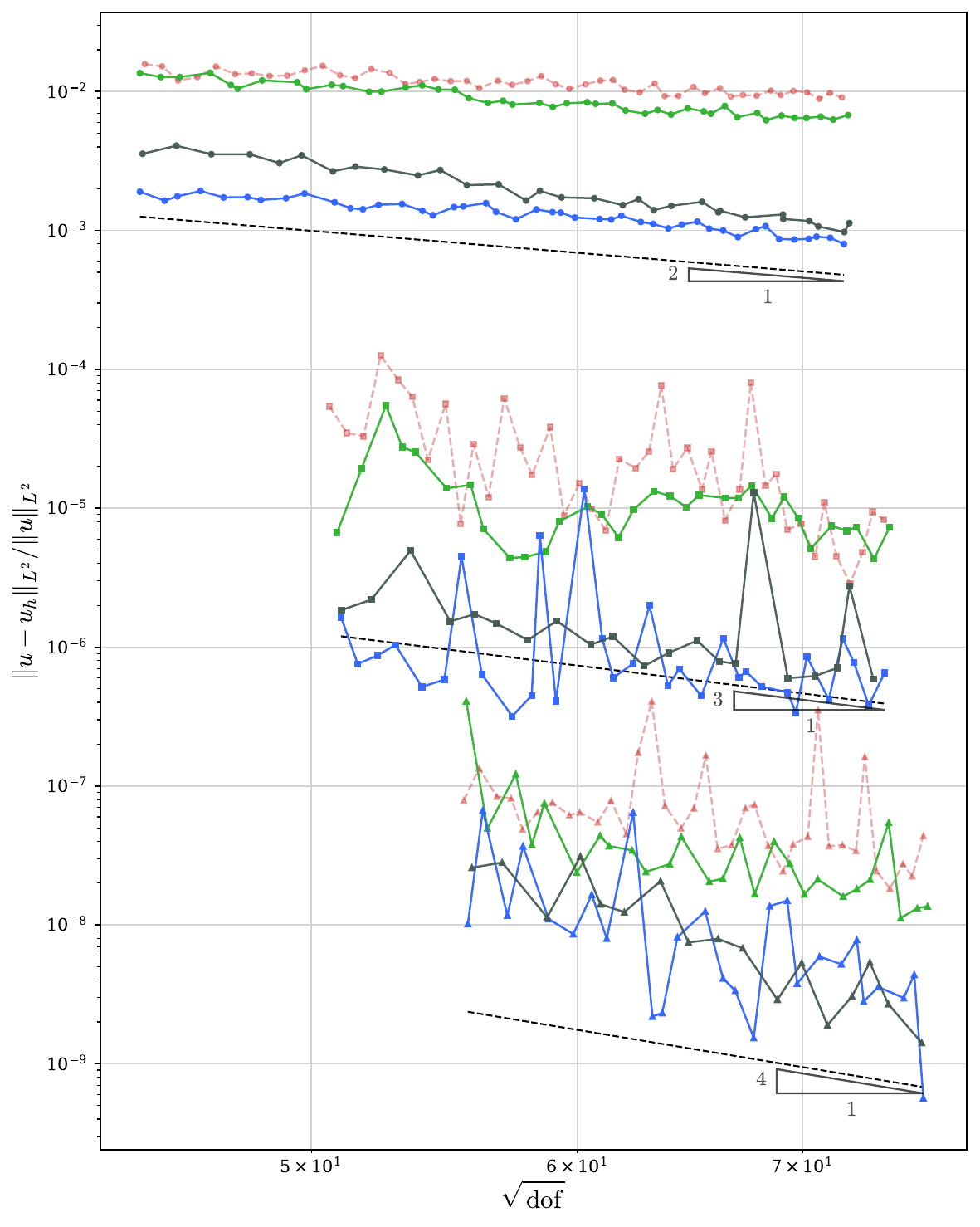}
    \caption{$L^2$ convergence study vs. DOFs.}
  \end{subfigure}
  \caption{
     Step-by-step convergence study for the geometry in Fig.~\ref{fig:pikachu}, using different local refinement strategies, internal Neumann boundary conditions, and external Dirichlet boundary conditions. The study spans from an initial mesh of $20\times20$ up to $100\times100$ knot spans.
     Legend: dashed red: non-refined SBM; dashed black: reference line for optimal convergence; green: $h$-refinement; blue: $p$-refinement; black: $k$-refinement. Circles: $p = 1$; squares: $p = 2$; triangles: $p = 3$.
  }
  \label{fig:pikachu convergence}
\end{figure}

Figs.~\ref{fig:squirtle convergence} and \ref{fig:pikachu convergence} compare the convergence under local $h$-, $p$-, and $k$-refinements. 
These last results highlight even more what we observed in the previous examples.
For the geometry in Fig.~\ref{fig:squirtle}, where the Neumann boundary is the interior one, both the $p$- and $k$-refinements maintain the optimal convergence until the end of the study, with the $k$-refinement showing a much more stable step-by-step accuracy improvement than the $p$-refinement. Again, in terms of DOFs cost, the $p$-refinement proves to be the most efficient choice.
For the geometry in Fig.~\ref{fig:pikachu} the biggest portion of the Neumann boundary allows a clearer differentiation of the local refinement techniques. Focusing on basis functions of order $p=1$ and $p=2$, the $k$-refinement is the only one to achieve optimal convergence, while $p$ and $h$ only improve the convergence constant - with $p$-refinement offering a significantly greater improvement than $h$-refinement. For $p=3$ basis functions, it is clear that none of the refinement techniques is able to recover the optimal convergence rate, but they still follow a consistent hierarchy from worst to best: standard SBM, $h$-refinement, $p$-refinement, and $k$-refinement.
What is remarkable is that for $p=3$ the convergence with standard SBM is almost coincident, but still more oscillatory than the $k$-refinement with $p=2$. Looking at the efficiency study against the DOFs, we observe that in this case, at least for $p=3$, the $k$-refinement does not scale much differently from the $p$-refinement. This means that the complexity of the Neumann boundary in this scenario can take full benefit of the higher cost of the $k$-refinement.

\newpage

\newpage

{\color{red}
\begin{rem}
Immersed discretization techniques may exhibit so-called cross-talk effects when geometric features such as narrow slits, holes, or cutouts intersect the support of spline basis functions. This phenomenon arises when the support of a basis function is fragmented into disjoint regions, which can introduce an artificial coupling between physically separated parts of the domain and thereby deteriorate the predictive quality of the simulation~\cite{coradello_hierarchically_2020, Lian2025CPD}. Owing to their comparatively large support, higher-order and higher-continuity spline bases can be particularly susceptible to this behavior.

A complete separation of the support is commonly referred to as \emph{Type~1 cross-talk}. In higher dimensions, even partially trimmed supports may generate non-physical fluxes across disconnected regions; this mechanism is typically classified as \emph{Type~2 cross-talk}.

Since the SBM is an immersed framework, cross-talk may in principle also arise in the present context. However, the manufactured solutions considered in this study are smooth and do not involve internal boundaries across which the flux vanishes. Accordingly, no indications of cross-talk were observed in the numerical experiments.

Local refinement near trimming features has been identified as an effective strategy for mitigating Type~1 cross-talk~\cite{coradello_hierarchically_2020}. While local $p$-refinement reduces the support of basis functions, the element size remains unchanged, implying that the support cannot become smaller than a single element. Consequently, $p$-refinement alone may be insufficient when trimming features are smaller than the element size, for instance when a slit partitions an element into disconnected regions. In such situations, refinement strategies that effectively decrease the element size, such as $h$- or $k$-refinement, become necessary. Alternatively, the control-point duplication approach proposed in~\cite{Lian2025CPD} offers another possibility for alleviating Type~1 cross-talk, whereas Type~2 cross-talk diminishes only in the asymptotic limit of successive $h$- (or $k$-) refinement.
\end{rem}
}


\section{Conclusions} \label{sec:conclusions}
The Shifted Boundary Method (SBM) provides an unfitted framework for
isogeometric analysis, but exhibits a reduction of the asymptotic
convergence rate in the $L^2$-norm under Neumann boundary conditions due to the reduced
admissible shift order for gradient quantities. In this work, the SBM is combined with locally refined THB-splines to investigate the effect of local refinement on the Neumann-induced boundary consistency error. In addition, an enhanced tensor-product shift operator is introduced and systematically compared to the standard operator. The study is structured in three steps:
\begin{enumerate}
    \item Implementation and Verification of local \( h \)-, \( p \)-, and \( k \)-refinement schemes for THB-splines within the SBM framework.
    \item Introduction of an enhanced shift operator and its evaluation for Dirichlet and Neumann boundary conditions.
    \item Numerical assessment of accuracy, stability, and computational cost on trimmed benchmark problems.
\end{enumerate}

\textbf{Key findings:}
\begin{itemize}

    \item The enhanced shift operator improves numerical robustness for Dirichlet boundary conditions and reduces oscillatory behavior.
    
    \item Local $h$-, $p$-, and $k$-refinement reduce the boundary-consistency error of the SBM and extend the regime in which convergence remains close to optimal for Neumann boundary conditions.

    \item Among the investigated strategies, $k$-refinement yields the
    highest overall accuracy and the most stable convergence behavior.
    In particular, it extends the regime of near-optimal convergence
    and reduces oscillations more effectively than pure $h$- or
    $p$-refinement, at the expense of increased computational cost.

    \item In terms of accuracy-to-cost ratio, $p$-refinement provides
    the most efficient strategy, achieving substantial error reduction
    with a moderate increase in degrees of freedom.

\end{itemize}

{\color{red}
The numerical investigations deliberately restrict all refinement strategies to a single local refinement step in order to isolate the effect of local refinement of the boundary on the convergence behavior under Neumann conditions and to enable a controlled comparison between the refinement types. For sufficiently fine meshes, the convergence rate eventually deviates from the optimal slope once the boundary consistency error becomes dominant. This indicates that additional local refinement near the surrogate boundary would be required to further extend the asymptotic regime.}

{\color{blue}
Although the numerical study is restricted to a scalar diffusion problem with smooth manufactured solutions, this setting allows us to clearly isolate the boundary-consistency effects inherent to the SBM. For problems with reduced regularity or strongly nonlinear behavior, the interior discretization error may dominate and refinement strategies should then target the principal error source within the domain. The proposed refinement framework is not limited to scalar diffusion problems. Since the SBM operates component-wise, the approach extends directly to vector-valued linear problems such as elasticity. For nonlinear applications of the SBM \cite{antonelli2025isogeometric}, where solution regularity may limit achievable convergence rates, local $h$-refinement is expected to provide a superior strategy than $p$-refinement.
}

\newpage

\section*{CRediT authorship contribution statement}
\noindent 
\textbf{Christoph Hollweck:} Software, Validation, Investigation, Conceptualization, Writing - Original Draft, Visualization.\\
\textbf{Andrea Gorgi:}  Software, Validation, Investigation, Conceptualization, Writing - Original Draft, Visualization.\\
\textbf{Nicol\`o Antonelli:} Software, Validation, Investigation.\\ 
\textbf{Marcus Wagner:} Review, Supervision.\\
\textbf{Roland Wüchner:} Review, Supervision, Funding acquisition.

\section*{Declaration of competing interest}
\noindent The authors declare that they have no known competing financial interests or personal relationships that could have appeared to influence the work reported in this paper.

\section*{Data availability}
\noindent Data will be made available on request.

\section*{Declaration of generative AI and AI-assisted technologies in the writing process}
\noindent During the preparation of this work, the authors used large language models (LLMs) solely to improve readability and language. The authors reviewed and edited the content as needed and take full responsibility for the content of the publication.

\section*{Acknowledgments}
\noindent The authors gratefully acknowledge the Design for IGA-type discretization workflows (GECKO) project.
The Design for IGA-type discretization workflows has received funding from the European Union's Horizon Europe research and Innovation programme under grant agreement No 101073106 Call: HORIZON-MSCA-2021-DN-01.

\noindent C. Hollweck acknowledges support from the PhD scholarship program of the Studienstiftung des deutschen Volkes and OTH Regensburg. The authors sincerely appreciate this support.
\newpage

\bibliographystyle{elsarticle-num} 
\bibliography{REFERENCES}

\end{document}